\documentclass[a4paper,11pt,reqno]{amsart}
\usepackage{amsmath}
\usepackage{amsthm}

\usepackage{amssymb}

\usepackage{color}

\usepackage{paralist}
  \usepackage{graphics} 
  \usepackage{epsfig} 
\usepackage{hyperref}

\usepackage[margin=2.35cm]{geometry}

\makeatletter
\def\cleardoublepage{\clearpage\if@twoside \ifodd\c@page\else%
		 \hbox{}%
	 \thispagestyle{empty}
	 \newpage%
	 \if@twocolumn\hbox{}\newpage\fi\fi\fi}
\makeatother

\let\cleardoublepage\clearpage

\newtheorem{theorem}{Theorem}[section]
\newtheorem{corollary}[theorem]{Corollary}
\newtheorem{lemma}[theorem]{Lemma}
\newtheorem{proposition}[theorem]{Proposition}
\newtheorem{definition}[theorem]{Definition}
\newtheorem{remark}[theorem]{Remark}
\numberwithin{equation}{section}

\begin{document}

\title[Porous media equations with two weights]{
Porous media equations with two weights: \\ smoothing and decay properties of energy solutions \\ via Poincar\'e inequalities}

\author {Gabriele Grillo, Matteo Muratori}

\address {Gabriele Grillo, Matteo Muratori: Dipartimento di Matematica, Politecnico di Milano, Piaz\-za Leonardo da Vinci 32, 20133 Milano, Italy}

\email {gabriele.grillo@polimi.it}

%

\email {matteo1.muratori@mail.polimi.it}

\author{Maria Michaela Porzio}

\address{Maria Michaela Porzio:
Dipartimento di Matematica, Universit\`a di Roma ``La Sapienza", Piaz\-zale A. Moro 2, 00185 Roma, Italy}

\email{porzio@mat.uniroma1.it}

\begin{abstract}
We study weighted porous media equations on domains $\Omega\subseteq{\mathbb R}^N$, either with Dirichlet or with Neumann  homogeneous boundary conditions when $\Omega\not={\mathbb R}^N$. Existence of weak solutions and uniqueness in a suitable class is studied in detail. Moreover,  $L^{q_0}$-$L^\varrho$ smoothing effects ($1\leq q_0<\varrho<\infty$) are discussed for short time, in connection with the validity of a Poincar\'e inequality in appropriate weighted Sobolev spaces, and the long-time asymptotic behaviour is also studied. In fact, we prove full equivalence between certain $L^{q_0}$-$L^\varrho$ smoothing effects and suitable weighted Poincar\'e-type inequalities. Particular emphasis is given to the Neumann problem, which is much less studied in the literature, as well as to the case $\Omega={\mathbb R}^N$ when the corresponding weight makes its measure finite, so that solutions converge to their weighted mean value instead than to zero. Examples are given in terms of wide classes of weights.
\end{abstract}

\maketitle
\footnotesize
\tableofcontents
\normalsize

\begin{section}{Introduction}
The goal of this paper is the study of both the Dirichlet and the Neumann homogeneous problem for weighted porous media equations (WPME for short) on Euclidean domains. In particular, given a domain $\Omega\subset{\mathbb R}^N$ and weights $\rho_\nu,\rho_\mu>0$  independent of time and  satisfying assumptions which will be stated below, we shall deal with the Dirichlet problem
\begin{equation} \label{eq: pmeDir}
\begin{cases}
\rho_{\nu}\,u_t =  \operatorname{div}\left(\rho_{\mu} \, \nabla{\left( u^m \right)} \right) & \textnormal{in} \  \Omega\times(0,\infty)  \\
 u =0  & \textnormal{on} \ \partial\Omega \times (0,\infty)  \\
 u(\cdot,0)=u_0(\cdot) & \textnormal{in} \ \Omega
 \end{cases}
\end{equation}
and with the Neumann problem
\begin{equation} \label{eq: pmeNeu}
\begin{cases}
\rho_{\nu}\,u_t = \operatorname{div}\left(\rho_{\mu} \, \nabla{\left( u^m \right)} \right) & \textnormal{in} \  \Omega\times(0,\infty)  \\
 \rho_{\mu} \, \frac{\partial{(u^m)}}{\partial{\mathbf{n}}} =0 & \textnormal{on} \ \partial\Omega \times (0,\infty) \\
 u(\cdot,0)=u_0(\cdot) & \textnormal{in} \ \Omega
 \end{cases} \,.
\end{equation}
For any $q>0$ and $x\in\mathbb R$ we use the convention $x^q=|x|^{q} \operatorname{sign}(x)$, and we shall assume throughout the paper that $m>1$. Precise meaning to the concepts of solution will be given in Section \ref{sec: wSS}. We mention from the beginning that, in fact, also the case $\Omega=\mathbb R^N$ can and shall be dealt with.

Problems \eqref{eq: pmeDir}, \eqref{eq: pmeNeu} are  generalizations of the well-known porous media equation, see e.g.\ the recent monograph \cite{Vaz07}. The appearance of the weights $\rho_\nu, \rho_\mu$ corresponds to spatial nonhomogeneity of the medium, either as concerns mass density and as concerns the diffusion coefficient. The one-weight case on ${\mathbb R}^N$ (i.e. $\rho_{\mu} \equiv 1$), mainly for a (mass) weight decaying as a negative power of $|\mathbf{x}|$ at infinity, is thoroughly analysed in the papers \cite{KR1, KR2, Eid90, EK, RV06, RV08, RV09, KRV10}.  Dealing with weights satisfying such asymptotic properties allows the Authors of some of those papers to perform an explicit and very detailed analysis of fine asymptotic properties of solutions in terms of suitably defined Barenblatt--type fundamental solutions (or, in some cases, of solutions with separated variables), much in the spirit of the unweighted case.

We shall first prove existence of solutions for the above mentioned problems and data in suitable $L^{q_0}(\Omega;\nu)$ spaces, both for the Dirichlet and the Neumann case, provided the weights are supposed to be strictly positive and sufficiently regular in $\Omega$ (see Section \ref{sec: wellP} for the details). We proceed by a careful generalization of the strategy of \cite{Vaz07}, which turns out to be applicable both to homogeneous Dirichlet and Neumann boundary conditions (for sufficiently smooth data and $\rho_\nu=\rho_\mu$ see also \cite{DG+08}). Key energy inequalities are shown and the fundamental $L^1$-contraction and comparison inequality \eqref{eq: princConf} is proved. When the data have lower integrability, we construct the corresponding solutions as \it limit solutions\rm, according to the terminology of \cite{Vaz07}. Uniqueness of solutions to the differential equations considered in general does not hold in $L^\infty$ even in the one-weight case, as already noticed in \cite{Eid90}. However it does indeed if it is sought for in a suitable class, which we identify as the class of \it weak energy solutions\rm, see Definitions \ref{wes-dirichlet} and \ref{wes-neumann}.

Two-weight operators are of common use in linear analysis (see e.g. \cite{Dav89}) and are widely studied in several different context: for example, every Riemannian Laplacian can be written, locally, as a two-weight linear operator  of second order. Still in the linear case, the validity of functional inequalities for the quadratic form associated to the generator of an evolution is well-known to be strictly related to regularizing and asymptotic properties of the evolution itself. In the nonlinear case, such connection is more subtle and less investigated (see e.g. \cite{BG05m, BG05p, BGV08, P09, BGV} and references quoted therein). Our main goal is contributing further to such an analysis here.

In fact, our aim will not be a full investigation of all the qualitative properties of solutions for some explicit class of weights but, rather, the comprehension of how the validity of $L^{q_0}$-$L^{\varrho}$ regularizing and asymptotic properties for solutions to \eqref{eq: pmeDir}, \eqref{eq: pmeNeu} is related to functional inequalities, naturally associated to the weights considered, for the largest possible choice of weights compatible with this approach. More precisely, our $L^{q_0}$-$L^{\varrho}$ estimates will follow by making use only of \emph{Poincar\'e inequalities} of the form
\begin{equation}\label{eq: poinIntro}
\left\| v\right\|_{2;\nu} \leq C_P \left\| \nabla{v} \right\|_{2;\mu} \ \ \ \forall v \in W^{1,2}_0(\Omega;\nu,\mu)
\end{equation}
or of the form
\begin{equation}\label{eq: poinIntroMed}
\left\| v-\overline{v} \right\|_{2;\nu} \leq M_P \left\| \nabla{v} \right\|_{2;\mu} \ \ \ \forall v \in W^{1,2}(\Omega;\nu,\mu)
\end{equation}
for suitable positive constants $C_P$ and $M_P$, where $W^{1,2}_0(\Omega;\nu,\mu)$ and $W^{1,2}(\Omega;\nu,\mu)$ are weighted Sobolev spaces, whose definition we shall give below, the indices $\nu,\mu$ mean that the corresponding norms are taken w.r.t.\ the measures ${\rm d}\nu=\rho_\nu\,{\rm d}\mathbf{x}$, ${\rm d}\mu=\rho_\mu\,{\rm d}\mathbf{x}$ and $\overline{v}$ denotes the mean value of $v$ w.r.t.\ $\nu$,  provided  (in this last case) $\rho_\nu$ is integrable.  We shall refer to \eqref{eq: poinIntroMed} as the \emph{zero-mean} Poincar\'e inequality. Actually, when dealing with the Neumann problem \eqref{eq: pmeNeu}, some of our results will hold starting from the validity of the inequality
\begin{equation}\label{eq: weakPoi}
\left\| v\right\|_{2;\nu} \leq W_P\left( \left\| \nabla{v} \right\|_{2;\mu}+\left\| v\right\|_{1;\nu} \right) \ \ \ \forall v \in W^{1,2}(\Omega;\nu,\mu),
\end{equation}
which is clearly weaker than \eqref{eq: poinIntroMed}, and which will turn out to be \it equivalent \rm to certain smoothing effects we shall prove.

It should be noted that inequalities like \eqref{eq: poinIntro} and \eqref{eq: poinIntroMed} have a clear spectral interpretation in terms of the differential operator formally given by \[\mathcal{L}_{\nu,\mu}=-\rho_\nu^{-1}\operatorname{div}(\rho_\mu\nabla) \, , \] which is formally self-adjoint and nonnegative in $L^2(\Omega;\nu)$: in fact, \eqref{eq: poinIntro} amounts to requiring that $\min\mathcal{S}(\mathcal{L}_{\nu,\mu})\ge 1/{C_P^2}>0$, where $\mathcal{S}(\mathcal{L}_{\nu,\mu})$ denotes the $L^2(\Omega; \nu)$ spectrum of $\mathcal{L}_{\nu,\mu}$, whereas \eqref{eq: poinIntroMed} says that $\min\left[\mathcal{S}(\mathcal{L}_{\nu,\mu})\setminus\{0\}\right]\ge 1/{M_P^2}>0$. It seems therefore quite unlikely, from the linear situation, that such inequalities can be related to anything but a long-time bound on the $L^2$ norm of the solutions. However in \cite{Gri10} it has been shown, in the context of $p$-Laplacian-type operators, that suitable Poincar\'e-type inequalities imply that $L^2$-$L^\varrho$ quantitative regularizing effects ($2<\varrho<\infty$) hold true, they being in turn equivalent to the Poincar\'e inequality one starts from. One of our goals is to investigate a similar connection in the porous media case, but having in mind the significantly greater difficulties which have to be expected for the WPME, especially in the Neumann case.

In fact, in our context, the standard smoothing effect (take for example the simpler Dirichlet case)
\[\|u(t)\|_\infty\le C\,t^{-\frac{\alpha}{m-1}}\,\|u_0\|_{q_0;\nu}^{1-\alpha} \] need not hold for any positive $\alpha$, in the strong sense that solutions can be \it unbounded at all times \rm (the same fact holding in the Neumann case, see the counterexamples in Sections \ref{sec: regAs} and \ref{sec: pmeNeuReg}). One can instead prove that the bound
\[
\left\| u(t) \right\|_{\varrho;\nu} \leq C\, t^{-\frac{\varrho-q_0}{\varrho(m-1)}}  \left\| u_0 \right\|_{q_0;\nu^{\phantom{a}}}^{\frac{q_0}{\varrho}} \!\!   \ \ \ \hbox{\rm for a.e.} \,\, t>0,\ \forall q_0\ge1,\ \forall \varrho\in[q_0,+\infty)
\]
holds true and, moreover, is \it equivalent \rm (when $\nu(\Omega)<\infty$) to the validity of the Poincar\'e inequality \eqref{eq: poinIntro}. A similar equivalence between suitable $L^{q_0}$-$L^\varrho$ smoothing effects and the functional inequality \eqref{eq: weakPoi} holds true in the Neumann case too, but it is technically harder. We want to stress that \it none \rm of our results depends upon the validity of weighted Sobolev inequalities, namely functional inequalities of the form
\[\begin{aligned} \left\| v\right\|_{q;\nu} &\leq C_S \left\| \nabla{v} \right\|_{2;\mu}\ \ \ \forall v \in W^{1,2}_0(\Omega;\nu,\mu) \, , \\  \left\| v -\overline{v} \right\|_{q;\nu} &\leq M_S \left\| \nabla{v} \right\|_{2;\mu}\ \ \ \forall v \in W^{1,2}(\Omega;\nu,\mu)  \end{aligned}\]
for a suitable $q$ \it strictly larger than 2\rm, nor on compactness of embeddings of the Sobolev space $W^{1,2}_0(\Omega;\nu,\mu)$ or $W^{1,2}(\Omega;\nu,\mu)$ into $L^2(\Omega;\nu)$. Also, none of the evolutions considered here, for whose generators only a Poincar\'e inequality is assumed to hold, seems to have been treated in the literature so far.

As a further comparison with the results of \cite{RV08,RV09,KRV10}, we notice that when $\rho_\nu(\mathbf{x})=(1+|\mathbf{x}|^2)^{\alpha/2}$, $\rho_\mu(\mathbf{x})=(1+|\mathbf{x}|^2)^{\beta/2}$ and $\beta>N-2$, radial solutions of finite (weighted) energy of the equations at hand are mapped, by a radial change of variable, into radial solutions of finite energy of the WPME corresponding to $\tilde\rho_\nu(\mathbf{y})=(1+|\mathbf{y}|^2)^{\gamma/2}$, $\tilde\rho_\mu(\mathbf{y})\equiv1$, for a suitable $\gamma$. It must be noted that this approach seems to be confined to radial solutions. Besides, this procedure only works in the range of parameters for which Sobolev-type inequalities hold (see Section \ref{sec: exa} for related examples). For these evolutions, one indeed expects $L^\infty$ regularizing effects along the lines of proof given e.g.\ in \cite{BG05m}, so that in a sense the class of equations considered in \cite{RV08, RV09, KRV10} and our present class are strictly disjoint. In fact, we shall show elsewhere that when suitable $L^\infty$ regularizing bounds hold for porous media type evolutions, then appropriate Sobolev inequalities must hold as well,  hence one cannot expect such $L^\infty$ bounds in cases (as the ones we shall deal with) in which no Sobolev inequality is valid.

As already mentioned, the Neumann problem for the WPME is much less studied and, in fact, constitutes the core of the paper. Asymptotic estimates on $u(\cdot,t)$ were already provided by the pioneering work of N.\ D.\ Alikakos and R.\ Rostamian \cite{AR81}, where equation \eqref{eq: pmeNeu} was studied in the case $\rho_\nu=\rho_\mu\equiv1$ on regular domains. Specifically, they proved that if $\overline{u}=0$ then $u(\cdot,t)$ converges uniformly to zero with the sharp rate  $t^{-1/(m-1)}$ \cite[Th. 3.1]{AR81}, while if $\overline{u} \neq 0$ convergence to the mean value (which is preserved along the evolution) is exponential \cite[Th. 3.3]{AR81}. However, such results were proved only for $u_0 \in L^\infty(\Omega)$, and an $L^{q_0}$-$L^\infty$ regularizing effect was proved much later in \cite{BG05m} in the context of evolution on Riemannian manifolds. The most important difference with respect to the present setting lies in the fact that in \cite{BG05m} the validity of the (classical) \emph{Sobolev} inequality was \emph{assumed} to hold in $W^{1,2}(\Omega;\nu,\mu)$, where ${\rm d}\nu={\rm d}\mu$ both coincide with the Riemannian measure on the manifold considered. The lack of any Sobolev embedding will make here, in general, smoothing into $L^\infty$ false. As for the asymptotic behaviour of solutions to \eqref{eq: pmeNeu}, we prove convergence to the mean value in $L^\varrho(\Omega;\nu)$ for all $\varrho\in[1,+\infty)$ and for general data in $L^1(\Omega;\nu)$. Moreover, we show that the rates of such $L^\varrho$ convergence are basically still those given in \cite{AR81}, provided we start from zero-mean $L^1$ data or from nonzero-mean bounded data. In this latter case we also show that global uniform convergence need not hold, but we prove \emph{local} uniform convergence to the mean value with exponential rate. This result is similar to what is shown in the one-dimensional, one-weight case studied in \cite{KR2}, though rates of convergence are not discussed there.

For a thorough analysis of smoothing and decay properties of solutions to large classes of nonlinear evolution equations on ${\mathbb R}^N$, see the monograph \cite{V} (see also \cite{GGS}).
Notice in addition that in \cite{DG+08, DNS, W} one can find decay bounds for weighted porous media equations in connection with functional inequalities, but no regularizing effect is dealt with there.

Other work on the Neumann problem for equations related to the porous media example can be found in \cite{ACLT, AT}, but the discussions there involve domains of infinite volume, so that convergence to zero, rather than to the mean value of the initial datum, takes place, thus giving rise to a situation which is in some sense closer to the Dirichlet case.
\\[3mm]
\noindent\bf Plan of the paper. \rm In Section \ref{sec: wSS} we provide the main definitions of the weighted Sobolev spaces used in the sequel and discuss briefly their main properties. In Section \ref{sec: wellP} we collect all our results concerning existence and uniqueness of solutions, see e.g.\ Theorems \ref{teo: teoFondDir}, \ref{teo: teoFondNeu} and the subsequent applications. Notice that it is well-known, in the weighted case, that non-uniqueness issues may arise even for bounded initial data, this being true even for the linear heat equation (see e.g.\ \cite{Gry} for examples on this latter fact in the manifold setting -- this is related to stochastic completeness of the underlying manifold): therefore uniqueness will have to be understood for solutions in a proper class (see Propositions \ref{pro: uniqPmeDir},  \ref{pro: uniqNeu}). Section \ref{sec: regAs} contains the results on the WPME with Dirichlet boundary conditions. Theorem \ref{teo: regPmeDir} gives the main $L^\varrho$ bounds (for all $\varrho \in [1,\infty)$) as a consequence of the Poincar\'e inequality only, whereas Theorem \ref{teo: pmeDirImpInv} gives the converse result (for $\nu(\Omega)<\infty$) and hence the equivalence between suitable $L^{q_0}$-$L^\varrho$ smoothing effects and \eqref{eq: poinIntro}. No $L^\infty$ regularization holds true in general. Section \ref{sec: pmeNeuReg} contains the results for the harder case, namely the WPME with Neumann boundary conditions. Theorem \ref{teo: regPmeNeu} gives the short-time $L^\varrho$ smoothing effect (for all $\varrho \in [1,\infty)$) as a consequence of \eqref{eq: weakPoi}, and suitable converse implications as well as equivalence results in the same spirit of the Dirichlet case are given in Theorem \ref{teo: impInvNeu} and Corollary \ref{equiv}. Again, no $L^\infty$ regularization holds true in general. As for the long-time asymptotics, the case of data with zero mean is studied in Theorem \ref{teo: absbEMnullaPme}, whereas Theorems \ref{teo: AsiMediaNN}, \ref{thm: convExp}, \ref{cor: convLoc} discuss the other cases. By means of a counterexample which is given in the end of the section, we show that in general one cannot expect uniform convergence to the mean value, even for bounded initial data. We comment that the short time behaviour given in Theorems \ref{teo: regPmeDir} and \ref{teo: regPmeNeu} is \it different \rm from the one valid in the non-weighted case, since the bound proved in such theorems must be valid for general $L^{q_0}$ data, possibly not compactly supported, so that the degeneracy or singularity of the weights at the boundary influences the resulting estimates at all times.

It is important to remark that both our analyses for the Dirichlet and the Neumann problem also apply when $\Omega=\mathbb{R}^N$: if $W^{1,2}_0(\mathbb{R}^N;\nu,\mu) \neq W^{1,2}(\mathbb{R}^N;\nu,\mu)$ the two problems, in general, will be different.

Section \ref{sec: exa} collects some examples of allowable weights both for the Dirichlet and the Neumann setting. Here we certainly do not even try to make a comprehensive account of the various conditions ensuring the validity of weighted Poincar\'e inequalities. Instead, we shall provide the reader with explicit examples of weights, hence of evolutions, for which weighted Poincar\'e inequalities hold but weighted Sobolev inequalities do not, thus showing cases in which only the present results can be used to investigate smoothing and decay properties of the associated evolutions. In Section \ref{list}, for the convenience of the reader, we shall list concisely some significant examples.
\end{section}

\begin{section}{Weighted Sobolev spaces and Poincar\'e inequalities}\label{sec: wSS}
We recall here some basic definitions and facts about weighted Sobolev spaces and related Poincar\'e inequalities. In the sequel $\Omega \subseteq \mathbb{R}^N$ is a domain and $\nu$ and $\mu$ are two measures absolutely continuous with respect to the Lebesgue measure: let $\rho_{\nu}$ and $\rho_{\mu}$ be the corresponding weights (or densities). We shall always assume
$$\rho_{\nu}(\mathbf{x}), \rho_{\mu}(\mathbf{x}) >0 \ \ \ \textnormal{for a.e. $\mathbf{x} \in \Omega \, ,$}$$
so that also the Lebesgue measure is absolutely continuous with respect to $\nu$ and $\mu$. For all $p \in [1,\infty)$ we denote as $L^p(\Omega;\nu)$ the Banach space of equivalence classes of Lebesgue-measurable functions $f$ such that
$$ \left\| f \right\|_{p;\nu}^p = \int_{\Omega} |f|^p \, \mathrm{d}\nu = \int_{\Omega} |f|^p \, \rho_{\nu}  \mathrm{d}\mathbf{x} < \infty \, . $$
We define the weighted Sobolev space $W^{1,p}(\Omega; \nu, \mu)$ (see e.g.\ \cite{KO84}) as the set of all (equivalence classes of) functions $v \in W^{1,1}_{loc}(\Omega)$ such that
\begin{equation*}
\left\| v \right\|^p_{p;\nu,\mu} = \left\| v \right\|^p_{p;\nu} + \left\| \nabla v \right\|^p_{p;\mu} < \infty \, .
\end{equation*}

Yet without further assumptions on $\rho_{\nu}$ and $\rho_{\mu}$ in general $W^{1,p}(\Omega; \nu, \mu)$ would not be complete.
\begin{definition} \label{def: Bp}
For all $p \in (1,\infty)$ we denote as $B^p(\Omega)$ the class of all measurable functions $f$ such that $f > 0$ a.e.\ and
\begin{equation*}
\left|f\right|^{\frac{1}{1-p}} \in L^1_{loc}(\Omega) \, .
\end{equation*}
\end{definition}
One can prove \cite[Th. 2.1]{KO84} that if $p \in (1,\infty)$ and $\rho_{\mu} \in B^p(\Omega)$ then $W^{1,p}(\Omega; \nu, \mu)$ is indeed complete. If $p=1$ the same result is true providing that the condition $\rho_{\mu} \in B^p(\Omega)$ is replaced by $\rho_{\mu}^{-1} \in L^{\infty}_{loc}(\Omega)$.\\
The fact that for any $\varphi \in C^{\infty}_c(\Omega)$ the quantity \mbox{$\| \varphi \|_{p;\nu,\mu} $} is finite is equivalent (see \cite[Lem. 4.4]{KO84}) to the local finiteness of $\nu$ and $\mu$, that is
\begin{equation}\label{eq: condizioneFin}
\rho_{\nu}, \rho_{\mu}  \in L^1_{loc}(\Omega)  \, .
\end{equation}
We then define, provided \eqref{eq: condizioneFin} holds, the space $W^{1,p}_0(\Omega; \nu, \mu)$ as the closure of $C^{\infty}_c(\Omega)$ with respect to the norm \mbox{$\| \cdot \|_{p;\nu,\mu} $}.

When dealing with $W^{1,p}(\Omega; \nu, \mu)$ [$W^{1,p}_0(\Omega; \nu, \mu)$] we shall always assume, without further comment, $\rho_{\mu} \in B^p(\Omega)$ [$\rho_{\mu} \in B^p(\Omega)$ and \eqref{eq: condizioneFin}].

In the following, we list some elementary properties of the spaces defined above.
\begin{proposition} \label{pro: inc1}
Let $p \in [1,\infty)$. The inclusion $W^{1,\infty}_c(\Omega) \subset W_0^{1,p}(\Omega;\nu,\mu)$ holds.
\begin{proof}
Given $v \in W^{1,\infty}_c(\Omega)$, thanks to \cite[Lemmas 2.18 and 3.15]{Ada75} we know that there exists a sequence of functions $\{v_n\} \subset C^{\infty}_c(\Omega)$ (the \emph{mollification} of $v$) such that, as $n \rightarrow \infty$,
$$ v_n \xrightarrow{\textnormal{a.e.}} v \, , \ \nabla{v_n} \xrightarrow{\textnormal{a.e.}} \nabla{v} \, , \ \left\| v_n \right\|_{{\infty}} \leq \left\| v \right\|_{{\infty}} \, , \  \left\| \nabla{v_n} \right\|_{{\infty}} \leq \left\| \nabla{v} \right\|_{{\infty}} $$
and $\operatorname{supp}(v_n) \subset \Omega^{\prime} \Subset \Omega $. Being $\nu$ e $\mu$ locally finite, $v_n \rightarrow v$ in $W^{1,p}(\Omega;\nu,\mu)$ by the dominated convergence Theorem.
\end{proof}
\end{proposition}
\begin{proposition} \label{pro: dens1}
Let $p \in [1,\infty)$. If $\rho_{\nu}, \rho_{\mu} , \rho^{-1}_{\nu} , \rho^{-1}_{\mu} \in L^\infty_{loc}(\Omega) $, then $C^{\infty}(\Omega) \cap W^{1,p}(\Omega;\nu,\mu)$ is dense in $W^{1,p}(\Omega;\nu,\mu)$.
\begin{proof}
The assumptions imply that the weighted norms $\| \cdot \|_{p;\nu}$ and $\| \cdot \|_{p;\mu}$ are \emph{locally} equivalent to the non-weighted norm $\| \cdot \|_p$. This is enough to reproduce the proof of \cite[Th. 3.16]{Ada75}.
\end{proof}
\end{proposition}
\begin{proposition} \label{pro: dens2}
For any $p \in [1,\infty)$ the space $L^{\infty}(\Omega) \cap W^{1,p}(\Omega; \nu, \mu)$ is dense in $W^{1,p}(\Omega; \nu, \mu)$.
\begin{proof}
Given $v \in W^{1,p}(\Omega; \nu, \mu)$, as well as for the non-weighted case, consider the approximating sequence of functions $v_n=\min\left(n,\max\left(-n,v \right)\right) $. By construction, $\{v_n\} \subset L^{\infty}(\Omega) \cap W^{1,p}(\Omega; \nu, \mu)$ and \mbox{$|v_n|\leq|v|$}; moreover, $\nabla{v_n}=\left(\nabla{v} \right) \chi_{\{-n < v < n \}}$. The assertion then follows by monotone convergence.
\end{proof}
\end{proposition}
Let us now introduce other useful weighted Sobolev spaces that we shall deal with throughout the discussion.
\begin{definition}\label{def: spaz0}
Given $p\in[1,\infty)$, $\rho_\nu,\rho_\mu \in L^1_{loc}(\Omega) $ and $\rho_\mu \in B^2(\Omega)$, let $V^p_0(\Omega;\nu,\mu)$ be the closure of $C^\infty_c(\Omega)$ with respect to the norm
\begin{equation*}
\left\| \varphi \right\|_{p,2;\nu,\mu} = \left\| \varphi \right\|_{p;\nu} + \left\| \nabla{\varphi} \right\|_{2;\mu}
\end{equation*}
and $V_0(\Omega;\mu)$ the space of all functions $v \in W^{1,1}_{loc}(\Omega)$ such that $\nabla v \in L^2(\Omega;\mu)$ and for which there exists a sequence $\{\varphi_n\} \subset C^\infty_c(\Omega)$ such that
$$\left\| \nabla{v}-\nabla{\varphi_n} \right\|_{2;\mu}\rightarrow 0 \, . $$
\end{definition}
\noindent Clearly, $V^p_0(\Omega;\nu,\mu)$ is a (reflexive if in addition $p>1$) Banach space.
\begin{definition}\label{def: spaz1}
Given $p\in[1,\infty)$ and $\rho_\mu \in B^2(\Omega)$, we denote as $V^p(\Omega;\nu,\mu)$ the space of all functions $v \in W^{1,1}_{loc}(\Omega)$ such that
\begin{equation*}
\left\| v \right\|_{p,2;\nu,\mu} = \left\| v \right\|_{p;\nu} + \left\| \nabla{v} \right\|_{2;\mu} < \infty \, .
\end{equation*}
\end{definition}
\noindent $V^p(\Omega;\nu,\mu)$ is also a (reflexive if in addition $p>1$) Banach space.

\noindent Finally, we mention two elementary properties of weighted Poincar\'e inequalities. If $\nu(\Omega)<\infty$ we shall denote as $\overline{f}$ the weighted mean value of any function $f\in L^1(\Omega;\nu)$, that is
\begin{equation}\label{eq: meanV}
\overline{f}=\frac{\int_\Omega f \, \mathrm{d}\nu}{\nu(\Omega)}  \, .
\end{equation}
\begin{proposition}\label{pro: Edm}
Suppose that $\nu(\Omega)<\infty$. The validity of the zero-mean $p$-Poincar\'e inequality
\begin{equation*}
\left\| v-\overline{v}  \right\|_{p;\nu} \leq M_P \left\| \nabla{v} \right\|_{p;\mu} \ \ \ \forall v \in W^{1,p}(\Omega;\nu,\mu)
\end{equation*}
for a suitable $M_P>0$ is equivalent to the validity of the inequality
\begin{equation*}
\inf_{c \in \mathbb{R}}{\left\| v-c \right\|_{p;\nu}} \leq M_I \left\| \nabla{v} \right\|_{p;\mu} \ \ \ \forall v \in W^{1,p}(\Omega;\nu,\mu) \,
\end{equation*}
for a suitable $M_I>0$.
\begin{proof}
See \cite[Lem. 3.1]{EO93}.
\end{proof}
\end{proposition}
\begin{proposition}\label{pro: OrdPoin}
Let $(\rho_{\nu_1},\rho_{\mu_1})$ and $(\rho_{\nu_2},\rho_{\mu_2})$ be two couples of weights. Suppose that there exist two constants $D_\nu >0$ and $D_\mu>0$ such that
$$ \rho_{\nu_2} \leq D_\nu \, \rho_{\nu_1} \, , \ \rho_{\mu_1} \leq D_\mu \, \rho_{\mu_2} \, .$$
Then if $W^{1,p}_0(\Omega;\nu_1,\mu_1)$ satisfies the $p$-Poincar\'e inequality (with $\nu=\nu_1$ and $\mu=\mu_1$)
$$ \left\| v \right\|_{p;\nu} \leq M_{P_1} \left\| \nabla{v} \right\|_{p;\mu} $$
so does $W^{1,p}_0(\Omega;\nu_2,\mu_2)$ (with $\nu=\nu_2$ and $\mu=\mu_2$, up to multiplicative constants). Similarly, if $W^{1,p}(\Omega;\nu_1,\mu_1)$ satisfies the zero-mean $p$-Poincar\'e inequality (with $\nu=\nu_1$ and $\mu=\mu_1$)
$$ \left\| v-\overline{v}  \right\|_{p;\nu} \leq M_{P_2} \left\| \nabla{v} \right\|_{p;\mu} $$
so does $W^{1,p}(\Omega;\nu_2,\mu_2)$ (with $\nu=\nu_2$ and $\mu=\mu_2$, up to multiplicative constants).
\begin{proof}
It is immediate to prove that the $p$-Poincar\'e inequality, with respect to $(\nu_2,\mu_2)$, holds in $C^\infty_c(\Omega)$. By density such property is extended to the whole $W^{1,p}_0(\Omega;\nu_2,\mu_2)$. For the zero-mean $p$-Poincar\'e inequality one argues likewise, taking $L^\infty(\Omega) \cap W^{1,p}(\Omega;\nu_2,\mu_2)$ as a dense space (Proposition \ref{pro: dens2}). In this case it is also convenient to exploit Proposition \ref{pro: Edm}.
\end{proof}
\end{proposition}
\end{section}

\begin{section}{Well-posedness of the problems}\label{sec: wellP}
In this section we provide some existence and uniqueness results for solutions to the previously mentioned evolutions, whose smoothing and asymptotic properties we shall study in detail in Sections \ref{sec: regAs} and \ref{sec: pmeNeuReg}.
\begin{subsection}{The {WPME} with Dirichlet boundary conditions}\label{sec: pmeDir}
We begin with giving our notion of weak solution to \eqref{eq: pmeDir}. Although we shall not mention it explicitly any further, we comment that the present results hold, with no modifications, in the case $\Omega={\mathbb R}^N$ as well.
\begin{definition} \label{den: solDebDir}
A function
$$u \in L^1((0,T);L^1_{loc}(\Omega;\nu)): \ u^m(t) \in V_0(\Omega;\mu)\, , \ \nabla{(u^m)} \in L^1((0,T);[L^2(\Omega;\mu)]^N) $$
$$\textnormal{for a.e.} \ t>0 \textnormal{ and } \forall T>0 \, ,$$
is a weak solution of \eqref{eq: pmeDir} with initial datum $u_0 \in L^1_{loc}(\Omega;\nu)$ if it satisfies:
\begin{equation}
\begin{aligned}\label{eq: solDebDir}
& \int_0^T \! \! \int_{\Omega} u(\mathbf{x},t) \eta_t(\mathbf{x},t) \, \mathrm{d}\nu \, \mathrm{d}t \\
= & -\int_\Omega u_0(\mathbf{x}) \eta(\mathbf{x},0) \,  \mathrm{d}\nu  + \int_0^T \! \! \int_{\Omega} \nabla{\left( u^m\right)}(\mathbf{x},t) \cdot \nabla{\eta}(\mathbf{x},t) \, \mathrm{d}\mu \, \mathrm{d}t
\end{aligned}
\end{equation}
$$ \forall \eta \in C^{1}(\Omega \times [0,T] ): \ \operatorname{supp}{\eta(\cdot,t)} \Subset \Omega \, , \ \eta(\mathbf{x},T)=0 \ \ \ \forall \mathbf{x} \in \Omega \, , \ \forall t \in [0,T]  \, . $$
\end{definition}
Such notion is very similar to the one given in \cite[Def. 5.4]{Vaz07} (non-weighted porous media equation on bounded domains). The main difference lies in the fact that, having to deal with general domains and weights, it seemed reasonable for us not to require any further a priori integrability property for $u^m$.

The next uniqueness result is the equivalent of Theorem 5.3 of \cite{Vaz07}.
\begin{proposition}\label{pro: uniqPmeDir}
There exists \emph{at most} one weak solution of \eqref{eq: pmeDir} satisfying the following additional hypotheses:
\begin{gather}\label{eq: hpEnerg}
u^m \in L^{\frac{m+1}{m}}((0,T); V^{\frac{m+1}{m}}_0(\Omega;\nu,\mu)) \, , \ \nabla{(u^m)} \in L^2((0,T);[L^2(\Omega;\mu)]^N)  \\
\forall T>0  \, .  \nonumber
\end{gather}
\begin{proof}
We use the method of proof of \cite[Th. 5.3]{Vaz07}. In particular, thanks to \eqref{eq: hpEnerg} and to a density argument, it is possible to choose in \eqref{eq: solDebDir} any test function $\eta$ such that
$$ \eta \in W^{1,\frac{m+1}{m}}((0,T); V^{\frac{m+1}{m}}_0(\Omega;\nu,\mu))\, , \  \nabla{\eta} \in L^2((0,T);[L^2(\Omega;\mu)]^N)  \, , \ \eta(T)=0 \, . $$
The assertion follows as in the quoted proof by plugging Ole\u{\i}nik's test function \cite{Ole1, Ole2}
$$ {\eta}(t)= \int_t^T \left( u_1^m(s) - u_2^m(s)  \right) \, \mathrm{d}s $$
into the weak formulation satisfied by the difference of two possible solutions $u_1$ and $u_2$ fulfilling \eqref{eq: hpEnerg} and performing analogous computations.
\end{proof}
\end{proposition}
According to a common terminology used in \cite{Vaz07}, we give the following definition.
\begin{definition}\label{wes-dirichlet}
 We shall call (weak) \emph{energy solutions} all weak solutions to \eqref{eq: pmeDir} that also satisfy \eqref{eq: hpEnerg}.
\end{definition}
In order to establish a suitable existence theorem, we first need to prove a fundamental lemma.  Hereafter, by saying that a domain $\Omega$ is smooth we shall mean, without further comment, that it is at least $C^{2,\alpha}$.
\begin{lemma} \label{lem: lemFondDir}
If one assumes that $\Omega $ is a smooth bounded domain of $\mathbb{R}^N$, $\rho_\nu \in C^{3,\alpha}(\overline{\Omega})$, $\rho_\mu \in C^{2,\alpha}(\overline{\Omega})$, $\rho_\nu^{-1}, \rho_\mu^{-1} \in L^{\infty}(\Omega)$ and $u_0 \in C^{2,\alpha}_c(\Omega)$, then there exists a weak solution $u$ of \eqref{eq: pmeDir} which satisfies, for almost every $T>0$ and every $q \geq 0$, the following estimates:
\begin{equation}\label{eq: stimeLemmaDir1}
\begin{aligned}
& \frac{4q(q+1)m}{(m+q)^2} \int_0^T \! \! \int_\Omega \left|\nabla{\left(u^{\frac{m+q}{2}} \right)}(\mathbf{x},t)\right|^2 \,   \mathrm{d}\mu  \, \mathrm{d}t + \int_\Omega \left|u(\mathbf{x},T)\right|^{q+1} \, \mathrm{d}\nu \\
 \leq & \int_\Omega \left|u_0(\mathbf{x})\right|^{q+1} \, \mathrm{d}\nu \, ,
\end{aligned}
\end{equation}
\begin{equation}\label{eq: stimeLemmaDir2}
\int_0^T \! \! \int_\Omega \zeta(t) \left[ \left(u^{\frac{m+1}{2}} \right)_t (\mathbf{x},t) \right]^2 \, \mathrm{d}\nu \, \mathrm{d}t  \leq \max_{t \in [0,T]}{ \zeta^\prime(t)} \frac{m+1}{8m} \int_{\Omega} |u_0(\mathbf{x})|^{m+1} \, \mathrm{d}\nu   \, ,
\end{equation}
where $\zeta \geq 0$ is any $C^1_c(0,T)$ function.\\
Moreover if $v$ is another weak solution, obtained with the same approximating scheme of the incoming proof (see \eqref{eq: pmeDirParab}), corresponding to an initial datum $v_0 \in  C^{2,\alpha}_c(\Omega)$, the inequality
\begin{equation}\label{eq: princConf}
\int_\Omega (u(\mathbf{x},T)-v(\mathbf{x},T))_{+} \, \mathrm{d}\nu \leq \int_\Omega (u_0(\mathbf{x})-v_0(\mathbf{x}))_{+} \, \mathrm{d}\nu \,
\end{equation}
holds for almost every $T>0$. In particular, the comparison principle holds.
\begin{proof}
We proceed along the lines of the proof of \cite[Lem. 5.8]{Vaz07}, where a first existence result for the non-weighted porous media equation is established. The essential idea is to approximate problem \eqref{eq: pmeDir} with non-degenerate problems. As a first step we pick a sequence $\Phi_n^\prime(x):\mathbb{R}\rightarrow \mathbb{R}$ of smooth functions such that:
\begin{itemize}
\item[$\bullet$] $\Phi_n^\prime(x) \rightarrow m \, |x|^{m-1}$ locally uniformly;
\item[$\bullet$] $\Phi_n^\prime(x) > 0 \ \ \ \forall x \in \mathbb{R}$;
\item[$\bullet$] $\Phi_n^\prime(x)=\Phi_n^\prime(-x) $, so that in particular $\Phi_n(0)=0$, where $\Phi_n(x)=\int_0^x \Phi_n^\prime(y) \, \mathrm{d}y $.
\end{itemize}
Now, consider the following non-degenerate (thanks to the properties of $\Phi_n^\prime$) quasilinear problem:
\begin{equation} \label{eq: pmeDirParab}
\begin{cases}
(u_n)_t =\rho_{\nu}^{-1}  \operatorname{div}\left(\rho_{\mu} \,  \nabla{\left( \Phi_n(u_n) \right) } \right) & \textnormal{in} \  \Omega\times(0,\infty)  \\
 u_n =0   & \textnormal{on} \ \partial\Omega \times (0,\infty) \\
 u_n(\cdot,0)=u_0(\cdot)  & \textnormal{in} \ \Omega
 \end{cases} \, .
\end{equation}
Performing the change of variable $w=\rho_\nu u_n$ it is convenient to write the latter in divergence form:
\begin{equation} \label{eq: pmeDirParabDiv}
\begin{cases}
w_t = \operatorname{div}\left(\frac{\rho_\mu}{\rho_\nu} \, \Phi_n^\prime\left(\frac{w}{\rho_\nu}\right) \nabla{w} - \frac{\rho_\mu}{\rho_\nu^2} \, \nabla{(\rho_\nu)} \, \Phi_n^\prime\left(\frac{w}{\rho_\nu}\right)w \right) & \textnormal{in} \  \Omega\times(0,\infty)  \\
 w =0  &  \textnormal{on} \ \partial\Omega \times (0,\infty) \\
 w(\cdot,0)=\rho_\nu(\cdot) u_0(\cdot) & \textnormal{in} \ \Omega
 \end{cases}\, .
\end{equation}
There is no loss of generality in assuming that, given $\epsilon > 0$, $\Phi_n^\prime(x)=c$ for $|x| \geq \|u_0\|_{\infty} + \epsilon$, $c$ being a suitable positive constant possibly depending on $n$. Under these hypotheses, Theorem V.6.1 of \cite{LSU68} is applicable, which provides us with a solution $w(\mathbf{x},t) \in C^{2,1}(\overline{\Omega} \times [0,T]) \  \forall T>0 $ of \eqref{eq: pmeDirParabDiv}; from standard parabolic regularity results \cite[Th. IV.5.2]{LSU68} we also have, in particular, $w_t(\mathbf{x},t) \in C^{1,0}(\Omega \times (0,T))$. Hence $u_n(\mathbf{x},t)$ is a solution of \eqref{eq: pmeDirParab} as regular as $w(\mathbf{x},t)$. Moreover, thanks to the parabolic maximum principle \cite[Th. I.2.9]{LSU68} we have
$$
\| u_n(T) \|_\infty \leq \| u_0 \|_\infty \quad \forall\,\, T>0 \, .
$$
Given a function $\eta$ as in the weak formulation \eqref{eq: solDebDir}, multiplying \eqref{eq: pmeDirParab} by $\rho_\nu \eta$ and integrating by parts in $\Omega \times (0,T)$, we get:
\begin{equation}\label{eq: QuasiSolDebDir}
\begin{aligned}
\int_0^T \! \! \int_{\Omega} u_n(\mathbf{x},t) \eta_t(\mathbf{x},t) \, \mathrm{d}\nu \, \mathrm{d}t = &-\int_\Omega u_0(\mathbf{x}) \eta(\mathbf{x},0) \,  \mathrm{d}\nu\\  &+ \int_0^T \! \! \int_{\Omega} \nabla{\left( \Phi_n(u_n) \right)}(\mathbf{x},t) \cdot \nabla{\eta}(\mathbf{x},t) \, \mathrm{d}\mu \, \mathrm{d}t \, .
\end{aligned}
\end{equation}
In order to pass to the limit in \eqref{eq: QuasiSolDebDir} as $n\rightarrow \infty$ we must obtain suitable estimates on $u_n$ and $\nabla{\left( \Phi_n(u_n) \right)}$ and afterwards identify weak limits. Setting
$$ \Psi_n(x)=\int_0^x \Phi_n(y) \, \mathrm{d}y \, , \ \Upsilon^1_n(x)=\int_0^x \sqrt{\Phi_n^\prime(y)} \, \mathrm{d}y \, ,  $$
through computations similar to the ones developed in \cite[Lem. 5.8]{Vaz07} and exploiting the spatial regularity of $u_t$ we arrive at:
\begin{equation}\label{eq: stimaEnerg1}
\int_\Omega \Psi_n(u_n(\mathbf{x},T)) \, \mathrm{d}\nu + \int_0^T \! \! \int_\Omega \left| \nabla{(\Phi_n(u_n))}(\mathbf{x},t) \right|^2  \, \mathrm{d}\mu \, \mathrm{d}t = \int_\Omega \Psi_n(u_0(\mathbf{x})) \, \mathrm{d}\nu \, ,
\end{equation}
\begin{equation}\label{eq: stimaEnerg2}
\frac{1}{2} \int_\Omega u_n^2(\mathbf{x},T) \, \mathrm{d}\nu + \int_0^T \! \! \int_\Omega \left| \nabla{(\Upsilon^1_n(u_n))}(\mathbf{x},t) \right|^2  \, \mathrm{d}\mu \, \mathrm{d}t =\frac{1}{2} \int_\Omega u_0^2(\mathbf{x}) \, \mathrm{d}\nu \, ,
\end{equation}
\begin{equation}\label{eq: stimaDerivataT}
\int_0^T \! \! \int_\Omega \zeta(t) \left[\left( \Upsilon^1_n(u_n)\right)_t (\mathbf{x},t) \right]^2 \, \mathrm{d}\nu \, \mathrm{d}t = \int_0^T \! \! \int_\Omega \frac{\zeta^\prime(t)}{2} \left| \nabla{(\Phi_n(u_n))}(\mathbf{x},t) \right|^2 \, \mathrm{d}\nu \, \mathrm{d}t \, ,
\end{equation}
where $\zeta \geq 0$ is any $C^1_c(0,T)$ function. From \eqref{eq: stimaEnerg2}--\eqref{eq: stimaDerivataT}, the maximum principle and the inner regularity of $\rho_\nu$ and $\rho_\mu$ (in particular, here it is crucial that the weights are locally equivalent to $1$), we deduce that $\{\Upsilon^1_n(u_n)\}$ is locally bounded in $H^1(\Omega \times (0,\infty))$. Therefore, up to a subsequence, $\{\Upsilon^1_n(u_n)\}$ converges a.e.\ in $\Omega \times (0,\infty)$. This easily implies (by the smoothness of the approximating sequence of functions $\{\Phi_n^\prime\}$) the existence of a function $u$ such that:
\begin{equation}\label{eq: stimaLimPunt}
\begin{aligned}
&u_n \rightarrow u \, , \  \Psi_n(u_n) \rightarrow \frac{1}{m+1}|u|^{m+1}  \,, \ \Phi_n(u_n) \rightarrow u^m \, , \\
&\Upsilon^1_n(u_n) \rightarrow \frac{2\sqrt{m}}{m+1} \, u^{\frac{m+1}{2}} \ \ \ \mathrm{a.e.\ in} \ \Omega \times (0,\infty) \, .
\end{aligned}
\end{equation}
The maximum principle, estimates \eqref{eq: stimaEnerg1}, \eqref{eq: stimaDerivataT} and the pointwise limits given in \eqref{eq: stimaLimPunt} permit to conclude (again along a subsequence) that
$$ u_n \rightarrow u \textnormal{ \ in \ }  L^2((0,T);L^2(\Omega)) \, , \ \Phi_n(u_n) \rightharpoonup u^{m} \textnormal{ \ in \ } L^2((0,T);H^1_0(\Omega)) \, , $$
$$ \Upsilon^1_n(u_n) \rightharpoonup \frac{2\sqrt{m}}{m+1}  \, u^{\frac{m+1}{2}} \textnormal{ \ in \ }  W^{1,2}((\tau,T);L^2(\Omega)) \ \ \ \forall \tau \in (0,T) \, ;$$
hence, passing to the limit in \eqref{eq: QuasiSolDebDir} as $n\rightarrow \infty$, we conclude that $u$ is a weak solution of \eqref{eq: pmeDir} (in the sense of Definition \ref{den: solDebDir}) with initial datum $u_0$.\\
Finally, we must obtain \eqref{eq: stimeLemmaDir1}, \eqref{eq: stimeLemmaDir2} and \eqref{eq: princConf}. The first one follows (at least for $q\geq1$) by  multiplying \eqref{eq: pmeDirParab} by $\rho_\nu \, u_n^{q}$, integrating in $\Omega \times (0,T)$ and suitably passing to the limit. If $q\in(0,1)$ things are slightly more technical. However, it is only a question of approximating $x^q$  with a sequence of regular functions. To retrieve the case $q=0$ one lets $q\downarrow 0$. See also \cite[Prop. 5.12]{Vaz07}. Estimate \eqref{eq: stimeLemmaDir2} is a direct consequence of \eqref{eq: stimaDerivataT}, \eqref{eq: stimaEnerg1} and the weak convergence of $\{\Upsilon^1_n(u_n)\}$ to $\frac{2\sqrt{m}}{m+1} \, u^{(m+1)/2} $ in $W^{1,2}((\tau,T);L^2(\Omega))$. Inequality \eqref{eq: princConf} can be obtained exactly as  in \cite[Prop. 3.5]{Vaz07} by using in addition an approximation procedure, see \cite[Prop. 6.1]{Vaz07}.
\end{proof}
\end{lemma}
Starting from the previous lemma, we are able to prove existence of weak energy solutions when the initial datum $u_0$, the domain $\Omega$ and the weights $\rho_\nu,\rho_\mu$ are less regular.
\begin{theorem}\label{teo: teoFondDir}
Let $\Omega \subset \mathbb{R}^N$ be a domain, and let $\rho_\nu,\rho_\mu$ be two weights such that
\begin{equation*}
\rho_\nu \in C^{3,\alpha}_{loc}(\Omega) \, , \ \rho_\mu \in C^{2,\alpha}_{loc}(\Omega) \, , \ \rho_\nu^{-1},\rho_\mu^{-1} \in L^{\infty}_{loc}(\Omega) \, .
\end{equation*}
If $u_0 \in L^1(\Omega;\nu) \cap L^{r}(\Omega;\nu) $, with $r \geq m+1$, then there exists a weak solution $u$ of \eqref{eq: pmeDir} which satisfies estimate \eqref{eq: stimeLemmaDir1} for all $0\leq q \leq r-1$, estimate \eqref{eq: stimeLemmaDir2} and it is the unique energy solution in the sense of Definition \ref{wes-dirichlet}. Moreover, if $v$ is the energy solution corresponding to another initial datum $v_0 \in L^{1}(\Omega;\nu) \cap L^{m+1}(\Omega;\nu)$, inequality \eqref{eq: princConf} and in particular the comparison principle still hold.
\begin{proof}
To extend the results of Lemma \ref{lem: lemFondDir} to general $L^{\infty}$ data, it suffices to approximate $u_0$ with a sequence $\{u_{0n}\}$ of regular data and check appropriate convergence of the corresponding (sub)sequence of solutions $\{u_n\}$. To weaken the hypotheses on the domain $\Omega$ and on the weights $\rho_\nu,\rho_\mu$ one can proceed similarly to the end of the proof of \cite[Th. 5.7]{Vaz07}: assuming $u_0 \in L^1(\Omega;\nu) \cap L^\infty(\Omega) $, one picks an increasing sequence of bounded smooth domains $\{\Omega_n\} $ approximating $\Omega$, with $\Omega_n \Subset \Omega$, and solves on each of them the Dirichlet problem with initial datum $u_{0n}=u_0|_{\Omega_n}$, letting $u_n$ be the relative solution extended to be zero outside $\Omega_n$. Estimates \eqref{eq: stimeLemmaDir1} and \eqref{eq: stimeLemmaDir2} read as follows:
\begin{equation*}
\begin{aligned}
& \frac{4q(q+1)m}{(m+q)^2} \int_0^T\hskip-7pt\int_\Omega \left|\nabla{\left(u_n^{\frac{m+q}{2}} \right)}(\mathbf{x},t)\right|^2\hskip-3pt\mathrm{d}\mu  \, \mathrm{d}t + \int_\Omega \left|u_n(\mathbf{x},T)\right|^{q+1}\hskip-3pt \mathrm{d}\nu \\
\leq & \int_{\Omega} \left|u_{0}(\mathbf{x})\right|^{q+1} \hskip-3pt\mathrm{d}\nu \, ,
\end{aligned}
\end{equation*}
\begin{equation*}
\int_0^T \,\int_\Omega \zeta(t) \left[ \left(u_n^{\frac{m+1}{2}} \right)_t (\mathbf{x},t) \right]^2 \, \mathrm{d}\nu \, \mathrm{d}t  \leq \max_{t \in [0,T]}{ \zeta^\prime(t)} \frac{m+1}{8m} \int_{\Omega} |u_{0}(\mathbf{x})|^{m+1} \, \mathrm{d}\nu  \, .
\end{equation*}
As in the proof of Lemma \ref{lem: lemFondDir}, no major difficulty arises in showing that $\{u_n \}$, up to subsequences, converges to a weak energy solution $u$ of problem \eqref{eq: pmeDir} (with initial datum $u_0$) in such a way that estimates \eqref{eq: stimeLemmaDir1}, \eqref{eq: stimeLemmaDir2} and inequality \eqref{eq: princConf} (taking in addition another sequence of solutions $\{v_n\}$) are preserved.

In order to remove the hypothesis $u_0 \in L^1(\Omega;\nu) \cap L^\infty(\Omega) $, one picks a sequence of initial data $\{u_{0n}\} \subset L^1(\Omega;\nu) \cap L^\infty(\Omega)$ converging to $u_0$ in $L^1(\Omega;\nu) \cap L^r(\Omega;\nu)$ and considers the corresponding sequence $\{u_n\}$ of solutions to \eqref{eq: pmeDir}. Thanks to inequality \eqref{eq: princConf}, $\{u_n\}$ is Cauchy in $L^\infty((0,\infty);L^1(\Omega;\nu))$, so that it converges to a function $u$ belonging to the same space. The stability of the inequality \eqref{eq: princConf} as $n\rightarrow\infty$ is trivial, while the stability of estimates \eqref{eq: stimeLemmaDir1}, \eqref{eq: stimeLemmaDir2} is proved by arguing similarly to the proof of Lemma \ref{lem: lemFondDir}. In fact one can show that
\begin{equation*}\label{eq: convDebLr}
u_n^{\frac{m+q}{2}} \rightharpoonup u^{\frac{m+q}{2}} \textnormal{ \ in \ } L^2((0,T);V_0^{\frac{2r_1}{m+q}}(\Omega;\nu,\mu))
\end{equation*}
for any $r_1 \in \left({(m+q)}/{2},r\right]$ and $q\in(0, r-1]$.
\end{proof}
\end{theorem}
Notice that from estimate \eqref{eq: stimeLemmaDir2} we know, for instance, that $u^{(m+1)/2}$ is absolutely continuous in $C([\tau,\infty);L^2(\Omega;\nu))$ (for any $\tau>0$), which in particular implies $u \in C([\tau,\infty);L^{m+1}(\Omega;\nu)) $. By that, it is not difficult to prove the validity of the so called \emph{semigroup property}: for any $\tau > 0$, $u|_{[\tau,\infty)}$ is the (weak energy) solution of \eqref{eq: pmeDir} with initial datum $u(\tau)$.

When $u_0$ is smooth enough, we are able to ensure that $u^{(m+1)/2}$ is continuous even down to $t=0$. In fact, we have the following (see also \cite[Sec. 5.6]{Vaz07})
\begin{corollary}\label{cor: contZero}
If, together with the hypotheses of Theorem \ref{teo: teoFondDir}, one assumes that
$$u_0^m \in V_0^{\frac{m+1}{m}}(\Omega;\nu,\mu) \, ,$$
then for almost every $T>0$ the estimate
\begin{equation}\label{eq: stimaRegIniziale}
\begin{aligned}
&\int_0^T \! \! \int_\Omega \left[ \left( u^{\frac{m+1}{2}} \right)_t (\mathbf{x},t) \right]^2  \mathrm{d}\nu  \mathrm{d}t  +  \frac{(m+1)^2}{8m} \int_\Omega \left| \nabla{\left( u^m \right)} (\mathbf{x},T) \right|^2  \mathrm{d}\mu \\ \leq & \frac{(m+1)^2}{8m} \int_\Omega \left| \nabla{\left( u_0^m \right)} (\mathbf{x}) \right|^2  \mathrm{d}\mu
\end{aligned}
\end{equation}
holds. In particular, $u^{\frac{m+1}2}$ is an absolutely continuous curve in $C([0,\infty);L^2(\Omega;\nu))$.
\begin{proof}
One can proceed along the lines of the proofs of Lemma \ref{lem: lemFondDir} and Theorem \ref{teo: teoFondDir}: here the fundamental estimate to exploit is \eqref{eq: stimaDerivataT} (up to choosing $\zeta(t)=\chi_{[0,T]}(t) $). We omit the details.
\end{proof}
\end{corollary}
The analysis of the $L^1$-continuity of solutions in $t=0$, when $u_0$ only belongs to $L^{1}(\Omega;\nu)\cap L^{m+1}(\Omega;\nu)$, is not straightforward. We just mention that in the non-weighted case it is proved by means of an explicit (and technical) initial barrier argument \cite[Th. 6.2 and Sec. 7.5.1]{Vaz07}. Nevertheless if $\nu(\Omega)<\infty$ it is a direct consequence of Corollary \ref{cor: contZero} and \eqref{eq: princConf}.

When dealing with initial data in $ L^1(\Omega;\nu)$ with no further integrability properties, we are not able to provide a weak solution of \eqref{eq: pmeDir} in the sense of Definition \ref{den: solDebDir}. However, from \eqref{eq: princConf} we trivially have that the map $L^1(\Omega;\nu) \cap L^{m+1}(\Omega;\nu) \rightarrow L^\infty((0,\infty);L^1(\Omega;\nu)) $, which associates to an initial datum $u_0$ the corresponding energy solution $u(\cdot)$, is Lipschitz and densely defined in $L^1(\Omega;\nu)$, therefore it admits a unique Lipschitz extension to the whole $L^1(\Omega;\nu)$. We shall call such extended elements, according to \cite[Sec. 6.1]{Vaz07}, \emph{limit solutions}.
\begin{proposition}\label{pro: solLim}
Let $u$ and $v$ be two limit solutions of \eqref{eq: pmeDir} corresponding to two initial data $u_0,v_0 \in L^1(\Omega;\nu)$. Then:
\begin{itemize}
\item[$\bullet$] if in addition $u_0 \in L^{m+1}(\Omega;\nu)$, $u$ is the energy solution;
\item[$\bullet$] for a.e.\ $\tau>0$, $u|_{[\tau,\infty)}$ is the limit solution corresponding to the initial datum $u(\tau)$ (semigroup property);
\item[$\bullet$] for a.e.\ $T>0$ inequality \eqref{eq: princConf} and in particular the comparison principle hold.
\end{itemize}
\begin{proof}
The first claim is obvious. The semigroup property and inequality \eqref{eq: princConf} follow by approximation from the corresponding properties valid for energy solutions.
\end{proof}
\end{proposition}
In Section \ref{sec: regAs} we shall see that under the sole hypothesis that the Poincar\'e inequality \eqref{eq: poinIntroMed} holds, evolution \eqref{eq: pmeDir} gives rise to an $L^{q_0}$-$L^\varrho$ regularizing effect for any $q_0 \in [1,\infty)$ and $\varrho \in (q_0,\infty)$: this, together with Proposition \ref{pro: solLim}, implies in particular that limit solutions are indeed weak energy solutions after an arbitrarily small time $\tau>0$.
\begin{remark}\label{remcontinui}\rm
We note that since all weak energy solutions belong to the space $C([\tau,\infty);L^{m+1}(\Omega;\nu))$, actually the statements and proofs of Lemma \ref{lem: lemFondDir}, Theorem \ref{teo: teoFondDir} and Corollary \ref{cor: contZero} hold ``for any $T>0$'' rather than only ``for a.e.\ $T>0$''.
\end{remark}
\noindent\bf Comparison to some previous results. \rm
In the particular context where $\Omega=\mathbb{R}^N$ ($N \geq 3$), $\rho_\mu\equiv1$ and $\rho_\nu$ is a weight which satisfies appropriate decay conditions as $|\mathbf{x}|\rightarrow \infty$, recent works provided some existence and uniqueness results for nonnegative solutions of the {WPME} (also called {Inhomogeneous PME} -- see e.g.\ \cite{RV08} and quoted references). Let us briefly compare such results to ours.

Given a nonnegative initial datum $u_0 \in L^{1}(\mathbb{R}^N;\nu)$, according to \cite[Def. 1.1]{RV08} any nonnegative function $u(\mathbf{x},t)$ is a weak solution of \eqref{eq: pmeDir} if it is continuous in $\mathbb{R}^N \times (0,\infty)$ and:
\begin{itemize}
\item[$\bullet$] $u \in C([0,\infty);L^1(\mathbb{R}^N;\nu)) \cap L^\infty(\mathbb{R}^N \times (\tau , \infty)) \ \ \ \forall \tau>0$;
\item[$\bullet$] $\nabla(u^m) \in \left[ L^2(\mathbb{R}^N \times (\tau, \infty)) \right]^N \ \ \ \forall \tau >0$;
\item[$\bullet$] for any $\varphi \in C^1_c(\mathbb{R}^N \times (0,\infty))$ the identity
\begin{equation}\label{eq: propDebRV}
\int_0^\infty \!\! \int_{\mathbb{R}^N} \left(  \nabla(u^m)(\mathbf{x},t) \cdot \nabla{\varphi}(\mathbf{x},t)- u(\mathbf{x},t) \varphi_t(\mathbf{x},t) \rho_\nu(\mathbf{x}) \right) \, \mathrm{d}\mathbf{x} \, \mathrm{d}{t} =0
\end{equation}
holds true;
\item[$\bullet$] $u(\cdot,0)=u_0(\cdot)$.
\end{itemize}
The most important difference between our definition of weak solution and the one just given lies in the space where $u$ is looked for. In fact note that in \cite[Def. 1.1]{RV08} $u^m$ is not related to the test function space chosen in \eqref{eq: propDebRV}: in other words, it is not imposed that $u^m(\cdot,t)$ belongs to the closure of $C^1_c(\mathbb{R}^N)$ with respect to a suitable norm. Indeed when $\rho_\nu(\mathbf{x})$ goes to zero sufficiently fast as $\mathbf{|x|}\rightarrow \infty$ some non-uniqueness issues arise (see, e.g., \cite{Eid90}).

The two most important well-posedness results proved in \cite{RV08} are the following:
\begin{itemize}
\item[$\bullet$] if $\rho_\nu \in C^1(\mathbb{R}^N)$ is bounded and strictly positive,  then \cite[Th. 3.1]{RV08} there exists a weak solution according to \cite[Def. 1.1]{RV08};
\item[$\bullet$] if, in addition, $\rho_\nu$ satisfies
\begin{equation*}
A_0 (1+|\mathbf{x}|)^{-N}  \leq  \rho_\nu(\mathbf{x}) \ \ \ \forall \mathbf{x} \in \mathbb{R}^N
\end{equation*}
for a suitable constant $A_0>0$, then such solution is also unique \cite[Th. 4.1]{RV08}.
\end{itemize}
The uniqueness result, in some sense, is not improvable: if $\rho_\nu(\mathbf{x})$ behaves like $|\mathbf{x}|^{-\gamma}$ at infinity, with $\gamma>N$, the finiteness of the $\nu$-measure of $\mathbb{R}^N$ implies that if the initial datum is $u_0 \equiv 1$ then $u(\mathbf{x},t) \equiv 1$ is a solution of \eqref{eq: pmeDir} according to \cite[Def. 1.1]{RV08}. Yet it is possible to prove \cite[Sec. 8]{RV08} that in this case (even for any $\gamma>2$) the solution built up in \cite[Th. 3.1]{RV08} necessarily satisfies the decay condition
\begin{equation}\label{eq: decInfRV}
\lim_{R \rightarrow \infty} R^{1-N} \int_{|\boldsymbol{\sigma}|=R} \int_0^T u^m(\boldsymbol{\sigma},t) \, \mathrm{d}t \, \mathrm{d}\boldsymbol{\sigma} =0 \ \ \ \forall T>0 \, ;
\end{equation}
since \eqref{eq: decInfRV} is trivially \emph{not} fulfilled by nonzero constants, this means we have at least two solutions. When the initial datum belongs to $L^{1}(\mathbb{R}^N;\nu) \cap L^{\infty}(\mathbb{R}^N)$ actually the solution from \cite[Th. 3.1]{RV08} satisfies the requirements of \cite[Def. 1.1]{RV08} down to $\tau=0$ and it is indeed an energy solution. Therefore it seems natural to wonder how such non-uniqueness problem matches with the uniqueness of energy solutions proved by Proposition \ref{pro: uniqPmeDir}. The answer is that in this case nonzero constants \emph{do not} belong to $V^{(m+1)/m}_0(\mathbb{R}^N;\nu,1)$. In fact whenever $\rho_\nu(\mathbf{x})$ behaves like $|\mathbf{x}|^{-\alpha}$ as $|\mathbf{x}|\rightarrow \infty$, with $\alpha>N$, in $W^{1,2}_0(\mathbb{R}^N;\nu,1)$ the Poincar\'e inequality holds (Section \ref{sec: exa0}), thus preventing any nonzero constant to lie in such space (as well as in $V^{(m+1)/m}_0(\mathbb{R}^N;\nu,1)$). Roughly speaking, the choice of test functions in the weak formulation \eqref{eq: propDebRV} corresponds to the one typical of a Dirichlet problem; however, no ``boundary condition'' is specified on $u$. Consequently, when the weight $\rho_\nu(\mathbf{x})$ goes to zero sufficiently fast as $|\mathbf{x}|\rightarrow \infty$, $\mathbb{R}^N$ behaves like a bounded domain, so that one expects to have to put boundary conditions at infinity to guarantee uniqueness. Indeed, condition \eqref{eq: decInfRV} turns out to be \emph{sufficient} for uniqueness (see \cite[Sec. 8]{RV08} and the references quoted therein).

Despite these non-uniqueness issues (for $\alpha>N$), it is not difficult to verify that the weak solution of \eqref{eq: pmeDir} (according to \cite[Def. 1.1]{RV08}) constructed in \cite[Th. 3.1]{RV08} coincides with the weak energy solution of the same problem (according to Definition \ref{den: solDebDir}) whose existence was proved by Theorem \ref{teo: teoFondDir}, at least for $u_0 \in L^1(\mathbb{R}^N;\nu) \cap L^\infty(\mathbb{R}^N)$ (and so for any $L^1$ datum thanks to inequality \eqref{eq: princConf}). As a consequence, in this context solutions are in fact $C([0,\infty);L^1(\mathbb{R}^N;\nu))$. Still from the results of \cite{RV08} we know that an $L^\infty$ regularizing effect takes place. This is consistent with the validity of the Sobolev inequality in $W^{1,2}_0(\mathbb{R}^N;\nu,1)$ (see \cite[Th. 1.5]{BG05m}).
\end{subsection}

\begin{subsection}{The {WPME} with Neumann boundary conditions}\label{sec: pmeNeu}
As we did for the Dirichlet problem, first of all we give a definition of weak solution for \eqref{eq: pmeNeu}. Again, it is important to comment that the present results hold if $\Omega={\mathbb R}^N$ as well. The only difference with respect to the Dirichlet problem are the weighted Sobolev spaces involved.
\begin{definition}\label{den: solDebNeu}
A function
$$u \in L^2((0,T);L^2(\Omega;\nu)): \ \nabla{(u^m)} \in L^2((0,T);[L^2(\Omega;\mu)]^N) \ \ \ \forall T>0  $$
is a weak solution of \eqref{eq: pmeNeu} with initial datum $u_0 \in L^2(\Omega;\nu)$ if it satisfies
\begin{equation}\label{eq: solDebNeu}
\begin{aligned}
& \int_0^T \! \! \int_{\Omega} u(\mathbf{x},t) \eta_t(\mathbf{x},t) \, \mathrm{d}\nu \, \mathrm{d}t \\
= & -\int_\Omega u_0(\mathbf{x}) \eta(\mathbf{x},0) \,  \mathrm{d}\nu + \int_0^T \! \! \int_{\Omega} \nabla{\left( u^m\right)}(\mathbf{x},t) \cdot \nabla{\eta}(\mathbf{x},t) \, \mathrm{d}\mu \, \mathrm{d}t
\end{aligned}
\end{equation}
$$ \forall \eta \in W^{1,2}((0,T);L^2(\Omega;\nu)): \  \nabla{\eta} \in L^2((0,T);[L^2(\Omega;\mu)]^N)  \, , \ \eta(T)=0 \, . $$
\end{definition}
Note that the weak formulation \eqref{eq: solDebNeu} is very similar to \eqref{eq: solDebDir}; in fact the boundary condition in \eqref{eq: pmeNeu} is purely formal, and what really changes with respect to \eqref{eq: pmeDir} are the underlying functional spaces.

Let us start now some well-posedness analysis. Most of the proofs of this section are driven from analogous ones already performed in the non-weighted context (see \cite[Sec. 11.2]{Vaz07}).
\begin{proposition}\label{pro: uniqNeu}
There exists \emph{at most} one weak solution of problem \eqref{eq: pmeNeu} which satisfies the following additional hypothesis:
\begin{equation}\label{eq: hpEnergNeu}
u \in L^{m+1}((0,T); L^{m+1}(\Omega;\nu)) \, .
\end{equation}
\begin{proof}
There is no major difference with respect to the proof of Proposition \ref{pro: uniqPmeDir}: one plugs Ole\u{\i}nik's test function into the weak formulation solved by the difference of two hypothetical solutions $u_1-u_2$ satisfying \eqref{eq: hpEnergNeu} and then argues likewise. The only relevant issue is to prove that such a test function is admissible, which is easily achievable by approximating it with a sequence of test functions as in the weak formulation \eqref{eq: solDebNeu}
\end{proof}
\end{proposition}
\noindent Again, we can give the  definition of (weak) \emph{energy solutions} as follows.
\begin{definition}\label{wes-neumann}
We shall call (weak) \emph{energy solutions} all weak solutions to \eqref{eq: pmeNeu} that also satisfy \eqref{eq: hpEnergNeu}.
\end{definition}
\noindent Existence of such solutions is ensured by the next theorem.
\begin{theorem}\label{teo: teoFondNeu}
Let $\Omega \subset \mathbb{R}^N$ be a domain, and let $\rho_\nu$ and $\rho_\mu$ be two weights such that
\begin{equation*}
\rho_\nu \in C^{3,\alpha}_{loc}(\Omega) \, , \ \rho_\mu \in C^{2,\alpha}_{loc}(\Omega) \, , \ \rho_\nu^{-1},\rho_\mu^{-1} \in L^{\infty}_{loc}(\Omega) \, .
\end{equation*}
If $u_0 \in L^1(\Omega;\nu) \cap L^{r}(\Omega;\nu) $, with $r \geq m+1$, then there exists a unique weak energy solution $u$ of \eqref{eq: pmeNeu} in the sense of Definition \ref{wes-neumann}, which satisfies estimates \eqref{eq: stimeLemmaDir1} and \eqref{eq: stimeLemmaDir2} for all $q \leq r-1$ and almost every $T>0$. If in addition $\nabla{(u_0^m)} \in [L^2(\Omega;\mu)]^N$, then  also estimate  \eqref{eq: stimaRegIniziale} holds true. Moreover, if $v$ is another energy solution corresponding to an initial datum $v_0 \in L^{1}(\Omega;\nu) \cap L^{m+1}(\Omega;\nu)$, inequality \eqref{eq: princConf} and in particular the comparison principle hold.
\begin{proof}
We proceed similarly to the proofs of Lemma \ref{lem: lemFondDir}, Theorem \ref{teo: teoFondDir} and Corollary \ref{cor: contZero}. That is, given a sequence of smooth functions $\{\Phi_n^\prime(x)\}$ approximating $m \, |x|^{m-1}$ as in Lemma \ref{lem: lemFondDir} and a \emph{fixed}  smooth domain $\Omega^\prime  \Subset \Omega $, one solves the following Neumann problems:
\begin{equation} \label{eq: pmeNeuParab}
\begin{cases}
(u_n)_t =\rho_{\nu}^{-1}  \operatorname{div}\left(\rho_{\mu} \,  \nabla{\left( \Phi_n(u_n) \right)} \right) & \textnormal{in} \  \Omega^\prime \times(0,\infty)  \\
 \frac{\partial{\Phi_n(u_n)}}{\partial\mathbf{n}} =0 & \textnormal{on} \ \partial\Omega^\prime \times (0,\infty) \\
 u_n(\cdot,0)=u_0(\cdot) & \textnormal{in} \ \Omega^\prime
 \end{cases}\, ,
\end{equation}
assuming in addition $u_0 \in C^{2,\alpha}\left(\Omega^\prime \right)$ and $\frac{\partial{(u_0)}}{\partial\mathbf{n}} =0$ on $\partial \Omega^\prime$. Setting $w=\rho_\nu u_n$, let us rewrite \eqref{eq: pmeNeuParab} in divergence form:
\begin{equation} \label{eq: pmeNeuParabDiv}
\begin{cases}
  w_t = \operatorname{div}\left(\frac{\rho_\mu}{\rho_\nu} \, \Phi_n^\prime\left(\frac{w}{\rho_\nu}\right) \nabla{w} - \frac{\rho_\mu}{\rho_\nu^2} \, \nabla{(\rho_\nu)} \, \Phi_n^\prime\left(\frac{w}{\rho_\nu}\right)w \right) & \textnormal{in} \  \Omega^\prime\times(0,\infty)  \\
 \Phi_n^\prime\left(\frac{w}{\rho_\nu}\right) \left( \nabla{w} - \frac{w}{\rho_\nu} \, \nabla{(\rho_\nu)}  \right) \cdot \mathbf{n} =0 & \textnormal{on} \ \partial\Omega^\prime \times (0,\infty) \\
 w(\cdot,0)=\rho_\nu(\cdot) u_0(\cdot) & \textnormal{in} \ \Omega^\prime
 \end{cases} \, .
\end{equation}
Quasilinear theory (see, for instance, \cite[Th. 13.24]{Lie96}) ensures that problem \eqref{eq: pmeNeuParabDiv} (and so \eqref{eq: pmeNeuParab}) admits a regular solution $w$ ($u_n$). From such solutions, proceeding as in the proof of Lemma \ref{lem: lemFondDir} and in the beginning of the proof of Theorem \ref{teo: teoFondDir}, one gets in turn a solution $u$ of \eqref{eq: solDebNeu} satisfying \eqref{eq: stimeLemmaDir1}, \eqref{eq: stimeLemmaDir2} (and \eqref{eq: stimaRegIniziale} when $\nabla{(u_0^m)} \in [L^2(\Omega;\mu)]^N$), at least if $u_0 \in L^\infty(\Omega^\prime)$ (also \eqref{eq: princConf} still holds). The crucial point is to extend this result to general domains: as in Theorem \ref{teo: teoFondDir}, one picks an initial datum $u_0 \in L^1(\Omega;\nu) \cap L^{\infty}(\Omega)$, an increasing sequence of domains $\Omega_n \in  C^{2,\alpha} $ approximating $\Omega$, with $\Omega_n \Subset \Omega$, solves in them the homogeneous Neumann problems \eqref{eq: pmeNeu} with initial data $u_{0n}=u_0|_{\Omega_n}$, denotes as $\{u_n\}$ the corresponding sequence of solutions and exploits analogous estimates. Now, if $u_n^m$ were extended to be zero outside $\Omega_n$ (what we actually do in Theorem \ref{teo: teoFondDir}), in general it would not belong to $W^{1,2}(\Omega;\nu,\mu)$. However, this does not matter: it suffices to extend to zero $u_n$ (so $u_n^m$) and $\nabla{(u_n^m)}$ independently from each other. That is, setting
$$ z_n=u_n \, \chi_{\Omega_n^{\phantom{A}}} , \ \mathbf{w}_n= \nabla{(u_n^m)} \, \chi_{\Omega_n^{\phantom{A}}} , $$
one has that $\{z_n\}$ and $\{\mathbf{w}_n\}$, up to subsequences, converge respectively pointwise and weakly in $L^2((0,T);$ $L^2(\Omega;\nu)) $ to $u$ and weakly in $L^2((0,T);$ $[L^2(\Omega;\mu)]^N)$ to $\mathbf{w}$. This is enough in order to pass to the limit in the weak formulation \eqref{eq: solDebNeu}. Thus it remains to show that $\nabla{(u^m)}=\mathbf{w}$. First of all observe that, given any $\Omega^\prime \Subset \Omega$, $\mathbf{w}_n|_{\Omega^\prime}=\nabla{(u_n^m|_{\Omega^\prime})}$ for any $n$ large enough; since $u_n^m|_{\Omega^\prime} \rightharpoonup u^m|_{\Omega^\prime} $ in $L^2((0,T); H^1(\Omega^\prime))$ and also $\mathbf{w}_n|_{\Omega^\prime} \rightharpoonup \mathbf{w}|_{\Omega^\prime}$ in $L^2((0,T);L^2(\Omega^\prime))$, necessarily $\nabla{(u^m|_{\Omega^\prime})}=\mathbf{w}|_{\Omega^\prime}$. The assertion follows from the arbitrariness of $\Omega^\prime$.

The validity of inequalities \eqref{eq: stimeLemmaDir1}--\eqref{eq: princConf} (and \eqref{eq: stimaRegIniziale} when $\nabla{(u_0^m)} \in [L^2(\Omega;\mu)]^N$), and their generalization to initial data in $L^1(\Omega;\nu)\cap L^{r}(\Omega;\nu)$ (for $r\geq m+1$), can be shown exactly as for the Dirichlet case.
\end{proof}
\end{theorem}
The theorem just proved provides us with the unique energy solution of \eqref{eq: pmeNeu}. Limit solutions for general $L^1$ data are defined exactly as in Section \ref{sec: pmeDir}. Let us also observe that, as a consequence of the method of proof of Theorems \ref{teo: teoFondDir} and \ref{teo: teoFondNeu}, we obtain the classical conservation of positivity: if $u_0$ is (essentially) nonnegative then $u(\cdot,t)$ is (essentially) nonnegative for a.e.\ $t>0$. Indeed such property is inherited from the solutions of the non-degenerate approximating problems. Actually, for nonnegative data, one can even set up a different and more ``natural'' proof of existence (see \cite[Sec. 5.4]{Vaz07}).

Of particular interest for our purposes is the case $\nu(\Omega)<\infty$, where it makes sense to speak about the weighted mean value \eqref{eq: meanV}. The next result is straightforward as well as classical and of great importance.
\begin{proposition}\label{pro: consMed}
Let $\nu(\Omega)<\infty$. If $u$ is a weak energy solution of \eqref{eq: pmeNeu} then
\begin{equation}\label{eq: consMed}
 \overline{u(t)}= \overline{u_0}=\overline{u} \ \ \ \textnormal{for a.e. } t>0 \, .
\end{equation}
\begin{proof}
Thanks to the hypotheses, one is allowed to plug into \eqref{eq: solDebNeu} the following test functions (independent of $\mathbf{x}$):
\begin{equation*}
\eta_h(s)=\chi_{[0,t-h/2)}(s) + \chi_{[t-h/2,t+h/2]}(s) \left( \frac{t-s}{h}+\frac{1}{2}  \right) \, ,
\end{equation*}
(where $0 < h < 2t$ is arbitrarily fixed)
which gives
\begin{equation*}
\frac{1}{h} \int_{t-h/2}^{t+h/2} \int_\Omega u(\mathbf{x},s) \, \mathrm{d}\nu \, \mathrm{d}s = \int_\Omega u_0(\mathbf{x}) \, \mathrm{d}\nu \, .
\end{equation*}
The assertion follows by letting $h\rightarrow 0$ and using Lebesgue's differentiation Theorem.
\end{proof}
\end{proposition}
Note that, since $\nu(\Omega)<\infty$,  any weak energy solution always belongs to the space $C((0,\infty);L^{1}(\Omega;\nu))$,  so that the ``a.e.''\ in \eqref{eq: consMed} can actually be removed (recall the brief discussion about continuity in Remark \ref{remcontinui}, which applies in this context too).

In Section \ref{sec: pmeNeuReg} we shall prove that the validity of the zero-mean Poincar\'e inequality \eqref{eq: poinIntroMed} implies an $L^{q_0}$-$L^{\varrho}$ regularization (provided $q_0 \in [1,\infty)$ and $\varrho \in (q_0,\infty)$) also for evolution \eqref{eq: pmeNeu}, which again means that, in particular, limit solutions are energy solutions after an arbitrarily small time $\tau>0$.
\end{subsection}
\end{section}

\begin{section}{Smoothing and asymptotic estimates: the Dirichlet problem} \label{sec: regAs}
In this and in the next section we investigate connections between the validity of weighted Poincar\'e inequalities and integrability properties of the solutions to \eqref{eq: pmeDir} and \eqref{eq: pmeNeu}, both for short and large times. In the sequel, we shall implicitly assume that the weights satisfy all the hypotheses of Theorem \ref{teo: teoFondDir} for well-posedness. Moreover, when referring to ``the solution" to the equation at hand we shall always mean, without further comment, the one constructed in Section \ref{sec: wellP} (see Theorem \ref{teo: teoFondDir} and the subsequent discussion about limit solutions for general $L^1$ data): this is particularly relevant in view of possible non-uniqueness issues which may arise. As already mentioned, the incoming results also hold, with no modifications, if $\Omega={\mathbb R}^N$.

Suppose that in the domain $\Omega\subset \mathbb{R}^N$, with respect to the weights $\rho_\nu,\rho_\mu$, the Poincar\'e inequality \eqref{eq: poinIntro} holds. By means of a Gross differential method we shall prove that solutions to the Dirichlet problem \eqref{eq: pmeDir} enjoy an $L^{q_0}$-$L^\varrho$ regularization for all $q_0 \in [1,\infty)$ and $\varrho \in (q_0,\infty)$. In fact regularization into $L^\infty$ need not hold in general, as we now show.
\\[2mm]
\noindent\bf Counterexample to the $\boldsymbol{{L}^\infty(\Omega)}$ regularization. \rm
Let $\Omega=(0,\infty)$. With respect to the weights $\rho_\nu(x)=\rho_\mu(x)=e^{-x}$, it is known that the Poincar\'e inequality \eqref{eq: poinIntro} holds (see Section \ref{list}). In this context, the {WPME} with Dirichlet boundary conditions reads
\begin{equation}\label{eq: pmeControEx}
\begin{cases}
u_t=e^x (e^{-x}(u^m )_x  )_x  & \mathrm{in \ } (0,\infty) \times (0,\infty) \\
u(0,t)=0   & \mathrm{for \ } t>0  \\
u(x,0)=u_0(x)  & \mathrm{in \ } (0,\infty)
\end{cases} \, .
\end{equation}
We want to prove that the solution $u(x,t)$ corresponding to the initial datum $u_0(x)=\log(x+1)$ remains unbounded for all $t \ge 0$. To this aim, consider the following family of functions:
\begin{equation*} 
v_B(x,t)=\frac{\log(x+1)}{\left(1 + B^{-1} (m-1) t \right)^{\frac{1}{m-1}} } \, .
\end{equation*}
We can show that for a suitable choice of the constant $B>0$, $v_B$ is a \emph{subsolution} to \eqref{eq: pmeControEx}. In fact, after some computations, one gets:
\begin{equation*} 
\begin{aligned}
e^x (e^{-x}(  \left[ \log(x+1) \right]^m )_x  )_x  = & -m \frac{[\log(x+1) ]^{m-1}}{x+1} -m \frac{[\log(x+1)]^{m-1} } {(x+1)^2} \\
&+  m(m-1) \frac{[\log(x+1) ]^{m-2} }{(x+1)^2} \, .
\end{aligned}
\end{equation*}
Clearly, there exists a constant $\widehat{B}>0$ such that
\begin{equation*} 
 \log(x+1) \geq -\widehat{B} e^x (e^{-x}(  [ \log(x+1) ]^m )_x  )_x  \,,
\end{equation*}
so that
\begin{equation*} 
\begin{aligned}
 \left(v_{\widehat{B}}\right)_t  = -\frac{\log(x+1)}{\widehat{B}\left(1+{\widehat{B}}^{-1} (m-1) t \right)^{\frac{m}{m-1}} } \leq & \frac{ e^x (e^{-x}(  [ \log(x+1) ]^m )_x  )_x }{\left(1 + {\widehat{B}}^{-1} (m-1) t  \right)^{\frac{m}{m-1}} } \\
 = & e^x \left(e^{-x} \left( \left[ v_{\widehat{B}} \right]^m \right)_x  \right)_x  \, .
\end{aligned}
\end{equation*}
As a consequence, $v_{\widehat{B}}$ is a subsolution to \eqref{eq: pmeControEx} for $u_0(x)=\log(x+1)$. From the comparison principle for sub-supersolutions (given the regularity of the data, one can argue as in \cite[Th. 8.10]{Vaz07} -- see also \cite{KRV10}) we have that $v_{\widehat{B}} \leq u$; in particular, $u(\cdot,t)$ is unbounded for all $t \geq 0$. Moreover, since $u_0 \in L^{q_0}((0,\infty);e^{-x}) \ \forall q_0 \in [1,\infty) $, we have also shown that in this case the $L^{q_0}$-$L^\infty$ regularization \emph{does not} take place for any $q_0 \in [1,\infty)$.

\begin{theorem}\label{teo: regPmeDir}
Let $q_0 \in [1,\infty)$ and $u_0 \in L^1(\Omega;\nu) \cap L^{q_0}(\Omega;\nu)$. If the Poincar\'e inequality \eqref{eq: poinIntro} holds then the solution $u$ of \eqref{eq: pmeDir}  with initial datum $u_0$ satisfies the estimate
\begin{equation}\label{eq: stimeRegPmeDir}
\left\| u(t) \right\|_{\varrho;\nu} \leq K_1 \, t^{-\frac{\varrho-q_0}{\varrho(m-1)}} \left\| u_0 \right\|_{q_0;\nu^{\phantom{a}}}^{\frac{q_0}{\varrho}} \ \ \ \hbox{\rm for a.e.}\,\, t>0 \, ,
\end{equation}
where $\varrho \in (q_0,\infty)$ and $K_1>0$ is a suitable constant depending only on $\varrho$, $m$ and $C_P$. \\
Moreover, if $\nu(\Omega) < \infty $ the absolute bound
\begin{equation}\label{eq: pmeDirAbsB}
\left\| u(t)  \right\|_{\varrho;\nu} \leq K_2 \, t^{-\frac{1}{m-1}} \ \ \ \hbox{\rm for a.e.}\,\, t>0
\end{equation}
holds as well, $K_2>0$ being another constant depending only on $\varrho$, $m$, $C_P$ and $\nu(\Omega)$.
\begin{proof}
The procedure is very similar to the one developed in \cite[Th. 1.3]{Gri10}, 
so we just point out the most significant differences. Upon defining the \emph{entropy functional}
\begin{equation*} 
J(r,f)=\int_{\Omega} \frac{\left| f \right|^r}{\left\| f \right\|_{r;\nu}^r} \log{ \left( \frac{\left|f \right|}{\left\| f \right\|_{r;\nu }} \right) } \,  \mathrm{d}\nu \, ,
\end{equation*}
the validity of the family of logarithmic Sobolev inequalities
\begin{gather}
\left(J(r,v) + \frac{1}{2-r} \log \varepsilon \right) \frac{(2-r)\left\| v \right\|_{r;\nu}^2}{\varepsilon \,  C_P^2 } \leq \left\| \nabla v \right\|_{2;\mu}^2 \label{eq: 2LogSobDir} \\
\forall r \in [1,2) \, , \ \forall \varepsilon >0 \, , \ \forall v \in L^r(\Omega;\nu) \cap W_0^{1,2}(\Omega;\nu,\mu) \nonumber
\end{gather}
was already established in \cite[Th. 1.3]{Gri10}. Given $u_0 \in L^1(\Omega;\nu) \cap L^\infty(\Omega)$, $t>0$, $q_0 \in (1,\infty)$, $\varrho \in (q_0,\infty)$, let $u$ be the solution of \eqref{eq: pmeDir} with initial datum $u_0$. Introducing an increasing, one-to-one and $C^1[0,t]$ function $q:[0,t] \rightarrow [q_0,\varrho] $, after explicit calculations one gets:
\begin{equation*}
\begin{aligned}
\frac{\mathrm{d}}{\mathrm{d}s} \log{\left\| u(s) \right\|_{q(s);\nu}}  = &  \frac{{q^\prime}(s)}{q(s)} J(q(s),u(s)) \\ & -\left( \frac{2}{q(s)+m-1} \right)^2 \frac{m(q(s)-1)}{\left\| u(s) \right\|_{q(s);\nu}^{q(s)}}  \left\| \nabla{\left( u^{\frac{q(s)+m-1}{2}} \right) }(s) \right\|_{2;\mu}^2 \, .
\end{aligned}
\end{equation*}
From now on, one applies \eqref{eq: 2LogSobDir} to $u^{(q+m-1)/2}$ in the equation above, sets $q(s)=q_0+\frac{s}{t}(\varrho-q_0)$, chooses $r$ and $\varepsilon$ appropriately and solves the resulting differential inequality in the variable $y(s)=\log{\left\| u(s) \right\|_{q(s);\nu}}$ along the lines of the proof given in \cite[Th. 1.3]{Gri10} (with respect to the notation used therein, it is enough to substitute $q-1$ with $m(q-1)$, $p-2$ with $m-1$, $p$ with $2$ and $C$ with $C_P^2$). So we get estimate \eqref{eq: stimeRegPmeDir}  for $q_0>1$. The case $q_0=1$ is obtained by taking limits since the constant $K_1$ in \eqref{eq: stimeRegPmeDir} can be shown to be bounded as $q_0\downarrow1$.

Concerning the absolute bound, first of all note that, thanks to \eqref{eq: poinIntro} and to the finiteness of $\nu(\Omega)$, one has:
\begin{equation*}
\begin{aligned}
\frac{\mathrm{d}}{\mathrm{d}s} \left\| u(s) \right\|_{\varrho;\nu}^\varrho & = - \left( \frac{2}{\varrho+m-1} \right)^2 m  \varrho ( \varrho-1 ) \left\| \nabla \left( u^{\frac{\varrho+m-1}{2}} \right) (s) \right\|_{2;\mu}^2 \\
 & \leq -\left( \frac{2}{\varrho+m-1} \right)^2 \frac{ m  \varrho ( \varrho-1 ) }{C_P^2} \left\| u(s) \right\|_{\varrho+m-1;\nu}^{\varrho+m-1}  \\
& \leq -\left( \frac{2}{\varrho+m-1} \right)^2 \frac{ m  \varrho ( \varrho-1 ) }{C_P^2 \, \nu(\Omega)^{\frac{m-1}{\varrho}} }  \left\| u(s) \right\|_{\varrho;\nu}^{\varrho+m-1}\\ &= -D \left(\left\| u(s) \right\|_{\varrho;\nu}^{\varrho} \right)^{\frac{\varrho+m-1}{\varrho} }  \, ,
\end{aligned}
\end{equation*}
where $D>0$ is a constant depending only on $\varrho$, $m$, $C_P$ and $\nu(\Omega)$. Solving the above differential inequality in the variable $y(s)=\|u(s) \|_{\varrho;\nu}^{\varrho} $ one arrives at:
\begin{equation*}
\left\| u(t) \right\|_{\varrho;\nu} \leq \frac{1}{\left( \left\| u_0 \right\|_{\varrho;\nu}^{1-m} + \frac{D(m-1)}{\varrho} t \right)^{\frac{1}{m-1}}} \ \ \ \forall t>0 \, ,
\end{equation*}
from which \eqref{eq: pmeDirAbsB} follows immediately.

Finally, the removal of the hypothesis $u_0 \in L^\infty(\Omega)$ is standard (see the end of the proof of \cite[Th. 1.3]{Gri10}). 
\end{proof}
\end{theorem}
Notice that the conclusions of Theorem \ref{teo: regPmeDir} hold for \emph{any} $t>0$ when weak energy solutions are considered, see Remark \ref{remcontinui}.
\begin{remark}
The calculations performed in the proof just given are formal, since the solution might not be regular enough. Nonetheless they can be justified by approximation, for instance developing analogous ones for the non-degenerate problems solved in Lemma \ref{lem: lemFondDir} (through a similar differential method) and passing to the limit.
\end{remark}
We also have, in some sense, the converse of Theorem \ref{teo: regPmeDir}.
\begin{theorem}\label{teo: pmeDirImpInv}
Suppose $\nu(\Omega)<\infty$. If there exist a constant $K_1>0$ and two given numbers $q_0 \in [1,m+1)$ and $\varrho\geq m+1$ such that, for all $u_0\in L^{q_0}(\Omega; \nu)$, the solution $u$ of \eqref{eq: pmeDir} corresponding to the initial datum $u_0$ satisfies the estimate
\begin{equation}\label{eq: pmeRegInv}
\left\| u(t) \right\|_{\varrho;\nu} \leq K_1 \, t^{-\frac{\varrho-q_0}{\varrho(m-1)}} \left\| u_0  \right\|_{q_0;\nu^{\phantom{a}}}^{\frac{q_0}{\varrho}} \ \ \ \hbox{\rm for a.e.}\,\, t>0 \, ,
\end{equation}
then the Poincar\'e inequality \eqref{eq: poinIntro} holds. In particular, the validity of \eqref{eq: pmeRegInv} for two given $\bar{q}_0 \in [1,m+1)$ and $\bar{\varrho}\ge m+1$ is equivalent to the validity of \eqref{eq: poinIntro}, and hence it implies the validity of \eqref{eq: pmeRegInv} for \emph{any} $q_0 \in [1,\infty)$ and $\varrho \in (q_0,\infty)$.
\begin{proof}
As in the second part of the proof of \cite[Th. 1.3]{Gri10}, we want to take advantage of the strong result \cite[Th. 3.1]{BC+95}. In order to do that, let us consider an initial datum $u_0 \in W^{1,\infty}_c(\Omega)$. First of all, one can prove the following inequality: 
\begin{equation}\label{eq: impInvPmeDir}
\left\| u(t) \right\|_{m+1;\nu}^{m+1} - \left\| u_0 \right\|_{m+1;\nu}^{m+1} \geq - (m+1) \, t \left\| \nabla{(u_0^m)}  \right\|_{2;\mu}^2  \ \ \ \hbox{\rm for almost every} \,\, t>0 \, .
\end{equation}
Formally, \eqref{eq: impInvPmeDir} is easily obtainable by multiplying equation \eqref{eq: pmeDir} by $\rho_\nu u^m$, integrating in $\Omega \times (0,t)$ and exploiting the fact that the quantity $\| \nabla{(u^m)} (\cdot) \|_{2;\mu}$ is nonincreasing (see Corollary \ref{cor: contZero}). However, in this case we must proceed more carefully. If estimate \eqref{eq: stimeLemmaDir1} (for $q=m$) were an \emph{equality} then \eqref{eq: impInvPmeDir} would be easily provable in a rigorous way. On the other hand recall that \eqref{eq: stimeLemmaDir1} was deduced by weak convergence, so in general it is just an inequality, with the wrong verse with respect to what we want to show. Yet if $\nu(\Omega)<\infty$ and the initial datum belongs to $W^{1,\infty}_c(\Omega)$, inequality \eqref{eq: impInvPmeDir} holds indeed. To prove it, we need to go back to the approximate problems of Lemma \ref{lem: lemFondDir}: from \eqref{eq: stimaEnerg1} and proceeding as in the proof of \eqref{eq: stimaDerivataT} we infer, in particular, that
\begin{equation*}
\int_{\Omega} \Psi_n(u_n(\mathbf{x},t)) \, \mathrm{d}\nu -\int_{\Omega} \Psi_n(u_0(\mathbf{x})) \, \mathrm{d}\nu \geq - t \, \int_{\Omega} \left| \nabla{(\Phi_n(u_0))}(\mathbf{x}) \right|^2 \, \mathrm{d}\mu  \ \ \ \forall t>0 \, .
\end{equation*}
Passing to the limit as $n \rightarrow \infty$ this last inequality continues to hold for a.e.\ $t > 0$ (provided the initial datum is regular enough) thanks to the pointwise a.e.\ convergence of $\Psi_n(u_n(\cdot,t))$, $\Psi_n(u_0(\cdot))$, $\Phi_n^\prime(u_0(\cdot))$ respectively to $\frac{1}{m+1}|u(\cdot,t)|^{m+1}$, $\frac{1}{m+1}|u_0(\cdot)|^{m+1}$, $m|u_0(\cdot)|^{m-1}$ and the fact that they are dominated in $L^\infty(\Omega)$. The hypotheses of Lemma \ref{lem: lemFondDir} can then be removed as in the proof of Theorem \ref{teo: teoFondDir}.  \\
Now, using \eqref{eq: impInvPmeDir} together with \eqref{eq: pmeRegInv} for $\varrho=m+1$ (by means of an interpolation between the norms $\| \cdot \|_{q_0;\nu}$, $\| \cdot \|_{m+1;\nu}$ and $\| \cdot \|_{\varrho;\nu}$ on the left hand side and by exploiting the non-expansivity of the $\| \cdot \|_{q_0;\nu}$ norm, one can deduce from \eqref{eq: pmeRegInv} the validity of the same estimate also for $\varrho=m+1$, possibly with a different constant $K_1$), we have:
\begin{equation*}
\left\| u_0 \right\|_{m+1;\nu}^{m+1} \leq  {K_1}^{(m+1)} \, t^{-\frac{(m+1-q_0)}{m-1}} \left\| u_0 \right\|_{q_0;\nu}^{q_0} +(m+1) \, t \left\| \nabla{(u_0^m)}  \right\|_{2;\mu}^2  \ \ \ \forall t>0 \, .
\end{equation*}
Minimizing explicitly (w.r.t.\ $t>0$) the right hand side of the inequality above, we obtain:
\begin{equation}\label{eq: ImpInvPmeDirPre3}
\left\| u_0 \right\|_{m+1;\nu} \leq  B \, \left\| \nabla{(u_0^m)} \right\|_{2;\mu}^{\frac{2(m+1-q_0)}{(m+1)(2m-q_0)}}    \left\| u_0 \right\|_{q_0;\nu}^{\frac{q_0(m-1)}{(m+1)(2m-q_0)}} \, ,
\end{equation}
where $B=B(q_0,m,K_1)>0$ is a suitable constant. In order to rewrite \eqref{eq: ImpInvPmeDirPre3} in a more convenient way for our purposes, we start considering the following sequence $\{\xi_n\}$ of real functions:
\begin{equation*}
\xi_n(x)=2 \left( x-\frac{1}{2n} \right)   \chi_{\left[\frac{1}{2n} , \frac{1}{n}\right)}(x) + 2 \left( x+\frac{1}{2n} \right)   \chi_{\left(-\frac{1}{n} , -\frac{1}{2n}\right]}(x) + x \, \chi_{\left[\frac{1}{n},\infty\right)}\left(|x|\right) \, ,
\end{equation*}
with $x \in {\mathbb R}$.
The regularized approximations
$$ v_n= \xi_n\left(u_0^{\frac{1}{m}} \right)  $$
of $u_0^{1/m}$ still belong to $W_c^{1,\infty}(\Omega)$. Moreover
\begin{equation*}
\left| v_n(\mathbf{x}) \right| \leq \left|u_0(\mathbf{x})\right|^{\frac{1}{m}} \, , \
\nabla{(v_n^m)}(\mathbf{x})= \xi_n^\prime \left( u_0^{\frac{1}{m}}(\mathbf{x}) \right) \left| v_n(\mathbf{x}) \right|^{m-1} \left|u_0(\mathbf{x})\right|^{\frac{1}{m}-1} \nabla{u_0}(\mathbf{x})
\end{equation*}
for a.e. $\mathbf{x} \in \Omega$. In view of the estimates above and from the properties of $\{\xi_n\}$, it is straightforward to check that
$$ \left|\nabla{(v_n^m)}\right|(\mathbf{x}) \leq 2 \left| \nabla{u_0}(\mathbf{x}) \right| \textnormal{\ \ \ for a.e.\ } \mathbf{x} \in \Omega  $$
and that $\{v_n\}$ and $\{\nabla{(v_n^m)}\}$ converge pointwise respectively to $u_0^{1/m}$ and $\nabla{(u_0)}$. Applying then \eqref{eq: ImpInvPmeDirPre3} to the sequence of initial data $\{v_n\}$ and passing to the limit as $n\rightarrow \infty$, by dominated convergence we infer that
\begin{equation}\label{eq: ImpInvPmeDirPre5}
\left\| u_0 \right\|_{\frac{m+1}{m};\nu} \leq  B^m \,  \left\| \nabla{(u_0)} \right\|_{2;\mu}^{\frac{2m(m+1-q_0)}{(m+1)(2m-q_0)}} \left\| u_0 \right\|_{\frac{q_0}{m};\nu}^{\frac{q_0(m-1)}{(m+1)(2m-q_0)}} \, .
\end{equation}
Setting
$$ \vartheta=\frac{2m(m+1-q_0)}{(m+1)(2m-q_0)} \, , \ r=\frac{m+1}{m}  \, , \ s=\frac{q_0}{m} \, , \ q=2 \, , \ \mathcal{W}(f)= \left\| \nabla{f} \right\|_{2;\mu}   \, ,$$
where $f$ is any nonnegative function belonging to $W_c^{1,\infty}(\Omega)$, inequality \eqref{eq: ImpInvPmeDirPre5} reads
\begin{equation}\label{eq: ImpInvPmeDirPreDis}
\left\| f \right\|_{r;\nu} \leq  \left(B^{\frac{m}{\vartheta}} \, \mathcal{W}(f) \right)^{\vartheta} \left\| f \right\|_{s;\nu}^{1-\vartheta} \, , \  \frac{1}{r}=\frac{\vartheta}{q} + \frac{1-\vartheta}{s}  \, ;
\end{equation}
Theorem 3.1 of \cite{BC+95} is now applicable, providing us with the existence of a nonnegative constant (that we keep denoting as $B$) such that \eqref{eq: ImpInvPmeDirPreDis} holds for $\vartheta=1$ and $q=2$ as well, which in this case means
\begin{equation}\label{eq: ImpInvPmeDirPostDis}
\left\| f \right\|_{2;\nu} \leq  B^m \, \left\| \nabla{f} \right\|_{2;\mu} \, ,
\end{equation}
namely the Poincar\'e inequality for nonnegative functions of $W_c^{1,\infty}(\Omega)$. The extension of \eqref{eq: ImpInvPmeDirPostDis} (up to multiplicative constants) to signed functions of $W_c^{1,\infty}(\Omega)$ is simply achieved by writing $f=f_{+}-f_{-}$, while the extension to the whole $W^{1,2}_0(\Omega;\nu,\mu)$ follows by density.
\end{proof}
\end{theorem}
Let us observe that estimate \eqref{eq: stimeRegPmeDir} shows an $L^{q_0}$-$L^{\varrho}$ regularization which \emph{does not} hold up to $\varrho=\infty$ (by direct calculations one verifies that $K_1$ diverges as $\varrho \rightarrow \infty$). If instead we assumed that a Sobolev-type inequality holds, namely that there exists $q>2$ such that $\| v \|_{q;\nu} \leq C_S \| \nabla{v} \|_{2;\mu} $ for all $v \in W^{1,2}_0(\Omega;\nu,\mu)$, then there would be no difficulty in repeating the proof of Theorem 1.5 of \cite{BG05m} and so conclude that an $L^{q_0}$-$L^{\infty}$ regularization takes place indeed. However, the explicit counterexample we constructed above shows that the sole validity of the Poincar\'e inequality in general prevents the $L^\infty$ regularization.

\section{Smoothing and asymptotic estimates: the Neumann problem}\label{sec: pmeNeuReg}
Given a domain $\Omega\subset\mathbb{R}^N$ and two weights $\rho_\nu,\rho_\mu$, assume that $\nu(\Omega)<\infty$. By means of a Gross differential method and Moser iterative techniques, we shall now analyse $L^{q_0}$-$L^{\varrho}$ smoothing and asymptotic properties of solutions to the Neumann problem \eqref{eq: pmeNeu} exploiting the validity of functional inequalities like \eqref{eq: poinIntroMed} or the weaker \eqref{eq: weakPoi}.

Similarly to the Dirichlet problem, the present results hold as well, with no modifications, if $\Omega={\mathbb R}^N$ and its $\nu$-measure is finite.  Notice once again that in the sequel, when referring to ``the solution" to the equation at hand, we shall always mean, without further comment, the one constructed in Section \ref{sec: wellP}, see in particular Theorem \ref{teo: teoFondNeu}.

\subsection{Smoothing estimates}\label{sec: stimeRegNeuPon}
Most of the smoothing results we shall obtain firstly will only hold for initial data which at least belong to $L^{1\vee(m-1)}(\Omega;\nu)$, where  $1\vee(m-1) \equiv \max \{ 1, m-1 \}$. In order to extend them to general $L^{q_0}$ data, the next two lemmas will turn out to be very useful.
\begin{lemma} \label{lem: Dolb}
Suppose that the zero-mean Poincar\'e inequality \eqref{eq: poinIntroMed} holds. Then for any $a \in (0,1]$ one has:
\begin{equation}\label{eq: PoinDol}
\left\| v^a-\overline{v^a}  \right\|_{2;\nu} \leq M_{P,a} \left\| \nabla{v} \right\|^a_{2;\mu} \ \ \ \forall v \in W^{1,2}(\Omega;\nu,\mu) \, ,
\end{equation}
where one can choose $M_{P,a}=2^{1-\frac{a}{2}} \, \nu(\Omega)^{\frac{1}{2}(1-a)} \, M^{a_{}}_P$.
\begin{proof}
See \cite[Prop. 2.2]{DG+08}.
\end{proof}
\end{lemma}
Starting from the previous lemma, we are able to prove a first regularization result.
\begin{lemma} \label{eq: regSottoNeuPoin1}
Suppose that the zero-mean Poincar\'e inequality \eqref{eq: poinIntroMed} holds. Given $k \in \mathbb{N}$, $q_0 \in   \left( 1,({m-1})/({2^k-1}) \right] $ and $m>2$, the solution $u$ of \eqref{eq: pmeNeu} with initial datum $u_0 \in L^{q_0}(\Omega;\nu)$ satisfies the following estimate:
\begin{equation} \label{eq: regSottoStep1}
\left\| u(t) \right\|_{2^k q_0;\nu} \leq D \left( t^{-\frac{1}{q_0+m-1}} \left\|  u_0 \right\|_{q_0;\nu^{\phantom{a}}}^{\frac{q_0}{q_0+m-1}} +  \left\| u_0 \right\|_{q_0;\nu}  \right) \ \ \ \hbox{\rm for a.e.} \,\, t>0 \, ,
\end{equation}
$D>0$ being a constant that depends only on $k$, $q_0$, $m$, $M_P$ and $\nu(\Omega)$.
\begin{proof}
It is convenient to proceed by induction. Let us first prove \eqref{eq: regSottoStep1} for $k=1$. We shall consider $L^\infty$ data, since the passage to general $L^{q_0}$ data is standard. Setting $T=t$ and $q=q_0-1$ in \eqref{eq: stimeLemmaDir1}, applying to the function $u^{(q_0+m-1)/{2}}$ on its left hand side inequality \eqref{eq: PoinDol} and recalling that $\| f \|_{2;\nu}-\nu(\Omega)^{-1/2}\| f \|_{1;\nu} \leq \|f-\overline{f} \|_{2;\nu}$, we obtain:
\begin{equation*} 
\begin{aligned}
& \frac{4(q_0-1)q_0m}{M_{P,a}^{{2}/{a}}(q_0+m-1)^2} \int_{0}^{t} \hskip-4pt\left( \left\| u^{a \frac{q_0+m-1}{2}} (s) \right\|_{2;\nu}  - \frac1{\sqrt{\nu(\Omega)}} \left\| u^{a \frac{q_0+m-1}{2}}(s)  \right\|_{1;\nu}    \right)^{\frac{2}{a}}  \mathrm{d}s  \\
\leq & \left\| u_0 \right\|_{q_0;\nu^{\phantom{a}}}^{q_0} \! \! .
\end{aligned}
\end{equation*}
Exploiting Jensen's inequality in the time integral of the inequality above and raising to the power of ${a}/{2}$, we get ($D$ will always represent a generic constant which possibly depends on $k$, $q_0$, $m$, $M_P$, $\nu(\Omega)$ and may differ from line to line):
\begin{equation*} 
t^{\frac{a}{2}-1} \int_{0}^{t} \left( \left\|  u(s)  \right\|_{a(q_0+m-1);\nu}^{\frac{a}{2}(q_0+m-1)} - \nu(\Omega)^{-\frac{1}{2}} \left\| u(s) \right\|_{\frac{a}{2}(q_0+m-1);\nu}^{\frac{a}{2}(q_0+m-1)}   \right) \, \mathrm{d}s \leq D  \left\| u_0 \right\|_{q_0;\nu^{\phantom{a}}}^{q_0\frac{a}{2}}  \! ;
\end{equation*}
by the non-expansivity of the norms (trivial consequence of \eqref{eq: stimeLemmaDir1}), we then deduce that
\begin{equation*} 
t^{\frac{a}{2}} \left\| u(t)  \right\|_{a(q_0+m-1);\nu}^{\frac{a}{2}(q_0+m-1)}  \leq  D  \left\| u_0 \right\|_{q_0;\nu}^{q_0\frac{a}{2}} + \frac{t^{\frac{a}{2}}}{\nu(\Omega)^{\frac{1}{2}}}  \left\| u_0 \right\|_{\frac{a}{2}(q_0+m-1)}^{\frac{a}{2}(q_0+m-1)}  \, ,
\end{equation*}
that is
\begin{equation} \label{eq: MoserSottoPasso4}
\left\| u(t)  \right\|_{a(q_0+m-1);\nu} \leq D \left( t^{-\frac{1}{q_0+m-1}} \left\| u_0 \right\|_{q_0;\nu^{\phantom{a}}}^{\frac{q_0}{q_0+m-1}}  +  \left\| u_0 \right\|_{\frac{a}{2}(q_0+m-1);\nu} \right) \, .
\end{equation}
Choosing $a={2q_0}/{(q_0+m-1})$ we recover \eqref{eq: regSottoStep1} for $k=1$ (the constraint $a \le 1$ reads $q_0\leq m-1$).\\
Suppose now that \eqref{eq: regSottoStep1} holds for some $k$. Set $a={2^{k+1}q_0}/{(q_0+m-1)}$ in \eqref{eq: MoserSottoPasso4}: this is feasible providing that $q_0 \leq {(m-1)}/{(2^{k+1}-1)}$, and it leads to
\begin{equation} \label{eq: MoserSottoPassoInd}
\left\| u(t)  \right\|_{2^{k+1}q_0;\nu} \leq   D \left( t^{-\frac{1}{q_0+m-1}}  \left\| u_0 \right\|_{q_0;\nu^{\phantom{a}}}^{\frac{q_0}{q_0+m-1}} +  \left\| u_0 \right\|_{2^k q_0;\nu} \right) \, .
\end{equation}
Since \eqref{eq: MoserSottoPassoInd} holds for a.e.\ $t>0$, thanks to the semigroup property we are allowed to take (for a.e.\ $t$) ${t}/{2}$ as the origin of the time axis (that is, we replace $u_0$ with $u(t/2)$ and $t$ with $t/2$ in the right hand side). Applying \eqref{eq: regSottoStep1} (which holds by inductive hypothesis) and the non-expansivity of the norm $\| \cdot \|_{q_0;\nu}$ to the resulting right hand side of \eqref{eq: MoserSottoPassoInd}, we obtain \eqref{eq: regSottoStep1} for $k+1$ and $q_0 \in \left(1,{(m-1)}/{(2^{k+1}-1)} \right]$ as well.
\end{proof}
\end{lemma}
The following lemma provides an elementary numerical inequality.
\begin{lemma} \label{eq: lemminoNumerico}
Given $\alpha, \beta \in (0,1)$, with $\alpha > \beta$, there exists a constant $c=c(\alpha, \beta)>0$ such that $\forall x,y \in \mathbb{R}^{+} $
\begin{equation} \label{eq: minoreLemmino}
x^{-\alpha}y^{1-\alpha}+x^{-\beta}y^{1-\beta}+y \leq c(\alpha,\beta)(x^{-\alpha}y^{1-\alpha}+y) \, .
\end{equation}
\begin{proof}
We need to show that
\begin{equation} \label{eq: R}
R(x,y)=\frac{x^{-\beta}y^{1-\beta}} {x^{-\alpha}y^{1-\alpha}+y }
\end{equation}
is bounded in $\mathbb{R}^{+} \! \times \mathbb{R}^{+}$ by a constant which depends only on $\alpha$ and $\beta$. In order to do that, we can fix $y$ and find the zeros $x^{\ast}(y)$ of $R_{x}(\cdot,y)$ (in fact $0<\beta<\alpha$ implies $R(0^+,y)=R(+\infty,y)=0$). Through an explicit calculation we get
$$x^\ast(y)=\left( \frac{\alpha-\beta}{\beta}  \right)^{\frac{1}{\alpha}} y^{-1} \, .$$
Substituting such value in \eqref{eq: R} we easily obtain \eqref{eq: minoreLemmino} with
$$c(\alpha, \beta)=1 + \left( \frac{\beta}{\alpha} \right)^{\frac{\beta}{\alpha}} \left( 1-\frac{\beta}{\alpha} \right)^{1-\frac{\beta}{\alpha} } \, .  $$
\end{proof}
\end{lemma}
Before proving the main result of this section we comment that, as in the Dirichlet case, $L^\infty$ regularization need not hold in general.
\\[2mm]
{\noindent\bf Counterexample to the $\boldsymbol{{L}^\infty(\Omega)}$ regularization.}
Consider the domain $\Omega=\mathbb{R}$ and the weights $\rho_\nu(x)=\rho_\mu(x)=e^{-|x|}$ (one can regularize them in $x=0$ without significant modifications). These weights (see Section \ref{list}) satisfy the zero-mean Poincar\'e inequality \eqref{eq: poinIntroMed}.  Moreover, it is easy to show that $W^{1,2}(\mathbb{R};e^{-|x|},e^{-|x|})=W^{1,2}_0(\mathbb{R};e^{-|x|},e^{-|x|})$, so that in this case the Neumann problem coincides with the Dirichlet one: in particular, in order to prove that a certain function is a (sub)solution to \eqref{eq: pmeNeu}, one can neglect its behaviour at infinity and just test it on compactly supported functions in the weak formulation.

By means of computations analogous to the ones performed in the counterexample associated to Dirichlet boundary conditions, one can check that there exists a constant $B>0$ such that the function
\begin{equation*}
v(x,t)=\frac{\log(x^2+2)}{(1+B^{-1}(m-1)t)^{\frac{1}{m-1}}}
\end{equation*}
is a subsolution to \eqref{eq: pmeNeu}, so that $u(\cdot,t)\not\in L^\infty(\mathbb{R})$ for all $t\geq0$. This proves that for the initial datum $u_0(x)=\log(x^2+2)$, which belongs to $L^{q_0}(\mathbb{R};e^{-|x|}) \ \forall q_0 \in [1,\infty)$, there is no $L^{q_0}$-$L^\infty$ regularization.

\begin{theorem}\label{teo: regPmeNeu}
Let $u_0 \in L^{q_0}(\Omega;\nu)$. If the zero-mean Poincar\'e inequality \eqref{eq: poinIntroMed} holds, then the solution $u$ of \eqref{eq: pmeNeu} with initial datum $u_0$ satisfies the estimate
\begin{equation}\label{eq: stimaRegPmeNeu1}
\left\| u(t) \right\|_{\varrho;\nu} \leq K_1 \, t^{-\frac{\varrho-q_0}{\varrho(m-1)}} \left\| u_0 \right\|_{q_0;\nu^{\phantom{a}}}^{\frac{q_0}{\varrho}} \!\! e^{H \, \left\| u_0 \right\|_{q_0;\nu}^{m-1}\,t } \ \ \ \hbox{\rm for a.e.} \,\, t>0
\end{equation}
for all $q_0 \in [1,\infty)$ and $\varrho \in (q_0,\infty)$, where $K_1=K_1(\varrho,m,M_P,\nu(\Omega))$ and $H=H(m,M_P,\nu(\Omega))$ are suitable positive constants. Moreover, the estimate
\begin{equation} \label{eq: stimaRegPmeNeu2}
\left\| u(t) \right\|_{\varrho;\nu} \leq K_2  \left( t^{-\frac{\varrho-q_0}{\varrho(m-1)}}  \left\| u_0 \right\|_{q_0;\nu^{\phantom{a}}}^{\frac{q_0}{\varrho}} \!\!  + \left\| u_0 \right\|_{q_0;\nu}    \right) \ \ \ \hbox{\rm for a.e.} \,\, t>0
\end{equation}
holds true for all $q_0 \in (1,\infty)$ and $\varrho \in (q_0,\infty)$, where $K_2=K_2(q_0,\varrho,m,M_P,\nu(\Omega))>0$.

If instead only the weaker inequality \eqref{eq: weakPoi} is assumed to hold, then the bound \eqref{eq: stimaRegPmeNeu1} holds true for all $q_0 \in [1\vee(m-1),\infty)$, while \eqref{eq: stimaRegPmeNeu2} holds true for all $q_0 \in (1,\infty)\cap[m-1,\infty)$, both upon replacing $M_P$ by $W_P$.
\begin{proof}
To prove estimate \eqref{eq: stimaRegPmeNeu1} we adopt the techniques of \cite[Th. 1.3]{Gri10} and \cite[Th. 1.1]{BG05p}. We also refer to the notations used in the quoted theorems. Again, we shall consider $L^\infty$ data without loss of generality, and $1<q_0<\varrho$.

First of all, from \eqref{eq: poinIntroMed} we have
\begin{equation}\label{eq: 2-poinMod}
\left\| v \right\|_{2;\nu}^2 \leq 2 \left( M_P^2 \left\| \nabla{v} \right\|_{2;\mu}^2  + \left\| \overline{v} \right\|_{2;\nu}^2  \right)  \ \ \ \forall v \in W^{1,2}(\Omega;\nu,\mu) \, ;
\end{equation}
from that, proceeding exactly as in \cite[Th. 1.3]{Gri10}, it is straightforward to obtain the following family of logarithmic inequalities:
\begin{gather}
\left( J(r,v) + \frac{1}{2-r} \log{\varepsilon} \right) \frac{(2-r) \left\| v \right\|_{r;\nu}^2 }{2\,\varepsilon \, M_P^2} - \frac{\left\| \overline{v} \right\|_{2;\nu}^2} {M_P^2} \leq \left\| \nabla{v} \right\|_{2;\mu}^2\label{eq: p-logSob1}\\
\forall \varepsilon>0 \, , \ \forall r \in [1,2) \, , \ \forall v \in W^{1,2}(\Omega;\nu,\mu)\, .\nonumber
\end{gather}
Introducing a real function $q$ as in the proof of Theorem \ref{teo: regPmeDir}, explicit calculations yield
\begin{equation*}
\begin{aligned}
\frac{\mathrm{d}}{\mathrm{d}s} \log{\left\| u(s) \right\|_{q(s);\nu}} = & \frac{{q^\prime}(s)}{q(s)} J(q(s),u(s)) \\
&- \left( \frac{2}{q(s)+m-1} \right)^2 \frac{m(q(s)-1)}{\left\| u(s) \right\|_{q(s);\nu}^{q(s)}}  \left\| \nabla{\left( u^{\frac{q(s)+m-1}{2}}(s) \right) } \right\|_{2;\mu}^2 \, .
\end{aligned}
\end{equation*}
By applying \eqref{eq: p-logSob1} to $u^{(q+m-1)/2}$ in the right hand side of the equation above, choosing suitably $r$, $\varepsilon$ and exploiting interpolation inequalities between the norms $\| \cdot \|_{1 \vee(m-1);\nu}$, $\| \cdot \|_{(q+m-1)/2;\nu}$, $\| \cdot \|_{q;\nu}$, we get:
\begin{equation*}
\begin{aligned}
\frac{\mathrm{d}}{\mathrm{d}s} \log{\left\| u \right\|_{q;\nu}}  \leq & - \frac{{q^\prime}}{q(m-1)} \log{\left[ \frac{2q(q-1)m(m-1)}{{q^\prime} (q+m-1)^2 M_P^2} \right] } -\frac{{q^\prime}}{q} \log{\left\| u \right\|_{q;\nu}} \\
&+ \left( \frac{2}{q+m-1} \right)^2 \frac{m(q-1) }{M_P^2 \, \nu(\Omega)^{1\wedge(m-1)}} \left\| u_0 \right\|_{1\vee(m-1);\nu}^{m-1} \, .
\end{aligned}
\end{equation*}
For any given $t>0$ let us set $q(s)=q_0+\frac{s}{t}(\varrho-q_0)$ and solve the resulting differential inequality in the variable $y(s)=\log{\left\| u(s) \right\|_{q(s);\nu}}$. Standard computations give estimate \eqref{eq: stimaRegPmeNeu1} (the case $q_0=1$ is handled by letting $q_0 \downarrow 1$) with $1\vee(m-1)$ instead of $q_0$ in the norm appearing in the exponential term. At the end of the proof we shall show how it is possible to replace there $1\vee(m-1)$ with $q_0$.

In order to prove estimate \eqref{eq: stimaRegPmeNeu2} it is enough to carry out a single step of the Moser iteration. Firstly we assume $q_0 \in (1,\infty) \cap [m-1,\infty)$, and with no loss of generality $\nu(\Omega)=1$. From \eqref{eq: 2-poinMod} we have, in particular,
\begin{equation}\label{eq: poinMedModif}
\frac{1}{2M_P^2} \left\| v \right\|^2_{2;\nu} - \frac{1}{M_P^2} \left\| v \right\|^2_{1;\nu} \leq \| \nabla{v} \|_{2;\mu}^2 \ \ \ \forall v \in W^{1,2}(\Omega;\nu,\mu) \, ;
\end{equation}
setting $T=t$ and $q=q_0-1$ in estimate \eqref{eq: stimeLemmaDir1} and applying \eqref{eq: poinMedModif} to $u^{(q_0+m-1)/2}$, we obtain:
\begin{equation}\label{eq: MoserPoin1}
\frac{4(q_0-1)q_0 m}{M_P^2(m+q_0-1)^2} \int_0^t \left(  \frac{1}{2} \left\| u(s) \right\|_{q_0+m-1;\nu}^{q_0+m-1} - \left\|u(s) \right\|_{\frac{q_0+m-1}{2};\nu}^{q_0+m-1} \right) \, \mathrm{d}s  \leq \left\| u_0 \right\|^{q_0}_{q_0;\nu}  \, .
\end{equation}
Since $ q_0 \geq m-1$ and $\nu(\Omega)<\infty$, the quantity $\| u \|_{(q_0+m-1)/{2};\nu}$ can be controlled from above by $\| u \|_{q_0}$; using this fact and the non-expansivity of the norms $\| \cdot \|_{q_0+m-1;\nu}$ and $\| \cdot\|_{q_0;\nu}$, after some calculations we arrive at
\begin{equation} \label{eq: MoserPoin2}
\left\| u(t) \right\|_{q_0+m-1;\nu} \leq D \left( t^{-\frac{1}{q_0+m-1}} \left\|u_0 \right\|_{q_0;\nu^{\phantom{a}}}^{\frac{q_0}{q_0+m-1}} +  \left\| u_0 \right\|_{q_0;\nu} \right) \, ,
\end{equation}
where $D>0$ is a constant possibly depending on $q_0$, $\varrho$, $m$, $M_P$ which may change from line to line. Clearly \eqref{eq: MoserPoin2} only provides a regularization from $L^{q_0}(\Omega;\nu)$ to $L^{q_0+m-1}(\Omega;\nu)$. However, by means of induction, interpolation inequalities and Lemma \ref{eq: lemminoNumerico} it is not difficult to get from \eqref{eq: MoserPoin2} the more general \eqref{eq: stimaRegPmeNeu2} for all $q_0 \in (1,\infty) \cap [m-1,\infty)$ and $\varrho \in (q_0,\infty)$.
To remove the constraint $q_0 \geq m-1$ we need to exploit Lemma \ref{eq: regSottoNeuPoin1}. Suppose $m>2$ (otherwise there is nothing to prove): given $q_0 \in (1,m-1)$, of course there exists an integer $k$ such that
\begin{equation}\label{eq: hpQ0Dolb}
\frac{q_0+m-1}{2} \leq 2^k q_0\, , \ \ \  q_0 \leq \frac{m-1}{2^k-1} \, .
\end{equation}
From \eqref{eq: MoserPoin1} and the first inequality in \eqref{eq: hpQ0Dolb} one deduces that
\begin{equation} \label{eq: hpQ0Dolb2}
\left\| u(t) \right\|_{q_0+m-1;\nu} \leq D \left( t^{-\frac{1}{q_0+m-1}} \left\|u_0 \right\|_{q_0;\nu^{\phantom{a}}}^{\frac{q_0}{q_0+m-1}} +  \left\| u_0 \right\|_{2^k q_0;\nu} \right) \, .
\end{equation}
Shifting the origin of the time axis to $t/2$ in \eqref{eq: hpQ0Dolb2} and applying to the so modified right hand side estimate \eqref{eq: regSottoStep1} evaluated at time $t/2$ (which is feasible in view of the second inequality in \eqref{eq: hpQ0Dolb}), we are able to conclude that \eqref{eq: MoserPoin2} holds for all $q_0 \in (1,\infty)$, and so \eqref{eq: stimaRegPmeNeu2} by arguing exactly as in the case $ q_0 \geq m-1$.

The initial assumption $\nu(\Omega)=1$ is removable by spatial scaling. In fact, if $u(\mathbf{x},t)$ is a solution of \eqref{eq: pmeNeu} on the domain $\Omega$ of measure $V=\nu(\Omega)$,  with respect to the weights $\rho_\nu(\mathbf{x}),\rho_\mu(\mathbf{x})$ and with initial datum $u_0(\mathbf{x})$, then
\begin{equation}\label{eq: scalingSp}
\widetilde{u}(\widetilde{\mathbf{x}},t)=V^{-\frac{2}{N(m-1)}}u\left(V^{\frac{1}{N}} \widetilde{\mathbf{x}},t\right)
\end{equation}
is also a solution of \eqref{eq: pmeNeu} on the domain $\widetilde{\Omega}=\Omega/V^{\frac{1}{N}}$ of measure $1$, with respect to the weights
$$\widetilde{\rho}_\nu(\widetilde{\mathbf{x}})=\rho_\nu\left(V^{\frac{1}{N}}\widetilde{\mathbf{x}}\right)\, , \ \widetilde{\rho}_\mu(\widetilde{\mathbf{x}})=\rho_\mu\left(V^{\frac{1}{N}}\widetilde{\mathbf{x}}\right)$$
and with initial datum
$$\widetilde{u}_0(\widetilde{\mathbf{x}})=V^{-\frac{2}{N(m-1)}}u_0\left(V^{\frac{1}{N}} \widetilde{\mathbf{x}}\right) \, . $$
From that, one applies \eqref{eq: stimaRegPmeNeu2} to $\widetilde{u}$ and then goes back to the original solution $u$ through \eqref{eq: scalingSp} and
$$ \left\| \widetilde{u} \right\|_{q;\widetilde{\nu}}= V^{-\frac{2}{N(m-1)}-\frac{1}{q}} \left\| u \right\|_{q;\nu} \, , \ M_P(\widetilde{\Omega})=V^{-\frac{1}{N}}M_P(\Omega) \, , $$
thus obtaining \eqref{eq: stimaRegPmeNeu2} for $u$ with a multiplicative constant that in general will depend on $\nu(\Omega)$ as well.

Finally, we are left to show that the estimate
\begin{equation}\label{eq: stimaRegPmeNeuMm1}
\left\| u(t) \right\|_{\varrho;\nu} \leq K_1 \, t^{-\frac{\varrho-q_0}{\varrho(m-1)}} \left\| u_0 \right\|_{q_0;\nu^{\phantom{a}}}^{\frac{q_0}{\varrho}} \!\! e^{H \, \left\| u_0 \right\|_{1 \vee(m-1);\nu}^{m-1}\,t } \ \ \ \hbox{\rm for a.e.} \,\, t>0 \, ,
\end{equation}
whose validity we proved above, implies \eqref{eq: stimaRegPmeNeu1}. Let us suppose $m>2$, otherwise one controls $\|u_0\|_{1;\nu}$ with $\|u_0\|_{q_0;\nu}$. Lemma \ref{lem: dolbPor}, which we shall prove below, gives as a byproduct the validity of the estimate
\begin{equation} \label{eq: lemmaUtile}
\left\| u(t) \right\|_{2;\nu} \leq C_1 \, \left( t^{-\frac{1}{m-1}} + \left\| u_0 \right\|_{1;\nu} \right) \, \ \ \ \hbox{\rm for a.e.} \, \, t>0
\end{equation}
for a suitable constant $C_1>0$ depending on $m$, $M_P$ and $\nu(\Omega)$. Hence from \eqref{eq: lemmaUtile} and the regularity estimate \eqref{eq: stimaRegPmeNeu2} we can infer in turn that
\begin{equation} \label{eq: lemmaUtile2}
\left\| u(t) \right\|_{m-1;\nu} \leq C_2 \, \left( t^{-\frac{1}{m-1}} + \left\| u_0 \right\|_{1;\nu} \right) \, \ \ \ \hbox{\rm for a.e.} \, \, t>0
\end{equation}
for another constant $C_2>0$ depending on the same quantities. Indeed, if $m \leq 3$ \eqref{eq: lemmaUtile2} clearly follows from \eqref{eq: lemmaUtile}, else one combines \eqref{eq: stimaRegPmeNeu2} with the choices $\varrho=m-1$, $q_0=2$ and \eqref{eq: lemmaUtile} by means of the usual $t/2$-shift argument, which again entails \eqref{eq: lemmaUtile2} (up to a different constant $C_2$). It is now plain that \eqref{eq: lemmaUtile2} and \eqref{eq: stimaRegPmeNeuMm1} give the desired result thanks to another $t/2$-shift argument. The dependence of the constants $K_1$ and $H$ on $q_0$ has been implicitly absorbed into $\varrho$, since $q_0 < \varrho$ and they remain bounded as $q_0$ varies in the interval $[1,\varrho]$.

The last statement of the theorem is a mere consequence of the fact that the validity of inequality \eqref{eq: poinMedModif} is sufficient in order to prove \eqref{eq: stimaRegPmeNeu2} at least for $q_0 \in (1,\infty) \cap [m-1,\infty)$. The same applies for \eqref{eq: stimaRegPmeNeu1} provided $q_0 \in [1\vee(m-1),\infty)$. The passage to data in $L^{q_0}(\Omega;\nu)$ with $q_0 \in (1,m-1)$, instead, needs Lemma \ref{eq: regSottoNeuPoin1} and Lemma \ref{lem: dolbPor}, which both require the original zero-mean Poincar\'e inequality \eqref{eq: poinIntroMed}.
\end{proof}
\end{theorem}

\begin{remark}\rm
 Our proof of the validity of \eqref{eq: stimaRegPmeNeu2} is not extendible to the case $q_0=1$, since the constant $D$ in \eqref{eq: MoserPoin2} diverges as $q_0 \downarrow 1$. Nevertheless estimate \eqref{eq: stimaRegPmeNeu1} also ensures an $L^{1}$-$L^{\varrho}$ regularization with the same short-time rate one would expect from \eqref{eq: stimaRegPmeNeu2} by letting $q_0\downarrow 1$.
\end{remark}
\vskip 2mm
\noindent{\bf Converse implications.}
In Section \ref{sec: regAs} we saw that the validity of a suitable estimate for solutions to the Dirichlet problem \eqref{eq: pmeDir} implies, in turn, the validity of the Poincar\'e inequality \eqref{eq: poinIntro}. For the Neumann problem we are able to prove, with analogous techniques, a slightly weaker but similar result.
\begin{theorem}\label{teo: impInvNeu}
Suppose $\nu(\Omega)<\infty$. If there exist a constant $K>0$ and a given $q_0 \in [m,m+1)$ such that, for all $u_0\in L^{q_0}(\Omega;\nu)$, the solution $u$ of \eqref{eq: pmeNeu} corresponding to the initial datum $u_0$ satisfies  the estimate
\begin{equation*}
\left\| u(t) \right\|_{m+1;\nu} \leq K  \left( t^{-\frac{m+1-q_0}{(m+1)(m-1)}}  \left\| u_0 \right\|_{q_0;\nu^{\phantom{a}}}^{\frac{q_0}{m+1}}  + \left\| u_0 \right\|_{q_0;\nu}    \right) \ \ \ \hbox{\rm for a.e.} \, \, t>0
\end{equation*}
(namely \eqref{eq: stimaRegPmeNeu2} when $\varrho=m+1$), then there exists a constant $B>0$ such that the functional inequality
\begin{equation}\label{eq: impInvNeuDF}
\left\|  v \right\|_{2;\nu} \leq B \left( \left\| \nabla{v}  \right\|_{2;\mu} + \left\| v \right\|_{\frac{q_0}{m};\nu} \right) \ \ \ \forall v \in W^{1,2}(\Omega;\nu,\mu)
\end{equation}
holds as well.
\begin{proof}
One starts considering an initial nonnegative datum $u_0$ belonging to $L^\infty(\Omega) \cap W^{1,2}(\Omega;\nu,\mu) $. Proceeding along the lines of the proof of Theorem \ref{teo: pmeDirImpInv}, no major difficulty arises in obtaining the following inequality:
\begin{equation*}\label{eq: impInvNeuBase0}
\left\| u_0 \right\|_{m+1;\nu} \leq  B \, \left(  \left\| \nabla{(u_0^m)} \right\|_{2;\mu^{\phantom{a}}}^{\frac{2(m+1-q_0)}{(m+1)(2m-q_0)}}    \left\| u_0 \right\|_{q_0;\nu^{\phantom{a}}}^{\frac{q_0(m-1)}{(m+1)(2m-q_0)}} + \left\| u_0 \right\|_{q_0;\nu}  \right) \, ,
\end{equation*}
where $B$ is a suitable positive constant. Upon setting
$$ \xi_n(x)=\frac{1}{n} \chi_{\left[0,\frac{1}{n}\right]}(x) + x \, \chi_{\left(\frac{1}{n},\infty\right)}(x) \, , \ \ \ \mathcal{W}(f)= \left\| \nabla{f} \right\|_{2;\mu} + \left\| f \right\|_{\frac{q_0}{m};\nu} $$
and exploiting Proposition \ref{pro: dens2}, the result follows again as in the proof of Theorem \ref{teo: pmeDirImpInv}.
\end{proof}
\end{theorem}
Note that \eqref{eq: impInvNeuDF} is equivalent to the fact that the spaces $V^{q_0/m}(\Omega;\nu,\mu)$ and $W^{1,2}(\Omega;\nu,\mu)$ coincide. Also, if $q_0=m$ \eqref{eq: impInvNeuDF} becomes \eqref{eq: weakPoi}, so that from Theorems \ref{teo: regPmeNeu} and \ref{teo: impInvNeu} we can get the following
\begin{corollary}\label{equiv}
Suppose $\nu(\Omega)<\infty$. Consider, for all $u_0\in L^{m}(\Omega;\nu)$, the solution $u$ of \eqref{eq: pmeNeu} corresponding to the initial datum $u_0$. The existence of constants $\varrho\ge m+1$, $K>0$ such that the estimate
\begin{equation}\label{eq: impInvNeuReg}
\left\| u(t) \right\|_{\varrho;\nu} \leq K  \left( t^{-\frac{\varrho-m}{\varrho(m-1)}}  \left\| u_0 \right\|_{m;\nu^{\phantom{a}}}^{\frac{m}{\varrho}}  + \left\| u_0 \right\|_{m;\nu}    \right) \ \ \ \hbox{\rm for a.e.} \, \, t>0
\end{equation}
holds, is \emph{equivalent} to the validity of inequality \eqref{eq: weakPoi}. In particular, the validity of \eqref{eq: impInvNeuReg} for a given $\bar{\varrho}\ge m+1$ implies the validity of the same estimate for \emph{all} $\varrho\in[ m+1,\infty)$ and, more generally, the validity of \eqref{eq: stimaRegPmeNeu2} for any $q_0 \in (1,\infty)\cap[m-1,\infty)$ and $\varrho\in (q_0,\infty)$.
\end{corollary}

\subsection{Asymptotic estimates}\label{sec: stimeAsiNeuPon}
As already mentioned, estimate \eqref{eq: stimaRegPmeNeu1} diverges as $t \rightarrow \infty$, so it prevents us from obtaining any information about the asymptotic behaviour of $u(\cdot,t)$. On the other hand, estimate \eqref{eq: stimaRegPmeNeu2} only allows us to deduce that
\begin{equation*}\label{eq: infoAsiPoin}
\limsup_{t\rightarrow \infty}{\left\| u(t) \right\|_{\varrho;\nu}} \leq K_2 \left\| u_0 \right\|_{q_0;\nu} \, .
\end{equation*}
In order to study in deeper detail such asymptotic behaviour, it is crucial to be able to suitably handle the quantity (recall by Proposition \ref{pro: consMed} that the mean value of $u_0$ is preserved)
$$\frac{\mathrm{d}}{\mathrm{d}s} \left\| u(s)-\overline{u} \right\|_{\varrho;\nu}^\varrho \, . $$
To this end, the next lemma will be fundamental.

\begin{lemma}\label{lem: GP}
Suppose that the zero-mean Poincar\'e inequality \eqref{eq: poinIntroMed} holds. Let $\Phi:\mathbb{R}\rightarrow \mathbb{R}$ be a continuous and increasing function with the following properties:
\begin{equation}\label{eq: GPMin}
 \lim_{x \rightarrow 0}{\frac{\Phi(x)}{x^r}}=l_0 \, , \ \lim_{x \rightarrow -\infty}{\frac{\Phi(x)}{x^r}}=l_{-} \, , \ \lim_{x \rightarrow +\infty}{\frac{\Phi(x)}{x^r}}=l_{+}
\end{equation}
for some constants $r \geq {1}/{2}$ and $l_0$, $l_{-}$, $l_{+} \in (0,+\infty)$. Then there exists a constant $C_{\Phi}>0$ such that for every function $\xi \in L^1(\Omega;\nu)$ such that $\overline{\xi}=0$ and $\Phi(\xi) \in W^{1,2}(\Omega;\nu,\mu)$ the inequality
\begin{equation}\label{diseqGP}
\left\| \Phi(\xi) \right\|_{2;\nu} \leq C_{\Phi} \left\| \nabla{\Phi(\xi)} \right\|_{2;\mu}
\end{equation}
holds.
\begin{proof}
We proceed by contradiction. Should the assertion be false, then there exists a sequence of functions $\{\xi_n\} \subset \{\xi \in L^1(\Omega;\nu): \ \overline{\xi}=0 \, , \ \Phi(\xi) \in W^{1,2}(\Omega;\nu,\mu) \}$ (not identically zero) such that
\begin{equation*}
\left\| \nabla{\Phi(\xi_n)} \right\|_{2;\mu} \leq  \frac{1}{n} \left\| \Phi(\xi_n) \right\|_{2;\nu} \, .
\end{equation*}
Let us set $a_n=\| \Phi(\xi_n) \|_{2;\nu} $ and
$$ \Psi_n(\xi_n)=\frac{\Phi(\xi_n)}{a_n} \, . $$
Clearly,
\begin{equation}\label{eq: GPpropPsi}
\left\| \Psi_n(\xi_n) \right\|_{2;\nu}=1 \, , \ \  \  \left\| \nabla{\Psi_n(\xi_n)} \right\|_{2;\mu} \leq \frac{1}{n} \, .
\end{equation}
Applying the zero-mean Poincar\'e inequality to the sequence $\{ \Psi_n(\xi_n) \}$ and exploiting the second inequality in \eqref{eq: GPpropPsi}, we have that
\begin{equation}\label{eq: GPstep1}
\left\| \Psi_n(\xi_n)-\overline{\Psi_n(\xi_n)} \right\|_{2;\nu} \leq \frac{M_P}{n} \, .
\end{equation}
The inequality just obtained and the normalization condition in \eqref{eq: GPpropPsi} together imply that the sequence of real numbers $\{ {\overline{\Psi_n(\xi_n)}} \}$ is bounded, hence up to subsequences it converges to some limit $c_0$. This and again \eqref{eq: GPstep1} allow us to deduce that
\begin{equation*}
\left\| \Psi_n(\xi_n) - c_0 \right\|_{2;\nu} \rightarrow 0 \, ,
\end{equation*}
that is, up to subsequences,
\begin{equation}\label{eq: GPstep3}
 \Psi_n(\xi_n(\mathbf{x})) \rightarrow c_0  \ \ \ \textnormal{for a.e.} \ \mathbf{x} \in \Omega \, .
\end{equation}
The normalization condition just mentioned prevents $c_0$ from being zero. Now we need to distinguish three cases according to the value of the quantity
$$ a_\infty=\lim_{n\rightarrow \infty} a_n  \, ,$$
the limit above existing possibly passing again to a subsequence.
If $a_\infty \in (0,+\infty) $, from the continuity of $\Phi$ (and so of $\Phi^{-1}$) it is easy to infer that
\begin{equation*}
  \xi_n(\mathbf{x}) \rightarrow \Phi^{-1}(a_\infty c_0)  \neq 0  \ \ \ \textnormal{for a.e.} \ \mathbf{x} \in \Omega \, .
\end{equation*}
When $a_\infty=0$ or $a_\infty=+\infty$ things are slightly more delicate. Let us begin with the case $a_\infty=0$. By the definition and the properties of $\Phi$, and in view of \eqref{eq: GPstep3}, it follows that
\begin{equation*}
\xi_n(\mathbf{x}) \rightarrow 0 \ \ \ \textnormal{for a.e.} \ \mathbf{x} \in \Omega \, ;
\end{equation*}
hence, exploiting the first equality in \eqref{eq: GPMin} and again \eqref{eq: GPstep3},
\begin{equation*}
 \mathcal{Z}_n(\mathbf{x})=\frac{\xi_n(\mathbf{x})}{a_n^{{1}/{r}}} = \left( \frac{[\xi_n(\mathbf{x})]^r}{\Phi(\xi_n(\mathbf{x}))} \Psi_n(\xi_n(\mathbf{x})) \right)^{\frac{1}{r}} \rightarrow \left(\frac{c_0}{l_0}\right)^{\frac{1}{r}} \neq 0       \ \ \ \textnormal{for a.e.} \ \mathbf{x} \in \Omega \, .
\end{equation*}
If instead $a_\infty = + \infty$ one argues likewise. In fact, suppose $c_0>0$. From the properties of $\Phi$ and \eqref{eq: GPstep3} we deduce that
\begin{equation*}
\xi_n(\mathbf{x}) \rightarrow +\infty \ \ \ \textnormal{for a.e.} \ \mathbf{x} \in \Omega \, ,
\end{equation*}
which, thanks to the third equality in \eqref{eq: GPMin}, implies
\begin{equation*}
 \mathcal{Z}_n(\mathbf{x})=\frac{\xi_n(\mathbf{x})}{a_n^{{1}/{r}}} = \left( \frac{[\xi_n(\mathbf{x})]^r}{\Phi(\xi_n(\mathbf{x}))} \Psi_n(\xi_n(\mathbf{x})) \right)^{\frac{1}{r}} \rightarrow \left(\frac{c_0}{l_{+}}\right)^{\frac{1}{r}} \neq 0       \ \ \ \textnormal{for a.e.} \ \mathbf{x} \in \Omega \, .
\end{equation*}
When $c_0<0$ one proves similarly that
\begin{equation*}
 \mathcal{Z}_n(\mathbf{x}) \rightarrow \left(\frac{c_0}{l_{-}}\right)^{\frac{1}{r}} \neq 0    \ \ \ \textnormal{for a.e.} \ \mathbf{x} \in \Omega \, .
\end{equation*}
Hence in any case the sequence $\{ \mathcal{Z}_n \}$ converges pointwise to a nonzero constant. Since obviously $\overline{\mathcal{Z}_n}=0$ and the mean value operator is continuous in $L^1(\Omega;\nu)$, we come to a contradiction as soon as we prove that $\{ \mathcal{Z}_n \}$ also converges in $L^1(\Omega;\nu)$ to such nonzero constant. To this end, note that from Egoroff's Theorem it is enough to show that the quantity
$$ \int_E \left| \mathcal{Z}_n(\mathbf{x}) \right|  \, \mathrm{d}\nu $$
converges to zero uniformly as $n\rightarrow \infty$ and $|E|\rightarrow 0$. First of all, observe that \eqref{eq: GPMin}, together with the continuity and the monotonicity of $\Phi$, imply the existence of a constant $D>0$ such that
\begin{equation*}
D^{-1} |x|^r \leq |\Phi(x)| \leq D |x|^r \ \ \ \forall x \in \mathbb{R} \, .
\end{equation*}
As a consequence,
\begin{equation*}
\begin{aligned}
\int_E \left| \mathcal{Z}_n(\mathbf{x}) \right|  \, \mathrm{d}\nu = \int_E \frac{\left| \xi_n(\mathbf{x}) \right|}{a_n^{1/r}}  \, \mathrm{d}\nu & \leq D^{\frac{1}{r}} \int_E |\Psi_n(\xi_n( \mathbf{x} ))|^{\frac{1}{r}}  \, \mathrm{d}\nu \\ & \leq D^{\frac{1}{r}} |E|^{1-\frac{1}{2r}} \left( \int_E |\Psi_n(\xi_n(\mathbf{x}))|^2 \, \mathrm{d}\nu  \right)^{\frac{1}{2r}}\end{aligned}
\end{equation*}
so that the quantity
\begin{equation*}
\int_E |\Psi_n(\xi_n(\mathbf{x}))|^2 \, \mathrm{d}\nu
\end{equation*}
indeed goes to zero uniformly as $n\rightarrow \infty$ and $|E|\rightarrow 0$ since
\begin{equation*}
\int_E |\Psi_n(\xi_n(\mathbf{x}))|^2 \, \mathrm{d}\nu  \leq  2 \left( \int_\Omega |\Psi_n(\xi_n(\mathbf{x}))-c_0|^2  \, \mathrm{d}\nu  + |E| c_0^2   \right) \, .
\end{equation*}
Therefore we conclude that $\{\mathcal{Z}_n\}$ converges in $L^1(\Omega;\nu)$ to a nonzero constant with zero mean, a contradiction.
\end{proof}
\end{lemma}
\begin{remark}\label{oss: GP}
When $\rho_\nu=\rho_\mu \equiv 1$ and $\Phi(x)=x^m$ ($m>1$), the result had already been proved in \cite[Lem. 3.2]{AR81} . However, the proof provided therein exploits the compactness of the embedding $H^1(\Omega) \hookrightarrow L^2(\Omega)$. The proof of Lemma \ref{lem: GP} does not need compactness. Note that it is essential that the behaviour of $\Phi(x)$ as $x \rightarrow 0$ and $x \rightarrow \pm \infty $ is given by the same power of $x$. If, for example, $\Phi(x) \sim x^{r_1}$ as $x \rightarrow 0$ and $\Phi(x) \sim x^{r_2}$ as $x \rightarrow \pm \infty$ with $r_1 \neq r_2$ our proof does not work (one loses control of $\int_E |\mathcal{Z}_n| \, \mathrm{d}\nu $ either when $a_\infty=0$ or $a_\infty=+\infty$).
\end{remark}
We are now ready to prove an asymptotic estimate for zero-mean solutions. With a slight abuse of notation, we shall indicate below by $C_{x^{a}}$ the value of $C_\Phi$ (see formula \eqref{diseqGP}) when $\Phi(x)=x^a$.
\begin{theorem}\label{teo: absbEMnullaPme}
Let $q_0 \in [1,\infty)$, $u_0 \in L^{q_0}(\Omega;\nu)$ and $\overline{u_0}=0$. If the zero-mean Poincar\'e inequality \eqref{eq: poinIntroMed} holds, then the solution $u$ of \eqref{eq: pmeNeu} with initial datum $u_0$ satisfies the following absolute bound:
\begin{equation}\label{eq: absbMNullaPme}
\left\| u(t) \right\|_{\varrho;\nu} \leq Q_2 \, t^{-\frac{1}{m-1}} \ \ \ \hbox{\rm for a.e.} \,\, t>0 \, ,
\end{equation}
where $\varrho \in [1,\infty)$ and $Q_2$ is a constant depending only on $\varrho$, $m$, $M_P$, $C_{x^{m}}$ and $\nu(\Omega)$. As a consequence, for initial data with zero mean \eqref{eq: stimaRegPmeNeu1} becomes
\begin{equation}\label{eq: stimaRegPmeNeuMnulla}
\left\| u(t) \right\|_{\varrho;\nu} \leq Q_1 \, t^{-\frac{\varrho-q_0}{\varrho(m-1)}} \left\| u_0 \right\|_{q_0;\nu^{\phantom{a}}}^{\frac{q_0}{\varrho}}  \ \hbox{\rm for a.e.} \,\,  t>0 \, ,
\end{equation}
for a suitable constant $Q_1$ depending only on $\varrho$, $m$, $M_P$ and $C_{x^{m}}$.
\begin{proof}
Given $\varrho \in (1,\infty)$ (and, as usual, assuming $u_0 \in L^\infty(\Omega)$), consider the (formal) identity
\begin{equation*}
\frac{\mathrm{d}}{\mathrm{d}s} \left\| u(s) \right\|_{\varrho;\nu}^\varrho = - \left( \frac{2}{\varrho+m-1} \right)^2 m  \varrho ( \varrho-1 ) \left\| \nabla \left( u^{\frac{\varrho+m-1}{2}} \right) (s) \right\|_{2;\mu}^2 \, .
\end{equation*}
In order to handle the right hand side, we can apply Lemma \ref{lem: GP} with the choice $\Phi(x)=x^{(\varrho+m-1)/2}$, which provides us with a constant $C_{x^{(\varrho+m-1)/2}}$ such that
\begin{equation*}
\left\| u \right\|_{\varrho+m-1;\nu}^{\varrho+m-1} \leq C_{x^{(\varrho+m-1)/2}}^2 \left\| \nabla\left(u^{\frac{(\varrho+m-1)}{2}}\right) \right\|_{2;\mu}^2  \, .
\end{equation*}
From now on, to obtain the absolute bound \eqref{eq: absbMNullaPme}, one proceeds exactly as in the proof of Theorem \ref{teo: regPmeDir}, replacing $C_P$ with $C_{x^{(\varrho+m-1)/2}}$ (the case $\varrho=1$ is recovered by the finiteness of $\nu(\Omega)$). If $\varrho < m+1$, since $\|u(t) \|_{\varrho;\nu} \leq \nu(\Omega)^{1/\varrho - 1/(m+1)} \| u(t) \|_{m+1;\nu} $, $C_{x^{(\varrho+m-1)/2}}$ can in turn be replaced by $C_{x^m}$ into constant $Q_2$.

Now we have to prove estimate \eqref{eq: stimaRegPmeNeuMnulla}. First of all, let us rewrite \eqref{eq: stimaRegPmeNeuMm1} with the time origin shifted to $t/2$ (this choice, together with the following, are allowed for almost every $t$) and writing explicitly the constant $H$ appearing there. Exploiting the non-expansivity of the norm $\| \cdot \|_{q_0;\nu}$, we have:
\begin{equation}\label{eq: stimaRegPmeNeu1Tmezzi}
\left\| u(t) \right\|_{\varrho;\nu} \leq K_1 \, \left( \frac{t}{2} \right)^{-\frac{\varrho-q_0}{\varrho(m-1)}} \left\| u_0 \right\|_{q_0;\nu^{\phantom{a}}}^{\frac{q_0}{\varrho}} \!\! e^{\frac{2}{M_P^2 \, \nu(\Omega)^{1\wedge(m-1)}} \left\| u\left(t/2\right) \right\|_{1\vee (m-1);\nu}^{m-1}\, t } \, ;
\end{equation}
applying to the exponential term in \eqref{eq: stimaRegPmeNeu1Tmezzi} the absolute bound just proved, we easily deduce \eqref{eq: stimaRegPmeNeuMnulla}. The fact that the constant $Q_1$ is independent of $\nu(\Omega)$ can be shown as follows. First one proves that, setting $\varrho=m+1$, the constant $Q_2$ in \eqref{eq: absbMNullaPme} depends on $\nu(\Omega)$ through a multiplication by $\nu(\Omega)^{1/(m+1)}$. Therefore
\begin{equation*}
\frac{2}{M_P^2 \, \nu(\Omega)^{1\wedge(m-1)}} \left\| u\left(t/2\right) \right\|_{1\vee (m-1);\nu}^{m-1}\, t \leq Q_0(m,M_P,C_{x^m}) \, .
\end{equation*}
Afterwards one notices, from the proof of Theorem \ref{teo: regPmeNeu}, that the constant $K_1$ appearing in \eqref{eq: stimaRegPmeNeu1Tmezzi} depends only on $\varrho$, $m$ and $M_P$.

Estimate \eqref{eq: stimaRegPmeNeuMnulla} can be also used to show that the constant $Q_2$ in \eqref{eq: absbMNullaPme} depends in turn on a constant coming from Lemma \ref{lem: GP} which is at most $C_{x^m}$ even if $\varrho > m+1$. To this end, it is enough to perform the usual $t/2$-shift in \eqref{eq: stimaRegPmeNeuMnulla} (with $q_0=m+1$) and use the absolute bound itself (with $\varrho=m+1$) on the right hand side.
\end{proof}
\end{theorem}
The informations we provided by the previous corollary only concern the asymptotic behaviour of zero-mean solutions. However, this does not allow us to infer anything about \emph{nonzero}-mean solutions (i.e.\ the solution of \eqref{eq: pmeNeu} corresponding to the initial datum $u_0+c$, for $c \in \mathbb{R}\setminus\{0\}$, is \emph{not} the solution corresponding to $u_0$ plus $c$). In order to obtain such informations also when $\overline{u}\neq0$, we begin with an important lemma.
\begin{lemma}\label{lem: dolbPor}
If the zero-mean Poincar\'e inequality \eqref{eq: poinIntroMed} holds, then there exists a constant $Q>0$, possibly depending on $\varrho \in [1,2]$, $m$, $M_P$ and $\nu(\Omega)$, such that for all solutions to \eqref{eq: pmeNeu} the following absolute bound holds:
\begin{equation}\label{eq: convL2Dolb}
\left\| u(t)-\overline{u} \right\|_{\varrho;\nu} \leq Q \, t^{-\frac{1}{m-1}} \, \ \ \ \hbox{\rm for a.e.} \, \, t>0 \, .
\end{equation}
\begin{proof}
The result had already been proved in \cite[Th. 4.5]{DG+08} when $\rho_\nu=\rho_\mu$. For the sake of completeness we repeat the main lines in the case $\rho_\nu \neq \rho_\mu$. To this end notice that, formally:
\begin{equation}\label{eq: derMediaNeu}\begin{aligned}
\frac{\mathrm{d}}{\mathrm{d}s} \left\| u(s) -\overline{u} \right\|_{2;\nu}^2 &= -\frac{8m}{(m+1)^2} \left\| \nabla{\left(u^{\frac{m+1}{2}}\right)}(s)  \right\|_{2;\mu}^2\\ &\leq -\frac{8m}{(m+1)^2 \, M_{P,\,{2}/{(m+1)}}^{m+1}} \left\| u(s)-\overline{u} \right\|_{2;\nu^{\phantom{a}}}^{{m+1}} \, ,\end{aligned}
\end{equation}
where we have used Lemma \ref{lem: Dolb} (with $a=2/(m+1)$) applied to the function $u^{(m+1)/2}$, $M_{P,\,{2}/{(m+1)}}$ being the constant appearing in the statement of such lemma. Solving the above differential inequality in the variable $y(s)=\| u(s)-\overline{u} \|_{2;\nu}^2$, we get \eqref{eq: convL2Dolb} for $\varrho=2$. The case $\varrho \in [1,2)$ follows from the finiteness of the measure.
\end{proof}
\end{lemma}
From Lemma \ref{lem: dolbPor} and the smoothing results provided by Theorem \ref{teo: regPmeNeu}, it is not difficult to deduce asymptotic estimates for $\| u(t)-\overline{u}\|_{\varrho;\nu}$ also when $\varrho \in (2,\infty)$.
\begin{corollary}\label{teo: AsiMediaNN}
Let $q_0 \in [1,\infty)$ and $\varrho \in (2\vee q_0,\infty)$. If the zero-mean Poincar\'e inequality \eqref{eq: poinIntroMed} holds, then the solution $u$ of \eqref{eq: pmeNeu} with initial datum $u_0 \in L^{q_0}(\Omega;\nu)$ satisfies the following asymptotic estimate:
\begin{gather}\label{eq: convLrhoDolb2}
\left\| u(t)-\overline{u} \right\|_{\varrho;\nu} \leq Q_1 \, t^{-\frac{2(1-\epsilon)}{\varrho(m-1)}} \, \left\| u_0 \right\|_{q_0;\nu^{\phantom{a}}}^{q_0 \frac{\epsilon}{\varrho}} \!\! e^{H_1 \, \left\| u_0 \right\|_{q_0;\nu}^{m-1}}     \\
\ \ \ \forall \epsilon \in (0,1) \, , \ \ \hbox{\rm for a.e.} \,\,  t>1 \, , \nonumber
\end{gather}
where $Q_1=Q_1(\epsilon,\varrho,m,M_P,\nu(\Omega))>0$ and $H_1=H_1(m,M_P,\nu(\Omega)) $.
\begin{proof}
We combine Lemma \ref{lem: dolbPor}, norm interpolation inequalities and the smoothing results proved in Theorem \ref{teo: regPmeNeu}. In fact, given $\varrho >2$, $\epsilon \in (0,1)$ and $t>1$, interpolating between the norms $\| \cdot\|_{2;\nu}$, $\| \cdot \|_{\varrho;\nu}$ and $\| \cdot\|_{(\varrho-2+2\epsilon)/{\epsilon};\nu}$ we obtain:
\begin{equation*}
\left\| u(t)-\overline{u} \right\|_{\varrho;\nu} \leq \left\| u(t)-\overline{u} \right\|_{2;\nu}^{\frac{2}{\varrho}(1-\epsilon)} \left\| u(t)-\overline{u}  \right\|_{\frac{\varrho-2+2\epsilon}{\epsilon};\nu}^{\frac{\varrho-2+2\epsilon}{\varrho}}  \, .
\end{equation*}
Applying \eqref{eq: convL2Dolb} to the first factor on the right hand side we get the time rate. The norm on the second factor can be handled in this way:
\begin{equation*}
\left\| u(t)-\overline{u} \right\|_{\frac{\varrho-2+2\epsilon}{\epsilon};\nu} \leq 2 \left\| u(t) \right\|_{\frac{\varrho-2+2\epsilon}{\epsilon};\nu} \leq 2 \left\| u(1) \right\|_{\frac{\varrho-2+2\epsilon}{\epsilon};\nu} \, .
\end{equation*}
Therefore we arrive at
\begin{equation*}
\left\| u(t)-\overline{u} \right\|_{\varrho;\nu} \leq 2^{\frac{\varrho-2+2\epsilon}{\varrho}} Q^{\frac{2}{\varrho}(1-\epsilon)} \, t^{-\frac{2(1-\epsilon)}{\varrho(m-1)}} \, \left\| u(1) \right\|_{\frac{\varrho-2+2\epsilon}{\epsilon};\nu}^{\frac{\varrho-2+2\epsilon}{\varrho}} \, ;
\end{equation*}
bounding from above the quantity $\left\| u(1) \right\|_{{(\varrho-2+2\epsilon)}/{\epsilon};\nu}$ by the smoothing estimate \eqref{eq: stimaRegPmeNeu1} (evaluated at the time $t=1$) of Theorem \ref{teo: regPmeNeu}, we then obtain \eqref{eq: convLrhoDolb2}.
\end{proof}
\end{corollary}
If in addition the initial datum is essentially bounded, it is easy to check that one can choose $\epsilon=0$ in estimate \eqref{eq: convLrhoDolb2} (indeed it is enough to interpolate between the norms $\| \cdot \|_{2;\nu}$, $\| \cdot \|_{\varrho;\nu}$ and $\| \cdot \|_\infty$). However, in this case we can prove a much stronger result.
\begin{theorem}\label{thm: convExp}
Let $\varrho \in (1,\infty)$, $u_0 \in L^\infty(\Omega)$ and $\overline{u}\neq 0$. If the zero-mean Poincar\'e inequality \eqref{eq: poinIntroMed} holds and $u$ is the energy solution of \eqref{eq: pmeNeu} with initial datum $u_0$, then $u(\cdot,t)$ converges at least exponentially to its mean value. More precisely:
\begin{equation}\label{eq: convExpNeu}
\left\| u(t) - \overline{u} \right\|_{\varrho;\nu} \leq e^{-C |\overline{u}|^{m-1}t} \, \left\| u_0 - \overline{u}  \right\|_{\varrho;\nu}   \ \ \  \hbox{\rm for a.e.}  \,\, t>0 \, ,
\end{equation}
where $C>0$ is a constant depending on $\varrho$, $m$ and $R>0$, the latter being any number such that
\begin{equation}\label{eq: convExpHpinf}
\frac{\left\|u_0 - \overline{u}\right\|_{\infty}}{|\overline{u}|} \leq R \, .
\end{equation}
\begin{proof}
Setting $w=u/\overline{u} - 1$, let us rewrite \eqref{eq: derMediaNeu} as follows:
\begin{equation}\label{eq: derMediaNeuRiscal}
\begin{aligned}
& \frac{\mathrm{d}}{\mathrm{d}s} \left\| w(s) \right\|_{\varrho;\nu}^\varrho  \\
= & -\varrho(\varrho-1)m |\overline{u}|^{m-1}  \int_{\Omega} |w(\mathbf{x},s)+1|^{m-1} |w(\mathbf{x},s)|^{\varrho-2}  |\nabla{w(\mathbf{x},s)} |^2  \, \mathrm{d}\mu \, .
\end{aligned}
\end{equation}
Upon defining
\begin{equation*}
\Phi(x)=\int_0^x |y|^{\frac{\varrho}{2}-1 }|y+1|^{\frac{m-1}{2}} \, \mathrm{d}y \, ,
\end{equation*}
\eqref{eq: derMediaNeuRiscal} becomes
\begin{equation}\label{eq: derMediaNeuRiscal0}
\frac{\mathrm{d}}{\mathrm{d}s} \left\| w(s) \right\|_{\varrho;\nu}^\varrho =  -\varrho(\varrho-1)m |\overline{u}|^{m-1}  \int_{\Omega}  |\nabla{\Phi(w)(\mathbf{x},s)} |^2  \, \mathrm{d}\mu \, .
\end{equation}
At this point need to apply Lemma \ref{lem: GP} to $\Phi$. Such function is certainly continuous and increasing, and by means of de l'H\^{o}pital's Theorem it is straightforward to verify that
\begin{equation*}
\lim_{x \rightarrow 0} \frac{\Phi(x)}{x^{\frac{\varrho}{2}}} = \frac{2}{\varrho} \, , \  \lim_{x \rightarrow \pm \infty } \frac{\Phi(x)}{x^{\frac{\varrho+m-1}{2}}} =\frac{2}{\varrho+m-1} \, .
\end{equation*}
Actually, since trivially ${\varrho}/{2} \neq {(\varrho+m-1)}/{2}$, as observed in Remark \ref{oss: GP} Lemma \ref{lem: GP} is not directly applicable to $\Phi$. However, exploiting the fact that $u_0 \in L^\infty(\Omega)$ and the quantity $\| u(t)-\overline{u} \|_\infty$ is non-expansive (immediate consequence of \eqref{eq: derMediaNeu}), we have that
\begin{equation*}
|w(\mathbf{x},t)| \leq \frac{ \left\| u_0 - \overline{u} \right\|_\infty}{|\overline{u}|} = R
\end{equation*}
for all $t>0$ and a.e.\ $\mathbf{x} \in \Omega$. Therefore the behaviour of $\Phi(x)$ for $|x|$ large has no effect on $\Phi(w)$. In view of that, we are allowed to modify $\Phi$, for instance, as follows:
\begin{equation*}
\Phi_{R}(x)=
\begin{cases}
\Phi(x) & \textnormal{for} \ x \in [-R-2,R]  \\
\Phi(R)+\int_{R}^x  y^{\frac{\varrho}{2}-1 }(R+1)^{\frac{m-1}{2}} \, \mathrm{d}y & \textnormal{for} \ x > R   \\
\Phi(-R-2) - \int_{x}^{-R-2}  |y|^{\frac{\varrho}{2}-1 }(R+1)^{\frac{m-1}{2}} \, \mathrm{d}y & \textnormal{for} \ x < -R-2
\end{cases}
;
\end{equation*}
the function $\Phi_{R}$ satisfies indeed all the hypotheses of Lemma \ref{lem: GP} since
\begin{equation*}
\lim_{x \rightarrow + \infty } \frac{\Phi_{R}(x)}{x^{\frac{\varrho}{2}}} =\lim_{x \rightarrow - \infty } \frac{\Phi_{R}(x)}{x^{\frac{\varrho}{2}}} =\frac{2(R+1)^{\frac{m-1}{2}}} {\varrho} \, .
\end{equation*}
Thus, being $\overline{w}=0$, we know that
\begin{equation}\label{eq: QuasiPhiDisGP}
\left\| \Phi_{R}(w)  \right\|_{2;\nu} \leq C_{\Phi_{R}} \left\| \nabla{\Phi_{R}(w)} \right\|_{2;\mu}
\end{equation}
for a suitable constant constant $C_{\Phi_{R}}>0$. Moreover, as $\Phi(w)=\Phi_{R}(w)$, \eqref{eq: QuasiPhiDisGP} together with \eqref{eq: derMediaNeuRiscal0} give
\begin{equation*} 
\frac{\mathrm{d}}{\mathrm{d}s} \left\| w(s) \right\|_{\varrho;\nu}^\varrho \leq  -\frac{\varrho(\varrho-1)m |\overline{u}|^{m-1}} {C_{\Phi_{R}}^2} \left\| \Phi_{R}(w(s)) \right\|_{2;\nu}^2 \, .
\end{equation*}
In view of the way $\Phi_{R}$ was defined, clearly there exists a constant $D=D(\varrho,m)>0$ such that
\begin{equation*}
D^{-1} |x|^{\frac{\varrho}{2}} \leq |\Phi_{0}(x)| \leq |\Phi_{R}(x)| \ \ \ \forall x \in \mathbb{R}  \, ,
\end{equation*}
so that
\begin{equation*} 
\frac{\mathrm{d}}{\mathrm{d}s} \left\| w(s) \right\|_{\varrho;\nu}^\varrho \leq  -\frac{\varrho(\varrho-1)m |\overline{u}|^{m-1}} {C_{\Phi_{R}}^2  D^2} \left\| w^{\frac{\varrho}{2}}(s) \right\|_{2;\nu}^2  = -{C \varrho}\,|\overline{u}|^{m-1} \left\| w(s) \right\|_{\varrho;\nu}^\varrho    \, .
\end{equation*}
Solving the above differential inequality in the variable $y(s)=\| w(s) \|_{\varrho;\nu}^{\varrho}$ and going back to the original function $u-\overline{u}$ we finally obtain \eqref{eq: convExpNeu}.
\end{proof}
\end{theorem}
Note that, from the proof above, the constant $C$ appearing in \eqref{eq: convExpNeu} depends on $\varrho$, $m$ and $R$ also through the constant $C_{\Phi_R}$ from Lemma \ref{lem: GP}: this, in particular, implies that it is impossible to deduce how $C$ depends on $\| u_0 \|_\infty$, therefore the result just proved in principle is not extendible to data which do not belong to $L^\infty(\Omega)$ (recall that the existence of the constant $C_\Phi$ in Lemma \ref{lem: GP} was established by means of an argument by contradiction).

As concerns \it local convergence \rm this last result can be improved. Notice that \emph{global} uniform convergence need not hold, as we shall see at the end of this section.

\begin{theorem}\label{cor: convLoc}
Let $\nu(\Omega)<\infty$ and suppose that the zero-mean Poincar\'e inequality \eqref{eq: poinIntroMed} holds. Let $u$ be the energy solution of \eqref{eq: pmeNeu} corresponding to an initial datum $u_0\in L^\infty(\Omega)$, with $\overline{u_0}=\overline{u} \neq 0$. Then for any given compact set $\Omega_K \Subset \Omega $ there exist two constants $G$ and $C$ such that
\begin{equation}\label{eq: convLocUnif}
\left\| u(t)-\overline{u} \right\|_{\infty,\Omega_K} \leq G e^{-C |\overline{u}|^{m-1}t} \ \ \ \hbox{\rm for a.e.} \  t>0 \, ,
\end{equation}
where $G$ depends on $m$, $C_P$, $\Omega_K$ and $u$, while $C$ depends on $m$, $N$, $\Omega$, $\rho_\nu$, $\rho_\mu$ but can be taken independent of $u_0$ varying in a set such that $\|u_0-\overline{u} \|_\infty / |\overline{u}| \leq R$, $R$ being any fixed positive constant.
\begin{proof}
In Theorem \ref{thm: convExp} we have shown exponential $L^\varrho$ convergence of $u$ to $\overline{u}$. The rate $C$ in \eqref{eq: convExpNeu} depends only on $\varrho$, $m$, $\Omega$, $\rho_\nu$, $\rho_\mu$ and $R$, with $R$ as in \eqref{eq: convExpHpinf}. From the local regularity results of \cite{DiB} and \cite{PV} we can infer that the solution $u(\cdot,t)$ is spatially H\"{o}lder continuous in any compact set $\Omega_K \Subset \Omega$ and, most importantly, with constants and exponents depending on $\Omega_K$ but not on $t \ge 1 $. In particular, in the case $\rho_\nu\equiv1$ one can apply Theorem 1.2 of \cite{DiB}: in order to do it one has to control, uniformly w.r.t.\ $t>1$, the quantities $\| u(t) \|_{2,\Omega_K}$ and  $\| \nabla(u^m) \|_{2,\Omega_K \times (1,t)} $, which in the present context are straightforward consequences of the energy estimates given in Theorem \ref{teo: teoFondNeu}. If $\rho_\nu\not\equiv1$ one can proceed by setting $w=\rho_\nu u$ and use carefully the results of \cite{DiB} and \cite{PV} under the present regularity assumptions. Arguing by contradiction one obtains that the just mentioned H\"{o}lder continuity and the $L^\varrho$ convergence to the mean value necessary imply \emph{uniform} convergence on $\Omega_K$: this implies in turn that the solution is eventually positive (or negative, depending on the sign of $\overline{u}$), so that \eqref{eq: pmeNeu} for $t$ large enough becomes a non-degenerate quasilinear parabolic equation (locally in space). Therefore, through standard bootstrap techniques, one gets that $| u(t) |_{C^2(\Omega_K)}$ is uniformly bounded in time (again, for $t$ large).

Given a regular function $f$ on $\Omega_K$, define $|f|_{C^0(\Omega_K)}:=\|f\|_\infty$ and, for any multi-index $\eta=(\eta_1,\ldots,\eta_N)$, the quantity $|\eta|=\eta_1+\ldots+\eta_N$ and the seminorms $\left|f\right|_{C^k(\Omega_K)}:=\max_{|\eta|=k}\left\|\partial^\eta f\right\|_\infty$, $k\in{\mathbb N}$. The following generalized interpolation inequalities
\begin{equation}\label{interpolation}
\left|f \right|^{\phantom{\frac 12}}_{C^j(\Omega_K)}\le C_{j,k,p}\, \left| f \right|_{C^k(\Omega_K)}^{\frac{N+jp}{N+kp}}\,\left\|f \right\|_{p^{\phantom{a}}}^{\frac{p(k-j)}{N+kp}}
\end{equation}
are known to hold for all integers $k>j\ge0$ and real $p\ge1$ (see \cite[p. 130]{N} or, for a short review, \cite[App. 3]{BGV}). Estimate \eqref{eq: convLocUnif} then follows by applying \eqref{interpolation} with $f=u(t)-\overline{u}$, $k=2$, $p=2$ (and so \eqref{eq: convExpNeu} for $\varrho=2$) and $j=0$.
\end{proof}
\end{theorem}
As already anticipated, in general one cannot expect \emph{global} uniform convergence to the mean value, so that in some sense the result of Theorem \ref{cor: convLoc} is sharp. We close this section by providing, in the one-dimensional context, classes of weights such that the zero-mean Poincar\'e inequality \eqref{eq: poinIntroMed} holds but convergence to the mean value for the corresponding solutions to \eqref{eq: pmeNeu} does not occur in $L^\infty(\Omega)$.
\\[2mm]
\noindent\bf Counterexamples to the uniform convergence to the mean value. \rm
By means of an explicit counterexample, we have seen in this section that, in general, if the sole zero-mean Poincar\'e inequality \eqref{eq: poinIntroMed} holds then the regularizing effect of equation \eqref{eq: pmeNeu} works up to $L^\varrho(\Omega;\nu)$ with $\varrho$ \emph{strictly} less than infinity. Also, by Corollary \ref{teo: AsiMediaNN} we know that convergence to the mean value for solutions to \eqref{eq: pmeNeu} always takes place in all the $L^\varrho(\Omega;\nu)$ spaces, provided again $\varrho < \infty$, and it is \emph{locally} uniform for bounded data (Theorem \ref{cor: convLoc}). Nonetheless it seems natural to ask whether, at least for regular initial data, \emph{global} uniform convergence of solutions to their mean value holds true. Actually, the answer in general is negative, and we shall now prove this fact through another counterexample.  We stress that the construction of such a counterexample works for all $m>1$.

Indeed, consider equation \eqref{eq: pmeNeu} with the following choices:
\begin{equation*}\label{eq: choices}
\Omega=(0,1) \, , \ \rho_\nu(x)=x^{\beta-2} \, , \ \rho_\mu(x)=x^{\beta} \, , \ \beta>1 \, , \ m \geq 2 \, ;
\end{equation*}
such weights satisfy the zero-mean Poincar\'e inequality \eqref{eq: poinIntroMed} (see Section \ref{list}). We look for a function $r:\mathbb{R}^{+} \rightarrow (0,1)$ regular, decreasing, with $\lim_{t\rightarrow\infty}r(t)=0$ and such that the function
\begin{equation}\label{eq: soprasol}
\hat{u}(x,t)=
\begin{cases}
0 & \textnormal{ for } x \in \left[0, \frac{r(t)}{2} \right] \\
\frac{2x}{r(t)}-1 & \textnormal{ for } x \in \left( \frac{r(t)}{2} , r(t) \right] \\
1 & \textnormal{ for } x \in \left(r(t), 1 \right] \\
\end{cases}
\end{equation}
is a supersolution to \eqref{eq: pmeNeu}. Since $\hat{u}_x(x,t)$ vanishes in neighbourhoods of $x=0$ and $x=1$ for all $t>0$, it is enough to check that $\hat{u}$ is a supersolution in the distributional sense: in other words, it is certainly a supersolution as regards the boundary conditions, so it remains to verify that it is a supersolution also as regards the equation. This amounts to asking that
\begin{gather}\label{eq: supersolEq}
\rho_\nu(x) \hat{u}_t(x,t) \geq \left(\rho_\mu \right)_x(x) \left( \hat{u}^m\right)_x(x,t) + \rho_\mu(x) \left( \hat{u}^m\right)_{xx}(x,t) \\
\quad \hbox{\rm in  } \ \mathcal{D}'((0,1)\times (0,\infty)) \, . \nonumber
\end{gather}
After some straightforward computations, one gets:
$$ \hat u_t(x,t)=-\frac{2r^\prime(t) x}{r^2(t)}\ \chi_{\left( \frac{r(t)}{2} , r(t) \right]}(x) \, , $$
$$ \left( \hat{u}^m\right)_x(x,t)=\frac{2m}{r(t)} \left( \frac{2x}{r(t)}-1 \right)^{m-1}\chi_{\left( \frac{r(t)}{2} , r(t) \right]}(x) \, , $$
\begin{equation}\label{eq: derxx1}
\left( \hat{u}^m\right)_{xx}(x,t) = v(x,t) - \frac{2m}{r(t)} \delta_{x=r(t)}(x,t) \, ,
\end{equation}
where
$$ v(x,t)=\frac{4m(m-1)}{r^2(t)} \left( \frac{2x}{r(t)}-1 \right)^{m-2}\chi_{\left( \frac{r(t)}{2} , r(t) \right]}(x) $$
and, of course, $\left(\rho_\mu \right)_x(x) = \beta x^{\beta-1}$. Since the contribution of the Dirac mass in \eqref{eq: derxx1} is negative, we can neglect it, so that \eqref{eq: supersolEq} holds true if
\begin{gather}
-x^{\beta-2} \frac{2r^\prime(t) x}{r^2(t)} \geq \frac{2m \beta x^{\beta-1}}{r(t)} \left( \frac{2x}{r(t)}-1 \right)^{m-1}
+ \frac{4m(m-1) x^\beta }{r^2(t)} \left( \frac{2x}{r(t)}-1 \right)^{m-2}  \label{eq: supersolEqExpl1}\\
\forall t>0 \, , \ \ \ \forall x \in \left( \frac{r(t)}{2} , r(t) \right) \, . \nonumber
\end{gather}
Dividing \eqref{eq: supersolEqExpl1} by $x^{\beta-1}$ we obtain
\begin{equation}\label{eq: supersolEqExpl2}
-\frac{2r^\prime(t)}{r^2(t)} \geq \frac{2m \beta}{r(t)} \left( \frac{2x}{r(t)}-1 \right)^{m-1}
+ \frac{4m(m-1) x}{r^2(t)} \left( \frac{2x}{r(t)}-1 \right)^{m-2} \, .
\end{equation}
Clearly, for all fixed $t>0$, \eqref{eq: supersolEqExpl2} holds true for all $x \in \left( {r(t)}/{2} , r(t) \right)$ if and only if it holds true at $x=r(t)$ (recall that $m\ge2$). Set then $x=r(t)$ in the right hand side of \eqref{eq: supersolEqExpl2} in order to get
\begin{equation*}\label{eq: supersolEqExpl3}
\frac{r^\prime(t)}{r(t)} \leq -{m (\beta + 2(m-1)) }=-C(m,\beta) <0   \, .
\end{equation*}
Integrating between $0$ and $t$ gives then
\begin{equation}\label{eq: supersolEqExpl5}
r(t) \leq r(0)e^{-C(m,\beta) t } \, .
\end{equation}
Providing that one chooses $r(0)$ small enough, a function $r(t)$ equal to the right-hand side of \eqref{eq: supersolEqExpl5} certainly has all the requirements listed in the beginning of the construction. So we have proved that there exists a supersolution to \eqref{eq: pmeNeu} of the type \eqref{eq: soprasol}. In particular, the \emph{solution} to \eqref{eq: pmeNeu} associated to any positive datum $u_0(x) \leq \hat{u}(x,0)$ will be less or equal than $\hat{u}(x,t)$ for all $t>0$, and as a consequence it will be prevented from converging uniformly to the constant function corresponding to its mean value $\overline{u}>0$ since it is bounded to be zero in $( 0,{r(t)}/{2} )$ for any $t>0$. It is easy to check (by a time scaling argument and by taking $r(0)$ small enough) that the same applies to solutions corresponding to any uniformly bounded initial datum which is less or equal than zero in a neighbourhood of $x=0$ and has positive mean value.

As the reader may note, the assumption $m\geq2$ we made in the beginning cannot be relaxed, since for $m \in (1,2)$ the right hand side of \eqref{eq: supersolEqExpl1} blows up as $x \rightarrow {r(t)}/{2} $. However, in that case we are still able to build a similar supersolution. Indeed, upon setting
\begin{equation*}\label{eq: soprasol2}
\tilde{{u}}(x,t)=
\begin{cases}
\frac{x}{r(t)} & \textnormal{ for } x \in \left[0 , r(t) \right] \\
1 & \textnormal{ for } x \in \left(r(t), 1 \right] \\
\end{cases}
\end{equation*}
and performing analogous computations as above one arrives at
\begin{equation}\label{eq: supersol2EqExpl2}
-\frac{r^\prime(t)}{r^2(t)} \geq \frac{m \beta x^{m-1}}{r^m(t)} + \frac{m(m-1) x^{m-1}}{r^m(t)} \, ,
\end{equation}
which must be valid for $x \in (0,r(t))$. The right hand side of \eqref{eq: supersol2EqExpl2} is clearly maximized at $x=r(t)$, so that by substituting such value in it and solving the resulting differential inequality one obtains again \eqref{eq: supersolEqExpl5} (up to a different positive constant $C(m,\beta)$). Unlike $\hat{u}$, the supersolution $\tilde{{u}}$ has \emph{not} zero derivative at $x=0$: nonetheless, this turns out not to matter. In fact, the space of absolutely continuous functions in $[0,1]$ which vanish in a neighbourhood of $x=0$ is dense in $W^{1,2}((0,1);x^\alpha,x^\beta)$, provided $\alpha \in \mathbb{R}$ and $\beta \geq 1$ (see the proof of \cite[Th. 2.11]{FMT07}). This means that in order to prove that $\tilde{{u}}$ is a supersolution to \eqref{eq: pmeNeu} it is enough to test it on functions which vanish in a neighbourhood of $x=0$, so that its behaviour at $x=0$ is not relevant. The fact that $\tilde{{u}}$ is zero at $x=0$ for all $t>0$ is then sufficient in order to prevent uniform convergence to the mean value $\overline{u} > 0$ for the class of data discussed above.
\begin{remark}
Notice that all the results stated in the sections above for a.e.\ $t>0$ do in fact hold true for all $t>0$ provided weak energy solutions are considered, since the continuity property stated in Remark \ref{remcontinui} holds true also for solutions to the Neumann problem.
\end{remark}
\end{section}

\begin{section}{Some examples of weighted Poincar\'e inequalities}\label{sec: exa}
In the following, we list actual examples of domains $\Omega\subset\mathbb{R}^N$ and weights $\rho_\nu,\rho_\mu$ with respect to which weighted Poincar\'e inequalities hold. The reader should keep in mind that we just aim at stating some significant and explicit results: an exhaustive description of the known theory is far beyond our purposes (as well as hopeless).
\subsection{Poincar\'e inequalities in \texorpdfstring{$\boldsymbol{{W^{1,2}_0(\Omega;\nu,\mu)}}$}{}}\label{sec: exa0}
We begin with the one-dimensional case, where necessary and sufficient conditions are available. Afterwards we shall consider also $N$-dimen\-sional domains.
\\[2mm]\bf{The case $\mathbf{N=1}$.}\label{sec: exa01} \rm
Let $\Omega=(a,b)$, the cases  $a=-\infty$ and $b=+\infty$ being allowed. We look for weights $\rho_\nu,\rho_\mu$ such that the Poincar\'e inequality (sometimes it is also called \emph{Hardy inequality}, having in mind the pioneering weighted inequalities originally proved by G.\ H.\ Hardy)
\begin{equation} \label{eq: poin01d}
\left\| \eta \right\|_{2;\nu} \leq C_P \, \| \eta^\prime \|_{2;\mu^{\phantom{A}}}
\end{equation}
holds for every $\eta$ belonging to a suitable functional space. We shall mainly refer to \cite{KO90}. Accor\-ding to the notation used therein, we indicate as $AC_{\mathcal{L}}(a,b)$ the space of all functions \mbox{$\eta:(a,b)\rightarrow \mathbb{R}$} which are locally absolutely continuous and such that $\lim_{x\rightarrow a^{+}}\eta(x)=0$. The space $AC_{\mathcal{R}}(a,b)$ is understood likewise, replacing $a^{+}$ with $b^{-}$, while $AC_{\mathcal{LR}}=AC_{\mathcal{L}} \cap AC_{\mathcal{R}}$. Since $C^\infty_c(a,b)$ is included in $AC_{\mathcal{LR}}(a,b)$ and it is dense in $W^{1,2}_0((a,b);\nu,\mu)$, the validity of \eqref{eq: poin01d} in $AC_{\mathcal{LR}}(a,b)$ implies in turn the validity of the same inequality in $W^{1,2}_0((a,b);\nu,\mu)$.

\noindent \bf Necessary and sufficient conditions. \rm From \cite[Th. 1.14]{KO90} we have that \eqref{eq: poin01d} holds in $AC_{\mathcal{L}}(a,b)$ if and only if the weights satisfy the following condition:
\begin{equation*} \label{eq: condFond1}
\mathcal{B}_{\mathcal{L}}(a,b,\nu,\mu)=\sup_{x\in(a,b)} \left( \int_x^b \rho_\nu(y)\, \mathrm{d}y \right) \left( \int_a^x \rho_\mu(y)^{-1}\, \mathrm{d}y  \right) < \infty \, .
\end{equation*}
Similarly, \eqref{eq: poin01d} holds in $AC_{\mathcal{R}}(a,b)$ if and only if
\begin{equation*} \label{eq: condFond2}
\mathcal{B}_{\mathcal{R}}(a,b,\nu,\mu)=\sup_{x\in(a,b)}  \left(\int_a^x \rho_\nu(y)\, \mathrm{d}y \right) \left(\int_x^b \rho_\mu(y)^{-1}\, \mathrm{d}y\right)  < \infty \, .
\end{equation*}
It is then possible to show \cite[Th. 8.8]{KO90} that the existence of a constant $c \in [a,b] $ such that, setting conventionally $\mathcal{B}_{\mathcal{L}}(a,a,\cdot,\cdot)=\mathcal{B}_{\mathcal{R}}(b,b,\cdot,\cdot)=0$,
\begin{equation} \label{eq: condFond3}
\mathcal{B}_{\mathcal{L}}(a,c,\nu,\mu)< \infty \, , \ \mathcal{B}_{\mathcal{R}}(c,b,\nu,\mu) < \infty
\end{equation}
is necessary and sufficient for the validity of \eqref{eq: poin01d} in $AC_{\mathcal{LR}}(a,b)$. Actually, the same result holds true replacing $AC_{\mathcal{LR}}(a,b)$ with the space of absolutely continuous functions in $(a,b)$ with compact support. For a more general discussion about inequality \eqref{eq: poin01d}, possibly involving weights which are not absolutely continuous with respect to the Lebesgue measure, see \cite[Sec. 1.3]{Maz85} or also the pioneering work of B.\ Muckenhoupt \cite{Muc72a, Muc72b}.


We refer the reader to Section \ref{list} for a collection of explicit weights for which, using the above results, a Poincar\'e inequality can be proved, whereas no Sobolev inequality holds. The latter fact can be shown through the lack of validity of the corresponding necessary and sufficient conditions given in \cite[Th. 8.8]{KO90} for Sobolev-type inequalities. See also \cite[Sec. 6]{KO90}.
\\[2mm] \bf The case $\mathbf{N\geq 1}$. \rm It is well-known, see e.g.\ \cite[Sec. 15]{KO90}, that the so called \emph{Muckenhoupt classes} of weights originally introduced in \cite{Muc72b} have an important role in weighted functional inequalities. In fact, Muckenhoupt weights are defined by a sort of generalization of \eqref{eq: condFond3} for $N$-dimensional domains. We shall not investigate further such theory: instead, we prefer to give some explicit examples, for which we refer again to \cite{KO90} (see also Section \ref{list}). 

\noindent \bf Bounded domains. \rm Let $\Omega\subset \mathbb{R}^N$ be a bounded Lipschitz domain. We indicate as $\delta(\mathbf{x})={\rm dist}\,(\mathbf{x},\partial\Omega)$ the distance function of $\Omega$. Consider a parameter $\beta < 1$. From \cite[Th. 21.5]{KO90} we have that the Poincar\'e inequality
\begin{equation}\label{eq: poinEse}
\left\| v \right\|_{2;\nu} \leq C_P \left\| \nabla{v} \right\|_{2;\mu}
\end{equation}
holds for all $v \in W^{1,2}_0(\Omega;\delta^\alpha,\delta^\beta)$ if and only if $\alpha\geq\beta-2$. Moreover, Sobolev-type inequalities hold if in addition $\alpha > \beta -2$. When $\beta\geq 1$ it is possible to prove that $W^{1,2}_0(\Omega;\delta^\beta,\delta^\beta)=W^{1,2}(\Omega;\delta^\beta,\delta^\beta)$ \cite[Th. 2.11]{FMT07}: this is enough to conclude that for any $\alpha \in \mathbb{R}$ and $\beta\geq1$ the Poincar\'e inequality cannot be valid in $W^{1,2}_0(\Omega;\delta^\alpha,\delta^\beta)$.

\noindent \bf Exterior domains. \rm If $\Omega\subset \mathbb{R}^N$ is an exterior domain (namely the complement of any compact set) such that $\inf_{\mathbf{x}\in \Omega}|\mathbf{x}|>0$, then the Poincar\'e inequality \eqref{eq: poinEse} holds in $W^{1,2}_0(\Omega;|\mathbf{x}|^\alpha,|\mathbf{x}|^\beta)$ if and only if $\beta \neq 2-N$ and $\alpha \leq \beta -2$ \cite[Ex. 21.10]{KO90}. Under these conditions, Sobolev-type inequalities are valid either if $\beta>2-N$ or if $\beta<2-N$ and $\alpha<\beta-2$. The same results apply with respect to the weights $(|\mathbf{x}|+1)^\alpha , (|\mathbf{x}|+1)^\beta$ replacing $\Omega$ with $\mathbb{R}^N$, provided $\beta>2-N$; if instead $\beta<2-N$ then the Poincar\'e inequality \eqref{eq: poinEse} in this case does not hold (indeed one can prove that constants belong to $W^{1,2}_0(\mathbb{R}^N;(|\mathbf{x}|+1)^\alpha,(|\mathbf{x}|+1)^\beta)$ if $\alpha\leq\beta-2<-N$). For the limiting case $\beta=2-N$ and the peculiar functional inequalities satisfied for a suitable value of $\alpha$ we refer the reader to \cite{BD+10} and \cite{BGV}.

Similar results hold for exponential weights of the form $e^{\alpha|\mathbf{x}|}$ (see \cite[Ex. 21.12]{KO90} or Section \ref{list} for an outline).
\\[2mm]
\noindent{\bf An example in Riemannian geometry.} For the sake of simplicity we consider the following example only in the one-dimensional setting. In fact, weighted one-dimensional inequalities on $(0,+\infty)$
admit in some cases a geometric interpretation which we now sketch; such construction can be performed e.g.\ when $\rho_\nu=\rho_\mu$.  We consider a $C^2$ Riemannian manifold $M$, of dimension $N$, with a pole $o$ given on it and whose metric is defined, in polar or spherical coordinates around $o$, as
\begin{equation*}\label{metric}
{\rm d}s^2={\rm d}r^2+\psi(r)^2{\rm d}\Theta^2.
\end{equation*}
Here ${\rm d}\Theta^2$ denotes the canonical metric on the Euclidean unit sphere ${\mathbb S}^{N-1}$, the function $\psi$ is assumed to be smooth and positive on $(0,+\infty)$, with $\psi(0)=\psi^{\prime\prime}(0)=0$, $\psi^\prime(0)=1$ (the prime indicating right derivative), and $r$ is by construction the Riemannian distance between a point whose coordinates are $(r,\Theta)$ and $o$. The conditions on $\psi$ ensure that $M$ is $C^2$ in a neighborhood of $o$. A manifold satisfying the above conditions is said to be a \it Riemannian model \rm with pole $o$. The running assumptions entail that $M$ is complete.

The Riemannian Laplacian of a scalar function $f$ on $M$ is given, in the above coordinates, by
\[
\begin{aligned}
\Delta f(r,\theta_1,\ldots,\theta_{N-1})= & \frac1{\psi(r)^{N-1}}\frac{\partial}{\partial r}\left[\psi(r)^{N-1} \frac{\partial f}{\partial r}(r,\theta_1,\ldots,\theta_{N-1})
\right] \\
 & + \frac1{\psi(r)^2}\Delta_{{\mathbb S}^{N-1}}f(r,\theta_1,\ldots,\theta_{N-1}) \, ,
\end{aligned}
\]
where $\Delta_{{\mathbb S}^{N-1}}$ is the Riemannian Laplacian on the unit sphere ${\mathbb S}^{N-1}$. In particular, for \it radial \rm functions, namely functions depending only on the geodesic distance $r$, one has
\[
\Delta f(r)=\frac1{\psi(r)^{N-1}}\left[\left(\psi(r)^{N-1}\right)f^\prime(r)\right]^\prime \, ,
\]
where the prime denotes derivative w.r.t.\ $r$. Consider the inequality
\begin{gather}\label{gap}
\int_0^\infty f(r)^2\psi(r)^{N-1}\,{\rm d}r\le C\int_0^\infty f^\prime(r)^2\psi(r)^{N-1}\,{\rm d}r   \\
\forall  f : {\rm \ supp}\,f\subset[0,b] \, , \ f^\prime(0)=0   \nonumber
\end{gather}
for a suitable $b=b(f)>0$. If \eqref{gap} holds, then it is easy to realize that it can be extended to \emph{all} radial functions in $W^{1,2}(M)$. By construction, the volume element on $M$ is $\psi(r)^{N-1}\,{\rm d}r\,{\rm d}\omega_{N-1}$, where ${\rm d}\omega_{N-1}$ is the volume element on the Euclidean unit sphere ${\mathbb S}^{N-1}$. Inequality \eqref{gap} (or its extension to all radial function in $W^{1,2}(M)$) is therefore equivalent to the fact that $\mathcal{S}(-\Delta_r)\subset \left[1/C,+\infty\right)$, where $\Delta_r$ is the Laplacian for radial functions and $\mathcal{S}(\mathcal{L})$ denotes the $L^2$ spectrum of an operator $\mathcal{L}$. Necessary and sufficient conditions for the existence of a spectral gap for $-\Delta$ are known, and of course their validity implies the existence of a spectral gap for $-\Delta_r$ as well. From \cite[Th. 2.10]{Gry06} one then easily argues that the spectral gap for $-\Delta$ holds if and only if there exists $Q>0$ (independent of $r,\xi$) such that
\begin{equation}\label{gap2}
\left(\int_0^\xi\psi(s)^{N-1}\,{\rm d}s\right)\left(\int_\xi^{r}\frac1{\psi(s)^{N-1}}\,{\rm d}s\right)\le Q \ \ \ \forall r>0 \, , \ \forall \xi\in(0,r) \, .
\end{equation}
Therefore under condition \eqref{gap2} the $L^2$ spectrum of $-\Delta_r$ is bounded away from zero. For example, a smooth function $\psi$ satisfying $\psi(r)=r$ for $r\in(0,1)$ and $\psi(r)= e^{ar}$ ($a\not=0$) for $r>2$ does fulfil the above condition, and it can be shown that in this case no radial Sobolev inequality holds true.

As a consequence of the above discussion, we stress that whenever $\psi(s)$ satisfies \eqref{gap2} then \emph{radial} solutions to the Porous Media Equation on the Riemannian model $M$ associated to such $\psi$ enjoy the $L^{q_0}$-$L^{\varrho}$ regularizing property discussed in Section \ref{sec: regAs}. We omit the details, and refer to \cite{Heb} for a discussion of the technical conditions concerning the validity, or lack of validity, of Sobolev inequalities on $M$.

\subsection{Zero-mean Poincar\'e inequalities in \texorpdfstring{${W^{1,2}(\Omega;\nu,\mu)}$}{}}\label{sec: exaM}
In the framework of zero-mean Poincar\'e inequalities, less results are available with respect to those known for the Poincar\'e inequalities. See for example \cite{BBCG} and references quoted for a clever approach generalizing the Bakry-Emery criterion. We shall anyway confine ourselves to a list of significant examples.
\\[2mm]
\bf The case $\mathbf{N=1}$. \rm We set $\Omega=(a,b)$, the cases  $a=-\infty$ and $b=+\infty$ being allowed. Given two weights $\rho_\nu,\rho_\mu$ defined on $(a,b)$, with $\nu(a,b)<\infty$, consider the quantities
\begin{equation*}
\mathcal{K}_{\mathcal{L}}(a,b,\nu,\mu)=\sup_{x\in(a,b)}  \left(\int_x^b \rho_\nu(y)\, \mathrm{d}y \right) \left(\int_a^x \left( \int_a^y \rho_\nu(t) \, \mathrm{d}t \right)^2 \rho_\mu(y)^{-1}\, \mathrm{d}y\right) \, ,
\end{equation*}
\begin{equation*}
\mathcal{K}_{\mathcal{R}}(a,b,\nu,\mu)=\sup_{x\in(a,b)}  \left(\int_a^x \rho_\nu(y)\, \mathrm{d}y \right) \left(\int_x^b \left( \int_y^b \rho_\nu(t) \, \mathrm{d}t \right)^2 \rho_\mu(y)^{-1}\, \mathrm{d}y\right)  \, .
\end{equation*}
From \cite[Th. 1.4]{CW99} we have that the zero-mean Poincar\'e inequality
\begin{equation*} \label{eq: poinM1d}
\left\| \eta -\overline{\eta} \right\|_{2;\nu} \leq M_P \, \| \eta^\prime \|_{2;\mu^{\phantom{A}}}
\end{equation*}
holds for all $\eta \in W^{1,2}((a,b);\nu,\mu)$ if and only if
\begin{equation*}\label{eq: condNSMed}
\mathcal{K}_{\mathcal{L}}(a,b,\nu,\mu)+\mathcal{K}_{\mathcal{R}}(a,b,\nu,\mu) < \infty \, .
\end{equation*}

For explicit examples of weights depending on elementary functions that satisfy condition \eqref{eq: condNSMed} see again Section \ref{list}. To infer other properties such as the lack of validity of Sobolev inequalities we still refer to the results provided by the general Theorem 1.4 of \cite{CW99}.
\\[2mm]
\bf The case $\mathbf{N\geq 1}$. \rm Now we present some specific examples of weighted zero-mean Poincar\'e inequalities in the $N$-dimensional context, both for bounded domains and the Euclidean space. \\
\noindent \bf Bounded domains. \rm If $\Omega\subset \mathbb{R}^N$ is a bounded star-shaped domain, $w:(0,+\infty)\rightarrow(0,+\infty)$ is any increasing function such that $w(sr)\geq sw(r) \ \forall s \in (0,1)$ and $k$ is any integer, then by \cite[Th. 1]{BK98} there exists a constant $M_P$ such that the following inequality holds:
\begin{equation}\label{eq: buck1}
 \int_{\Omega} 	\left|\eta(\mathbf{x}) -\overline{\eta} \right|^2  \, w(\delta(\mathbf{x}))^k \mathrm{d}\mathbf{x}  \leq M_P^2 \int_{\Omega} \left|\nabla{\eta}(\mathbf{x})\right|^2 \, w(\delta(\mathbf{x}))^k \mathrm{d}\mathbf{x}  \ \ \ \forall \eta \in C^1(\Omega) \, .
\end{equation}
In this case $\nu$ and $\mu$ satisfy the hypotheses of Proposition \ref{pro: dens1}, so that by density \eqref{eq: buck1} can be extended to the whole $W^{1,2}(\Omega;(w\circ\delta)^k,(w\circ\delta)^k )$.



If $\Omega \subset \mathbb{R}^N$ is a bounded convex domain, from \cite[Th. 1.1]{CW10} we have that the inequality
\begin{equation}\label{eq: eseConvx}
\left\| v - \overline{v}\right\|_{2;\delta^{\beta-2}} \leq M_P \left\| \nabla{v} \right\|_{2;\delta^\beta}
\end{equation}
holds in $W^{1,\infty}_{loc}(\Omega) \cap W^{1,2}(\Omega;\delta^{\beta-2},\delta^\beta) $, that is in the whole $W^{1,2}(\Omega;\delta^{\beta-2},\delta^\beta) $ (again by Proposition \ref{pro: dens1}), provided $\beta \geq 2$. The boundedness of the domain implies in turn that \eqref{eq: eseConvx} continues to hold if one replaces $\delta^{\beta-2}$ by $\delta^{\gamma}$, for any $\gamma \geq \beta-2$ (Proposition \ref{pro: OrdPoin}). Sobolev-type inequalities hold providing that  $\gamma > \beta - 2$, while at the limit value $\gamma=\beta-2$ there is no Sobolev embedding (see \cite[Th. 1.1]{CW10} and \cite[Ex. 18.15, Th. 19.9, Th. 19.11]{KO90}). For analogous theory in less regular domains we refer to \cite{Hur90}.\\
\noindent \bf The Euclidean space. \rm We consider power-type weights defined on $\mathbb{R}^N$, where for simplicity we assume $N\ge3$ (see \cite{BD+10} for the cases $N=1,2$):
\begin{equation*}
\rho_\nu(\mathbf{x})=(1+|\mathbf{x}|^2)^{\alpha-1} \, , \ \rho_\mu(\mathbf{x})=(1+|\mathbf{x}|^2)^{\alpha} \, .
\end{equation*}
The papers \cite{BB+07}, \cite{BD+10} show that the zero-mean Poincar\'e inequality $\left\| v -\overline{v} \right\|_{2;\nu} \leq M_P \left\| \nabla{v} \right\|_{2;\mu}$ is valid in the appropriate Sobolev space if and only if $\alpha<1-N/2$. From \cite[Ex. 20.6]{KO90} it follows that for such weights no Sobolev-type embedding holds.

As concerns the Gaussian weights $\rho_\nu(\mathbf{x})=\rho_\mu(\mathbf{x})=e^{-d |\mathbf{x}|^2}$ ($d>0$) we refer the reader to \cite{F}, \cite{Gro75}, \cite[Ch. 2, Sec. 4.3]{Dav89} and \cite{DF92}.

\begin{remark}\rm
The distance function $\delta$ of a domain $\Omega \subset \mathbb{R}^N$, which we often used above, in general need not be more regular than Lipschitz. However, thanks to a well-known theorem due to Whitney \cite[Th. VI.2]{Ste70}, it is always possible to construct a function $\widetilde{\delta}$ which is $C^\infty(\Omega)$ and equivalent to $\delta$, so that in many of the the examples we have seen the weights considered can be \emph{locally} regularized without affecting the validity of the mentioned results just by replacing $\delta$ with $\widetilde{\delta}$.
\end{remark}

\subsection{Explicit weights}\label{list}
For the reader's convenience, in the sequel we list concisely known cases of couples of weights for which Poincar\'e inequalities hold but Sobolev-type inequalities do not, so that the results we provided for the WPME are in fact new.
\vskip 2mm
\noindent\emph{For the following weights, Poincar\'e inequalities hold in ${W^{1,2}_0}$ but no Sobolev-type inequality holds in the same space:}
\begin{itemize}
\item[$\bullet$] \emph{Intervals}:
\begin{itemize}
\vskip 1mm
\item[$\circ$] $(x^{\beta-2},x^\beta)$ for $\beta \neq 1$ on $(0,+\infty)$, $(0,b)$ or $(a,+\infty)$ [let $a,b>0$];
\vskip 1mm
\item[$\circ$] $\left(\frac1{x}|\log{x}|^{\beta-2}, {x}|\log{x}|^{\beta}\right)$ for $\beta \neq 1$ on $(0,1)$;
\vskip 1mm
\item[$\circ$] $(e^{\alpha x},e^{\alpha x})$ for $\alpha \neq 0$ on $\mathbb{R}$;
\end{itemize}
\vskip 1mm
\item[$\bullet$] \emph{Bounded Lipschitz domains}:
\begin{itemize}
\vskip 1mm
\item[$\circ$] $(\delta^{\beta-2},\delta^\beta)$ for $\beta < 1$;
\end{itemize}
\vskip 1mm
\item[$\bullet$] \emph{Exterior domains}:
\begin{itemize}
\vskip 1mm
\item[$\circ$] $(|\mathbf{x}|^{\beta-2},|\mathbf{x}|^\beta)$ for $\beta<2-N$;
\vskip 1mm
\item[$\circ$] $(e^{\alpha|\mathbf{x}|},e^{\alpha|\mathbf{x}|})$ for $\alpha < 0$.
\end{itemize}
\end{itemize}
\vskip 2mm
\noindent\emph{For the following weights, zero-mean Poincar\'e inequalities hold in ${W^{1,2}}$ but no Sobolev-type inequality holds in the same space:}
\begin{itemize}
\item[$\bullet$] \emph{Intervals}:
\begin{itemize}
\vskip 1mm
\item[$\circ$] $(x^{\beta-2},x^\beta)$ for $\beta > 1$ on $(0,b)$ or for $\beta<1$ on $(a,+\infty)$ [let $a,b>0$];
\vskip 1mm
\item[$\circ$] $\left(\frac1{x}|\log{x}|^{\beta-2}, {x}|\log{x}|^{\beta}\right)$ for $\beta \neq 1$ on $(0,c)$, with $c\in (0,1)$;
\vskip 1mm
\item[$\circ$] $(e^{\alpha |x|},e^{\alpha |x|})$ for $\alpha < 0$ on $\mathbb{R}$;
\end{itemize}
\vskip 1mm
\item[$\bullet$] \emph{Bounded convex domains}:
\begin{itemize}
\vskip 1mm
\item[$\circ$] $(\delta^{\beta-2},\delta^\beta)$ for $\beta \geq 2$;
\end{itemize}
\vskip 1mm
\item[$\bullet$] \emph{The Euclidean space $\mathbb{R}^N$}:
\begin{itemize}
\vskip 1mm
\item[$\circ$] $((1+|\mathbf{x}|^2)^{\alpha-1},(1+|\mathbf{x}|^2)^{\alpha})$ for $\alpha<1-\frac{N}{2}$;
\vskip 1mm
\item[$\circ$] $(e^{-d|\mathbf{x}|^2},e^{-d|\mathbf{x}|^2})$ for $d > 0$.
\end{itemize}
\end{itemize}

\end{section}

\section*{Acknowledgments} We thank J.\ L.\ V\'azquez for helpful discussions. In particular, part of this work was done when M.\ M.\ was visiting him at the Department of Mathematics of the Universidad Aut\'onoma de Madrid.


\end{document}